\documentclass[UTF-8,reqno]{amsart}
\usepackage{enumerate}
\usepackage{mhequ}
\usepackage{amssymb,url,color, booktabs,nccmath}
\usepackage[left=3.2cm,right=3.2cm,top=4cm,bottom=4cm]{geometry}
\usepackage{mathrsfs}
\usepackage{enumitem,dsfont}

\usepackage{color}
\usepackage[colorlinks=true]{hyperref}
\hypersetup{
    linkcolor=blue,          
    citecolor=red,        
    filecolor=blue,      
    urlcolor=cyan
}

\definecolor{darkergreen}{rgb}{0.0, 0.5, 0.0}

\numberwithin{equation}{section}
\def\theequation{\arabic{section}.\arabic{equation}}
\newcommand{\be}{\begin{eqnarray}}
\newcommand{\ee}{\end{eqnarray}}
\newcommand{\ce}{\begin{eqnarray*}}
\newcommand{\de}{\end{eqnarray*}}
\newtheorem{theorem}{Theorem}[section]
\newtheorem{lemma}[theorem]{Lemma}
\newtheorem{proposition}[theorem]{Proposition}
\newtheorem{Examples}[theorem]{Example}
\newtheorem{corollary}[theorem]{Corollary}
\newtheorem{assumption}{Assumption}[section]
\newtheorem{definition}[theorem]{Definition}
\theoremstyle{definition}
\newtheorem{remark}[theorem]{Remark}

\def\${|\!|\!|}
\def\lv{\left\Vert}
\def\rv{\right\Vert}

\def\eps{\varepsilon}

\def\a{\alpha}

\def\d{\delta}
\def\de{D((-A)^{\delta})}
\def\C{C_{t,x}}
\def\g{\gamma}

\def\<{{\langle}}
\def\>{{\rangle}}
\def\({{\Big(}}
\def\){{\Big)}}

\def\bx{{\mathbf{x}}}

\def\dif{{\mathord{{\rm d}}}}

\def\={&\!\!=\!\!&}

\def\k{{\kappa}}

\def\cB{{\mathcal B}}

\def\cR{{\mathcal R}}
\def\cS{{\mathcal S}}
\def\cT{{\mathcal T}}

\def\mN{{\mathbb N}}

\def\mP{{\mathbb P}}

\def\mR{{\mathbb R}}

\def\mT{{\mathbb T}}

\def\mZ{{\mathbb Z}}

\def\1{{\mathbf{1}}}

\def\E{\mathbf E}

\def\geq{\geqslant}
\def\leq{\leqslant}
\def\ge{\geqslant}

\def\div{\mathord{{\rm div}}}

\def\eps{\varepsilon}

\def\a{\alpha}

\def\g{\gamma}

\def\<{{\langle}}
\def\>{{\rangle}}
\def\({{\Big(}}
\def\){{\Big)}}

\def\z{\mathring{z}}

\def\bx{{\mathbf{x}}}

\def\dif{{\mathord{{\rm d}}}}

\def\={&\!\!=\!\!&}
\def\bt{\begin{theorem}}
\def\et{\end{theorem}}
\def\bl{\begin{lemma}}
\def\el{\end{lemma}}
\def\br{\begin{remark}}
\def\er{\end{remark}}
\def\bx{\begin{Examples}}
\def\ex{\end{Examples}}
\def\bd{\begin{definition}}
\def\ed{\end{definition}}
\def\bp{\begin{proposition}}
\def\ep{\end{proposition}}
\def\bc{\begin{corollary}}
\def\ec{\end{corollary}}

\def\geq{\geqslant}
\def\leq{\leqslant}
\def\ge{\geqslant}

\def\div{\mathord{{\rm div}}}

\def\Id{\textrm{Id}}

 \def\R{\mathbb R}
 \def\R{\mathbb R}    
\def\N{\mathbb N}  
   
\def\<{\langle} \def\>{\rangle}

\allowdisplaybreaks

\begin{document}

\title[Non-uniqueness in law of stochastic 3D NSE with general multiplicative noise]{Non-uniqueness in law of stochastic 3D Navier--Stokes equations with general multiplicative noise}

\author{Huaxiang L\"u}
\address[H. L\"u]{Academy of Mathematics and Systems Science,
Chinese Academy of Sciences, Beijing 100190, China}
\email{lvhuaxiang22@mails.ucas.ac.cn}

\author{Yichun Zhu}
\address[Y.Zhu]{Academy of Mathematics and Systems Science,
Chinese Academy of Sciences, Beijing 100190, China}
\email{steven00931002@hotmail.com}

\thanks{
Y.Zhu is grateful for the funding from China Postdoctoral Science Foundation (Certificate number: 2024T170981)
}

\begin{abstract}
We are concerned with the three dimensional Navier-Stokes equations driven by a general multiplicative noise. For every divergence free and mean free initial condition in $L^2$, we establish existence of infinitely many global-in-time probabilistically strong and analytically weak solutions, which implies non-uniqueness in law. Moreover, we prove the existence of infinitely many ergodic stationary solutions. Our results are based on a stochastic version of the convex integration and the Itô calculus.
 \end{abstract}
 
\subjclass[2010]{60H15; 35R60; 35Q30}
\keywords{stochastic Navier--Stokes equations, probabilistically strong solutions, non-uniqueness in law, (ergodic) stationary solutions, convex integration}

\date{\today}

\maketitle

\tableofcontents


\section{Introduction}
In this paper, we  study the following stochastic Navier-Stokes equation on $\mathbb{T}^3$ with a general multiplicative noise:
\begin{equation}
\label{1}
\aligned
 \dif u+\div(u\otimes u)\,\dif t+\nabla P\,\dif t&= \Delta u \,\dif t+G(u)\dif W,
\\
\div\ u&=0,
\endaligned
\end{equation}
where $u$ is the velocity of fluid, $P$ is the associated pressure and $G(u)\dif W$ represents a stochastic forcing acting on the fluid. 

The research focused on the well-posedness of the three-dimensional incompressible Naiver-Stokes equations has experienced an immense breakthrough recently. Buckmaster and Vicol \cite{BV19b} established non-uniquenss of weak solutions with finite kinetic energy.  Buckmaster, Colombo and Vicol \cite{BCV18} were able to connect two arbitrary strong solutions via a weak solution.  A number of  other ill-posed results, including non-uniqueness of weak solutions to the Naiver-Stokes equations for every given divergence free initial condition \cite{BMS21} and sharp non-uniqueness for the Navier-Stokes equation in dimension $d \geq 2$ \cite{CL20,CL23} were obtained. All these groundbreaking results are made through the method of convex integration introduced by De Lellis and L. Szekelyhidi \cite{DS09,DS10,DS13}.  We also mention that convex integration  leads to  the proof of Onsager's conjecture for the incompressible Euler equations \cite{Ise18,BDLSV19}.  We refer to the reviews \cite{BV19, BV21} for further details on this method.

It is  well-known that for a given dynamic system, random perturbation has a regularization effect on the system. Classical results  can be found in \cite{DZ14,DFPR13,DF10} as well as the references there. 
The regularization effect of  noise has aroused much attention in the community. For the study of transport noise, we refer readers to the results by Flandoli, Gubinelli, Luo and Priola \cite{FGP10,FL21} and the reference there. So, it's natural to ask 
whether the stochastic Navier-Stokes equation \eqref{1} is well-posed or not, and whether the non-uniqueness in law holds for the stochastic Navier-Stokes equation \eqref{1}. 

A lot of results on this subject has  been made in recent years using the convex integration method.  The first work on convex integration  in the stochastic setting is due to Breit, Feireisl and Hofmanov\'a \cite{BFH20}.   M.Hofmanov\'a, R.Zhu and X.Zhu first proved the non-uniqueness in law and the strong non-uniqueness of probabilistically strong solutions to the Naiver-Stokes system when the random perturbation is either additive or multiplicative  \cite{HZZ23a,HZZ24}. The space-time white noise case was studied in \cite{HZZ21b,LZ25}. Later, W.Chen, Z.Dong and X.Zhu also proved the sharp non-uniqueness of stochastic Navier-Stokes equation with additive noise \cite{CDZ22}. Recently, A.Cheskidov, Z.Zeng and D.Zhang constructed solutions to stochastic Navier-Stokes equation with additive noise with both prescribed energy and initial condition \cite{CZZ24}.  All these results only work for additive noise or multiplicative noise with very regular coefficients for the covariance operator, $G(u)$ in \eqref{1}, of the random perturbation.
 Recently, M.Hofmanová, T. Lange and U.Pappalettera proved the global existence and non-uniqueness of probabilistically strong solutions for 3D Euler/Navier-Stokes equations with transport noise \cite{HLP23,Pap24}.  Their method, based on a geometric flow, is innovative but only works typically for the transport noise and fails to cover other types of multiplicative noise.

One of the main purposes of this paper is to study whether the  non-uniqueness in law holds for a more general family of the random perturbation. Previously, the method adapted by M.Hofmanov\'a, R.Zhu and X.Zhu in \cite{HZZ24} is a combination of the convex integration scheme and the pathwise control from  rough path theory which requires high regularity for the noise part. In this paper, a general framework is provided which enables us to prove the non-uniqueness of probabilistically strong solution of \eqref{1} with more singular coefficients, by developing the stochastic convex integration technique first introduced in \cite{CDZ22,HZZ22}. 

Another important topic in fluid mechanics is the existence and uniqueness of stationary solution. For additive trace class noise and 2D space-time white noise, the existence and non-uniqueness of stationary solutions have been established in \cite{HZZ22,LZ24}. However, in the general case, the existence and uniqueness of  stationary solutions remain unresolved. In particular, for the multiplicative noise $G(u)=u$, it seems that only the trivial stationary solution $0$ is known  (see Remark \ref{r:lm} below). As an extension of the work \cite{HZZ22}, we also proved that there are infinitely many ergodic stationary solutions for general multiplicative noise. 

\subsection{Main results}\label{sec:main:result}

Let $(\Omega, \mathbf{P}, \mathcal{F})$ be a probability space and $W_t$, $t\in [0, \infty)$ be a cylindrical Wiener process in a Hilbert space $U$ defined on this probability space (For a precise definition of cylindrical Wiener process, see \cite[Section 4.1.2]{DZ14}). First of all, we intend to prescribe an arbitrary  divergence-free, mean-zero random initial condition $u_0 \in L^2$, $\mathbf{P}$-a.s., independent of the given Wiener process $W$. Let $(\mathcal{F}_t)_{t\geq 0}$ be the augmented joint canonical filtration on $(\Omega,\mathcal{F})$ generated by $W$ and $u_0$. 

In this paper, we make the following assumption on $G(u)$, where we refer readers to Section \ref{sec:notations} for the precise definition of $\mathcal{L}^2_0(U_1;U_2)$ and $H^{-2\d_0}$.
\begin{assumption}\label{a:G:1}
Let $p_0 \in [1,2)$ and $\d_0 \in [0,1/2)$, we assume $G$ is a linear operator from $L^{p_0}$ to $\mathcal{L}_0^2(U;H^{-2\d_0})$. 
Moreover, there exists a constant $L>0$ such that 
\begin{equation}\label{eq:G:lip}
\|G(u_1)-G(u_2)\|_{\mathcal{L}^2_0(U; H^{-2\d_0})} \leq L \|u_1-u_2\|_{L^{p_0}},
\end{equation}
and
\begin{equation}\label{eq:G:lg}
\|G(u)\|_{\mathcal{L}^2_0(U;H^{-2\d_0})}\leq L (1+ \|u\|_{L^{p_0}}).
\end{equation}
\end{assumption}

First, we consider the analytically weak and probabilistically strong solution to the stochastic Naiver--Stokes equation \eqref{1} in the following sense:
\begin{definition}\label{d:sol}
We say that $u$ is an analytically weak and probabilistically strong solution to the Navier--Stokes system \eqref{1} provided
\begin{enumerate}
\item the velocity $u \in L^2_{\rm{loc}}([0,\infty);L^2_{\sigma})\cap C([0,\infty);H^{-\delta})$ $\mathbf{P}$-a.s. for some $\delta>0$ and  is $(\mathcal{F}_{t})_{t \geq 0}$-adapted (see Section \ref{sec:notations} for the precise definition of $L^2_{\sigma}$);
\item for every $0\leq s\leq t<\infty$ it holds  $\mathbf{P}$-a.s.
$$
\begin{aligned}
&\langle u(t),\psi \rangle + \int_{s}^{t}\langle \div(u\otimes u),\psi \rangle \dif r\ =\langle u(s),\psi \rangle +\int_{s}^{t} \langle \Delta u, \psi\rangle \,\dif r +\langle \int_s^t G(u) \dif W_r,\psi\rangle
\end{aligned}
$$
for all $\psi\in C^{\infty}(\mathbb{T}^{3}),$ $\div\psi=0$.
\end{enumerate}
\end{definition}

The first main result of the present paper is the existence of infinitely many probabilistically strong and analytically weak solutions  for every given divergence-free initial condition in $L^2$.

\begin{theorem}\label{thm:main:3}
Suppose Assumption \ref{a:G:1} holds. Let $u_0 \in L^2_{\sigma}$, $\mathbf{P}-a.s.$, be independent of the cylindrical Wiener process $W_t$, $t \geq 0$. There exist infinitely many analytically weak and probabilistically strong solutions to \eqref{1} in the sense of Definition~\ref{d:sol}. The  solutions belong to  $C([0,\infty);L^{p_0}) \cap C((0,\infty);W^{\vartheta,p_0}) \cap L^2_{\rm{loc}}([0,\infty);L^2)$, $\mathbf{P}$-a.s.  for some $\vartheta>0$.
 \end{theorem}
 Prior to this work, the global existence of non-unique probabilistically strong solutions was known only for certain specific types of noise, as demonstrated in \cite{HZZ24,HZZ23a,Pap24}. Our main result establishes global existence for equations driven by a broad class of general multiplicative noises.
 
We have thus proved the  non-uniqueness in law holds for stochastic Navier-Stokes equation with general  multiplicative noise.
\begin{corollary}\label{coro:noninlaw}
Under Assumption \ref{a:G:1}, non-uniqueness in law holds for stochastic Navier-Stokes equation \eqref{1} for every given initial law supported on divergence free and mean free vector fields in $L^2$.
\end{corollary}

 {In the case of  additive noise and linear multiplicative noise, the pathwise non-uniqueness, non-uniqueness in law and even non-uniqueness of Markov-selections have been established in \cite{HZZ24} and \cite{HZZ23a}. In the general multiplicative case as considered in this paper, our first main result has proved the non-uniqueness in law of this system, thus} we understand stationarity in the sense of shift invariance of laws of solutions on the space of trajectories, see also \cite{BFHM19,FFH21，BFH24} and \cite{HZZ23}.

 More precisely, we define the joint trajectory space for the solution and the driving Wiener process by
$$
\cT = C(\mathbb{R};L^2)\times C(\mathbb{R};U),
$$
 and let $S_t$, $t\in \mathbb{R}$, be  shifts on trajectories given by
$$
S_t(u,W)(\cdot)=(u(\cdot+t),W(\cdot+t)-W(t)),\quad t\in\mR,\quad (u,W)\in\cT.
$$
We note that the shift in the second component acts differently in order to guarantee that for a Wiener process $W$, the shift $S_{t}W$ is again a Wiener process.

Stationary solutions to the stochastic Navier--Stokes equations \eqref{1} are defined as follows.
\begin{definition}\label{d:1.1}
	We say that  $((\tilde{\Omega},\tilde{\mathcal{F}},(\tilde{\mathcal{F}}_{t})_{t\in \mathbb{R}},\tilde{\mathbf{P}}),\tilde{u},\tilde{W})$ is a stationary solution to the stochastic Navier--Stokes equations \eqref{1} provided it satisfies  \eqref{1} in the  following sense:
\begin{enumerate}
\item
$(\tilde{\Omega},\tilde{\mathcal{F}},(\tilde{\mathcal{F}}_{t})_{t\in \mathbb{R}},\tilde{\mathbf{P}})$ is a stochastic basis with complete right-continuous filtration;
\item
$\tilde{W}$ is a $U$-valued, two sided cylindrical Wiener process with respect to the filtration $(\tilde{\mathcal{F}})_{t \in \mathbb{R}}$;
\item the velocity $\tilde{u} \in  C(\mathbb{R};H^{\theta})\ \tilde{\mathbf{P}}$-a.s. for some $\theta>0$ and  is $(\tilde{\mathcal{F}})_{t \in \mathbb{R}}$-adapted,  see Section \ref{sec:notations} for the precise definition of $L^2_{\sigma}$;
\item for every $-\infty < s\leq t<\infty$ it holds  $\tilde{\mathbf{P}}$-a.s.
$$
\begin{aligned}
&\langle \tilde{u}(t),\psi \rangle + \int_{s}^{t}\langle \div(\tilde{u}\otimes \tilde{u}),\psi \rangle \dif r\ =\langle \tilde{u}(s),\psi \rangle +\int_{s}^{t} \langle \Delta \tilde{u}, \psi\rangle \,\dif r +\langle \int_s^t G(\tilde{u}) \dif \tilde{W}_r,\psi\rangle
\end{aligned}
$$
for all $\psi\in C^{\infty}(\mathbb{T}^{3}),$ $\div\psi=0$.
\end{enumerate}
and its law is shift invariant, that is,
	$$\mathcal{L}[S_{t}(\tilde{u},\tilde{W})]=\mathcal{L}[\tilde{u},\tilde{W}]\qquad\text{for all}\quad t\in \mathbb{R}.$$
\end{definition}

 Every stationary solution defines a dynamical system $(\mathcal{T}, \mathcal{B}(\mathcal{T}),(S_t)_{t \in \mathbb{\mathbb{R}}}, \mathcal{L}(u,W))$ in the sense of \cite[Chapter 1]{DZ96}, where $\mathcal{B}(\mathcal{T})$ denotes the Borel $\sigma$-algebra on $\mathcal{T}$. As in \cite{HZZ22}, we give the following definition of ergodic  stationary solutions in the sense of the ergodicity of the associated dynamical system.
\begin{definition}
A stationary solution $((\tilde{\Omega},\tilde{\mathcal{F}},(\tilde{\mathcal{F}}_{t})_{t\in \mathbb{R}},\tilde{\mathbf{P}}), \tilde{u},\tilde{W})$ is ergodic provided 
\begin{equation*}
\mathcal{L}[\tilde{u},\tilde{W}](A)=1\ \ \ \ \text{or}\ \ \ \ \mathcal{L}[\tilde{u},\tilde{W}](A)=0,
\end{equation*}
for all $A \in \mathcal{B}(\mathcal{T})$ and shift invariant.
\end{definition}

As an ingredient in the construction of ergodic stationary solution on the path space, we can construct a probabilistically strong and analytically weak solution with a prescirbed energy profile $e(t)$, $t \in [0, +\infty)$.
\bt\label{thm:main:1}
Suppose Assumption \ref{a:G:1}  holds. For any $r_0 \geq 2$, let $e: [0, +\infty) \to(0,\infty)$ be a smooth function satisfying $\bar e\geq e(t)\geq \underline{e}>0$ where $\underline{e}=\underline{e}_{r_0,\d_0,p_0,L}$ is some large enough constant. There exists an $(\mathcal{F}_t)_{t\geq 0}$-adapted process $u$ which belongs to $ C([0,\infty),H^{\vartheta})\cap C^{\vartheta}([0,\infty),L^{2})$ $\mathbf{P}$-a.s.  for some $\vartheta>0$  and is an analytically weak and probabilistically strong solution to \eqref{1} in the sense of Definition~\ref{d:sol}.
Moreover, the solution satisfies
\begin{align}\label{est:u1}
\sup_{t\ge0}{\E}\sup_{s\in[t,t+1]}\| u(s)\|^{{r_0}}_{H^\vartheta}+\sup_{t\ge0}{\E}\|{u}\|^{{r_0}}_{C^\vartheta([t,t+1], L^2)}<\infty,
\end{align}
and  there holds
\begin{align}\label{eq:K1}
\begin{aligned}
\mathbf{E}\|u(t)\|_{L^2}^2=e(t),\ \ t \in [0,+\infty).
\end{aligned}
\end{align}
There are infinitely many such solutions by choosing different energy profiles $e$.
\et

Using Theorem \ref{thm:main:1} and the Krylov-Bogoliubov argument, we are able to construct ergodic stationary solutions of stochastic Navier-Stokes equation with general multiplicative noise.
\bt\label{thm:main:2}
Suppose Assumption \ref{a:G:1} holds. For any $r_0 \geq 2$, let $K=K_{r_0,\d_0,p_0,L}>0$ be some constant  large enough.  There exists some constant $C>0$ and  an ergodic stationary solution $((\tilde{\Omega},\tilde{\mathcal{F}},\tilde{\mathbf{P}}),\tilde{u},\tilde{W})$ to \eqref{1} in the sense of  Definition \ref{d:1.1} satisfying
\begin{align}\label{eq:s55}
	\begin{aligned}
		\tilde{\mathbf{E}}\| \tilde{u}\|_{L^2}^2= K,
	\end{aligned}
\end{align}
and for some $\vartheta>0$ and for every $N\in\mathbb{N}$, 
\begin{align}\label{eq:s56}
\tilde{\E}\sup_{t\in[-N,N]}\| \tilde{u}(t)\|^{{r_0}}_{H^\vartheta}+\tilde{\E}\|\tilde{u}\|^{{r_0}}_{C^\vartheta([-N,N], L^2)}\leq C N,
\end{align}
where $\tilde{\E}$ is the expectation taken with respect to $\tilde{\mathbf{P}}$. In particular, there exist infinitely many ergodic stationary solutions to \eqref{1} by choosing different $K$.
\et

The paper is organized as follows. In Section \ref{sec:pre}, we first give the precise definitions of the functional spaces we shall work with. Then, the truncation operator $\Pi_{\zeta}$ is introduced to truncate the stochastic part so that the trajectories are uniformly bounded and preserve their Hölder norms. 
In Section \ref{sec:stochastic}, we provide a priori estimate for the stochastic convolution which plays a key role in determining the convex integration parameters in the iteration. 
In Section \ref{sec:stationary}, we prove  Theorem \ref{thm:main:1} and Theorem \ref{thm:main:2} where the presentation of the stochastic convex integration iteration is a little bit different from the traditional one. We refer to Section \ref{sec:setup} for details. In Section \ref{sec:iterative:cut}, we slightly modify the building block in Section \ref{sec:stationary} and Section \ref{sec:iteration}  by introducing a cut-off function in the perturbation and prove our main result Theorem \ref{thm:main:3}.

\subsection{Examples}\label{sec:example}

 We close this section by providing some examples which satisfy Assumption \ref{a:G:1}.

 \subsubsection{Linear multiplicative noise}
One typical example of Assumption \ref{a:G:1} is the linear multiplicative noise. Indeed, we are interested in the following case:
\[
G(u) \dif W_t=u\dif B_t,
\]
 where $B$ is a one-dimensional  standard Brownian motion. We take $U$ in Assumption \ref{a:G:1} to be $\mathbb{R}$ equipped with the usual inner product.  Then the  Hilbert--Schmidt norm of $G$ is this case is $\|G(u)\|_{\mathcal{L}^2_0(U: H^{-2\d_0})}\lesssim \|u\|_{H^{-2\d_0}}$. By duality and the Sobolev embedding, there holds
\begin{equation}\label{ex:1}
\|u\|_{H^{-2\d_0}} \lesssim \|u\|_{L^{6/(3+4\d_0)}} .
\end{equation}
Therefore, we obtain that
\begin{equation*}
\|G(u_1) -G(u_2)\|_{\mathcal{L}^2_0(U:H^{-2\d_0})}\lesssim   \|u_1-u_2\|_{L^{6/(3+4\d_0)}}
\end{equation*}
 for all $\d_0 \in (0,1/2)$. Hence $G$ satisfies Assumption \ref{a:G:1}  with $p_0=\frac6{3+4\d_0}\in[1,2)$.
  \begin{remark}\label{r:lm}
 The linear multiplicative noise is an important type of noise. As implied in \cite{HZZ24}, by performing the transform $v:=u \exp(-B)$, one can derive non-uniqueness in law of solutions to Stochastic Naiver-Stokes equation driven by linear multiplicative noise. However, existence  of probabilistically strong solution with given $L^2$-initial condition is not known in the literature. Also by using such transformation, it seems only a trivial stationary solution can be obtained. Under the framework here, we can prove not only existence of probabilistically strong solutions with given $L^2$-initial condition but also the existence of infinitely many non-trivial stationary solutions. 
 \end{remark}
\subsubsection{Matrix-valued linear multiplicative noise}
Another typical example of Assumption \ref{a:G:1} is the linear Matrix-valued multiplicative noise (see e.g. \cite{HR23}). Indeed, we are interested in considering the following case:  for $k=1,2,3$
\[
(G(u) \dif W_t)^k= \sum_{i=1}^3u^i\dif B^{i,k}_t,
\]
 where $\{B^{i,k}\}_{i,k=1,2,3}$ is a family of one-dimensional independent standard Brownian motion.  This system is difficult to be transformed into random PDEs through conventional methods. We take $U$ in Assumption \ref{a:G:1} to be $\mathbb{R}^{3\times 3}$ equipped with the usual inner product.   By the same calculation as \eqref{ex:1},  $G$ satisfies Assumption \ref{a:G:1}  with $p_0=\frac6{3+4\d_0}\in[1,2)$.

\subsubsection{Multiplicative noise with Lipschitz continuous coefficients}
Let $U:=L^2(\mathbb{T}^3;\mathbb{R}^3)$ and let $\{e_k\}_{k \in \mZ^3}$ be the eigenfunctions of the Lapalcian on $\mathbb{T}^3$ with periodic boundary condtions. Then the eigenfunctions are equi-bounded in the sense that $\sup_{k \in \mZ^3} \|e_k\|_{L^{\infty}} < \infty.$ The cylindrical Wiener process $W(t,x)$ is defined by 
\begin{equation*}
W(t,x):= \sum_{k \in \mZ^3} e_k(x) \beta_k(t),
\end{equation*}
where $\{\beta_k\}_{k \in \mathbb{Z}^3}$ is family of independent real-valued Brownian motions on $[0,+\infty)$. 
 Let $\{g_k\}_{k \in \mathbb{Z}^3}$ be a sequence of functions such that for all $k \in \mathbb{Z}^3$, $g_k:\mathbb{R}^3 \to \mathbb{R}^3$ and there exists a sequence $\{L_k\}_{k \in \mathbb{Z}^3}$ such that
\begin{equation*}
|g_k(x)-g_k(y)|\leq L_k |x-y|,
\end{equation*}
and 
\begin{equation}\label{eq:gk:lg}
|g_k(x)| \leq L_k(1+ |x|).
\end{equation}
Let $\{\sigma_k\}_{k \in \mathbb{Z}^3}$ be a sequence such that
\begin{equation*}
\sum_{k \in \mathbb{Z}^3} \sigma_k^2 L_k^2 <\infty.
\end{equation*}

 For  $\d_0\in (0,1/2),p_0 \in [1,2)$,  the  Nemytskii operator ${G}:L^{p_0}\to\mathcal{L}_0^2(U; H^{-2\d_0})$ is defined as follows: for $u \in L^{p_0}$ and $f\in C^\infty(\mathbb{T}^3;\mathbb{R}^3)$, we define 
\begin{equation*}
({G}(u)f)(\xi):=\sum_{k \in \mZ^3}\sigma_kg_k(u(\xi))\langle f, e_k \rangle e_k(\xi),\ \ \  \xi \in \mathbb{T}^3.
\end{equation*}
Here $\langle\cdot,\cdot\rangle$ denotes the $L^2$-inner product.

{\begin{proposition}
 For any $\d_0\in (0,1/2)$, the operator $G$ defined above 
 satisfies Assumption \ref{a:G:1}  with $p_0=\frac6{3+4\d_0}\in[1,2)$.  
\end{proposition}}
\begin{proof}
It is east to see that
\begin{equation*}
\|{G}(u)\|_{\mathcal{L}^2_0(U;H^{-2\d_0})}^2 \leq \sum_{k \in \mZ^3} \|G(u)e_k\|_{H^{-2\d_0}}^2= \sum_{k \in \mZ^3} \|\sigma_k g_k(u)e_k\|_{H^{-2\d_0}}^2.
\end{equation*}
By \eqref{ex:1},  \eqref{eq:gk:lg} and the equiboundedness of the eigenfunctions, we have
\begin{equation*}
\begin{aligned}
\|G(u)\|_{\mathcal{L}^2_0(U;H^{-2\d_0})}^2   \lesssim  \sum_{k \in \mZ^3} \|\sigma_k g_k(u)\|_{L^{6/(3+4\d_0)}}^2 \lesssim \sum_{k \in \mZ^3} \sigma_k^2 L_k^2 (\|u\|_{L^{6/(3+4\d_0)}}^2+1).
\end{aligned}
\end{equation*}
Similarly, we obtain
\begin{equation*}
\|G(u)-G(v)\|_{\mathcal{L}^2_0(U;H^{-2\d_0})}^2 \leq \sum_{k \in \mZ^3} \sigma_k^2 L_k^2 \|u-v\|_{L^{6/(3+4\d_0)}}^2,
\end{equation*}
which implies Assumption \ref{a:G:1}.
\end{proof}

\section{Preliminaries}\label{sec:pre}
\subsection{Notations}\label{sec:notations}
We let $\mN_{0}:=\mN\cup \{0\}$. By $\mathbb{P}$ we denote the Helmholtz projection. In the following, we employ the nontation $a\lesssim b $  if $a \leq cb$ for some harmless constant $c>0$. Let
$\mathcal{L}^2_0(U_1;U_2)$ be the space of Hilbert-Schmidt operators from Hilbert spaces $U_1$ to $U_2$. We define $L^{2}_{\sigma}$ by
\[
L^{2}_{\sigma}=\{f\in L^2; \int_{\mathbb{T}^{3}} f\,\dif x=0,\div f=0\},
\] 
with the usual $L^2$ norm.  We use $L^p$ to denote the set of  standard $L^p$-integrable functions from $\mathbb{T}^3$ to $\mathbb{R}^3$ with the norm $\|f\|_{L^p}$. For $s>0$, $p>1$ we set 
\[
W^{s,p}:=\{f\in L^p; \|(I-\Delta)^{{s}/{2}}f\|_{L^p}< \infty\}
\] 
with the norm  $\|f\|_{W^{s,p}}=\|(I-\Delta)^{{s}/{2}}f\|_{L^p}$. 
For $s>0$, we define {$H^s:=W^{s,2}$}. For $s<0$ we define $H^s$ to be the dual space of $H^{-s}$.

Let $E$ be a Banach space with norm $\|\cdot\|_E$. For a domain  $D \subseteq \R$, $\a\in(0,1)$ and $t\in\mR$, we write $C_D^{\a}E=C^{\a}(D;E)$ for the space of $\a$-H\"older continuous functions from $D$ to $E$, equipped with the norm
\[
\|f\|_{C^\alpha_DE}=\sup_{s,r\in D,s\neq r}\frac{\|f(s)-f(r)\|_E}{|r-s|^\alpha}+\sup_{s\in D}\|f(s)\|_{E}.
\] 
In particular, for $\a \in (0,1)$ and $t \in \R$, we denote $C_{[t,t+1]}^{\a}E$ by $C_t^{\a}E$ for short.

Let  $D\subset\R$ be a domain and $N\in\N_{0}$, we denote by  $C^{N}_{D,x}$ the space of $C^{N}$-functions on $D\times\mathbb{T}^{3}$. The spaces are equipped with the norms
$$
 \|f\|_{C^N_{D,x}}=\sum_{\substack{0\leq n+|\alpha|\leq N\\ n\in\N_{0},\alpha\in\N^{3}_{0}}}\sup_{t\in D}\|\partial_t^n D^\alpha f\|_{L^\infty}.
$$
In particular, we denote  $C_{[t,t+1],x}^N$ by $C_{t, x}^N$, $t \in \mathbb{R}$. 

Next, we introduce norms for function spaces of random variables on $\Omega$ taking values in $C(\R;L^{2})$ and {$C^{1}(\R\times\mathbb{T}^{3})$}, respectively, with  bounds in $L^{p}(\Omega;C(I,L^{2}))$ and $L^{p}(\Omega;{C^{1}(I\times\mathbb{T}^{3}))}$ for any bounded interval $I\subset \mR$.

 For $p\in[1,\infty)$ and $t_0 \in \mathbb{R}$, we denote
\begin{equation*}
\$u\$^p_{L^2,p,t_0}:=\sup_{t \geq t_0}\mathbf{E}\left[\sup_{t\leq s\leq t+1}\|u(s)\|^p_{L^2}\right],\quad \$u\$^p_{C^1_{t,x},p,t_0}:=\sup_{t \geq t_0}\mathbf{E}\left[\|u(s)\|^p_{C^1_{[t,t+1],x}}\right],
\end{equation*}
and in particular we write $\$u\$^p_{L^2,p}$ and $\$u\$^p_{C^1_{t,x},p}$ for short when $t_0=0$. Similarly, we define the corresponding norms with $L^2$ replaced by $L^p$, $H^\vartheta$ and $C^1_{t,x}$ replaced by $C_t^\vartheta L^2$, $C_tW^{1,p}$.  

For any interval ${[a,b]} \subseteq \R$ and $p \in [1, \infty)$, we denote
 \[
 \$u\$_{E,p,{[a,b]}}^p= \E[\sup_{t \in {[a,b]}}\|u(t)\|_E^p ],
 \]
 where $E$ is some Banach space which can be  $L^p$, $H^\vartheta$  or $C^1_{t,x}$.
 
  It's well-known that $H^{2\d}$, $\d \in \mathbb{R}$, has equivalent norm with $D((-\Delta)^{\d})$. Let $S(t)$ be the analytic semigroup generated by $\Delta$. For $\d_1< \d_2 $,  $S(t)$ satisfies the following heat kernel smoothing effect 
\begin{equation}\label{eq:asg}
 \lv S_{}(t)\rv_{\mathcal{L}(H^{\d_1}; H^{\d_2})} \leq C_{\d_1,\d_2} \ t^{-(\d_2-\d_1)/2}\ \ ,\ \ \ \forall t>0.
\end{equation}

\subsection{The operator $\Pi_{\zeta}$}\label{sec:operator}
 Let $f \in L^2$ and $\{e_k\}_{k \in \mathbb{Z}^3}$ be the sequence of eigenvectors of the Laplacian on $\mathbb{T}^3$ with periodic boundary condition whose eigenvalues are $\{-|k|^2\}_{k \in \mathbb{Z}^3}$.  The Fourier series decomposition for $f$ reads as $f= \sum_{k \in \mathbb{Z}^3} f_k e_k$, and the truncation operator $\Pi_{\zeta},\zeta\in\mN$ is defined by
\begin{equation*}
\begin{aligned}
\Pi_{\zeta}f:= \sum_{|k|^2 \leq \zeta,k\in\mZ^3} f_{k,\zeta} \ e_k,\ \ \ \ \  f_{k,\zeta}:=(f_k\vee(-\zeta))\wedge\zeta.
\end{aligned}
\end{equation*}
 It's easy to prove that $\Pi_{\zeta}$ satisfies the following properties.
\begin{lemma}\label{l:truncate}
For any $\g \in [0,1)$, $f,g\in {H^\gamma}$ and $\zeta \geq 1$, it holds that
\begin{equation*}
{\| \Pi_{\zeta} f\|_{L^\infty}}\lesssim \zeta^4\ \ \text{and}\ \ \|\Pi_{\zeta} f-\Pi_{\zeta} g \|_{H^{\g}} \leq \|f-g\|_{H^{\g}}.
\end{equation*}
As a result, for any $\gamma_1 \in (0,1), a<b\in\mathbb{R}$ and $f \in C^{\gamma_1}((a,b) ;H^{\gamma})$, there holds
\begin{equation*}
\|\Pi_{\zeta} f \|_{C^{\gamma_1}((a,b) ;H^{\gamma})} \leq  \|f\|_{C^{\gamma_1}((a,b) ;H^{\gamma})}.
\end{equation*}
\end{lemma}
\begin{proof}
Since $\|e_k\|_{\infty} \leq 1$ for any $k \in \mathbb{Z}^3$ and $|f_{k,\zeta}|\leq \zeta$, we have
\begin{equation*}
\| \Pi_{\zeta} f\|_{L^\infty} =\big\|\sum_{|k|^2 \leq \zeta} f_{k,\zeta} \ e_k\big\|_{L^\infty} \leq \zeta  \sum_{|k|^2 \leq \zeta}\big\| e_k\big\|_{L^\infty}\lesssim \zeta^4,
\end{equation*}
where in the last inequality we used the number of $k$'s such that $|k|^2 \leq \zeta$ is no greater than $27\zeta^3$.

For the second statement, by the Plancherel identity, we have
\begin{equation*}
\|\Pi_{\zeta} f-\Pi_{\zeta} g\|_{H^{\g}}^2= \sum_{|k|^2 \leq {\zeta}} |k|^{2\g}|f_{k,\zeta}- g_{k,\zeta}|^2\leq \|f-g\|_{H^{\g}}^2,
\end{equation*}
where we used  $ |f_{k,\zeta}- g_{k,\zeta}| \leq |f_k -g_k|$.

\end{proof}

\subsection{Inverse divergence operators $\mathcal{R}$ and $\mathcal{B}$}\label{sec:div-1}
 We recall the inverse divergence operator $\mathcal{R}$ from \cite[Section 5.6]{BV19}, which acts on vector fields $v$ with $\int_{\mathbb{T}^3}v\dif x=0$ as
\begin{equation*}
	(\mathcal{R}v)^{kl}=(\partial_k\Delta^{-1}v^l+\partial_l\Delta^{-1}v^k)-\frac{1}{2}(\delta_{kl}+\partial_k\partial_l\Delta^{-1})\div\Delta^{-1}v,
\end{equation*}
for $k,l\in\{1,2,3\}$. Then $\mathcal{R}v(x)$ is a symmetric trace-free matrix for each $x\in\mathbb{T}^3$, and $\mathcal{R}$ is a right inverse of the div operator, i.e. $\div(\mathcal{R} v)=v$. By \cite[Theorem B.3]{CL20} we know for $\sigma\in\mathbb{N},p\in[1,\infty]$
\begin{align}\label{eR}\|\mathcal{R}f(\sigma\cdot)\|_{L^p}\lesssim \sigma^{-1}\|f\|_{L^p}.\end{align}

By $\cS^{3\times 3}$ we denote the set of symmetric $3\times 3$ matrices and by $\cS_0^{3\times 3}$  the set of  symmetric trace-free matrices. Let  $C_0^\infty(\mT^3,\mathbb R^{3\times 3})$ be the set of periodic smooth matrix valued functions with zero mean. We also introduce the bilinear version $\cB:C^\infty(\mT^3,\mR^3)\times C_0^\infty(\mT^3,\mathbb R^{3\times 3})\to C^\infty(\mT^3,\cS_0^{3\times 3})$ as in \cite[Section~B.3]{CL20} by
\begin{align*}
	\cB(v,S)=v\cR S-\cR(\nabla v\cR S).
\end{align*}
Then by \cite[Theorem B.4]{CL20} we have $\div(\mathcal{B} (v,S))=vS-\frac1{(2\pi)^3}\int_{\mT^3} vS\dif x $ and for $p\in[1,\infty]$
\begin{align}\label{eB}
	\|\cB(v,S)\|_{L^p}\lesssim \|v\|_{C^1}\|\mathcal{R}S\|_{L^p}.
\end{align}

 \section{Stochastic convolution}\label{sec:stochastic} 
To deal with the Navier-Stokes equation \eqref{1}, we first divide \eqref{1} into two equations:  the following stochastic differential equation:
\begin{equation}\label{eq:sc}
\begin{aligned}
\dif z(t)-(\Delta - cI) z(t) \dif t+\nabla p_z\dif t& =  G(v(t)+ z(t)) \dif W_t,\\ \div z&=0,\\ 
z(0)&=0, 
\end{aligned}
\end{equation}
and  the following random partial differential equation:
\begin{equation*}
\aligned
\partial_tv-c z-\Delta v +\div((v+z)\otimes (v+z))+\nabla p_v&=0,
\\
\div  v&=0,\\
v(0)&=u_0,
\endaligned
\end{equation*} 
where the constant $c>0$ is the dissipative constant so that a uniform-in-time moment estimate for $z$ is possible.
Here, the first equation \eqref{eq:sc} will be estimated  using the It$\rm \hat{o}$ calculus, while the latter will be handled by the convex integration.

In this section, we assume that $G$ satisfies Assumption \ref{a:G:1} with some $(\delta_0,p_0)$. Since $\mP$ and $\Delta$ commute on torus, the operator  $\mathbb{P} \circ G$ 
 also satisfies Assumption \ref{a:G:1} with the same $(\delta_0,p_0)$.

We { provided a priori estimate of \eqref{eq:sc} by assuming} that $v \in L^{r}(\Omega; C([0, \infty);L^2))$ is  a given $(\mathcal{F}_t)_{t \geq 0}$-adapted process with some  $r>2$.  By a standard contraction argument \cite[Chapter 7]{DZ14}, there exists a unique mild solution $z_c \in L^{\infty}_{\rm loc}([0,+\infty); L^2(\Omega;L^2_{\sigma}))$ to \eqref{eq:sc} which has the following form
\begin{equation*}
z_c(t)=  \int_0^t e^{-c(t-s)} S_{}(t-s) {\mathbb{P}}\circ G(v(s)+ z_c(s)) \dif W_s.
\end{equation*}
where we recall that $\mathbb{P}$ denotes the Helmholtz projection and $S(t)$ represents the heat semigroup.

\bl \label{l:z:holder}
Suppose Assumption \ref{a:G:1}  holds with some $(\delta_0,p_0)$ and suppose $v \in L^{r}(\Omega; C([0, \infty);L^2))$ with $r>2$.  Let  $\gamma $, $\d$ be  positive parameters such that
\[ 
0<\gamma+ \delta + \d_0< \frac12-\frac 2r.
\]
 There exists $c>0$ large enough and some constant $C_{c,r,\d,\d_0,\g,L}$  such that
\begin{equation*}
\sup_{t \geq 0}\mathbf{E}[\|z_c\|^r_{C^{\gamma}_{t}H^{2\d}}] \leq C_{c,r,\d,\d_0,\g,L}  (1+ \$v\$_{L^{p_0},r}^r).
\end{equation*}
\el

\begin{proof}

We first prove that there exists some constant $C_{r,\d,\d_0,L}$ such that 
\begin{equation}\label{eq:z:bound1}
\sup_{t \geq 0}\mathbf{E}[\|z_c(t)\|^r_{H^{2\d}}]  \leq C_{r,\d,\d_0,L}\  c^{r(\d+\d_0-1/2)}  (1+\$ v\$_{L^{p_0},r}^r ).
\end{equation}
By the B\"urkholder-Davis-Gundy inequality, \eqref{eq:asg} and  \eqref{eq:G:lg}, there exists some constant $C_r$ such that
\begin{equation*}
\aligned
&\mathbf{E}[\|z_c(t)\|^r_{H^{2\d}}] = \mathbf{E}\left\Vert \int_0^t  e^{-c(t-s)}S_{}(t-s){\mathbb{P}} \circ G(v(s)+ z_c(s))  \dif W_s\right\Vert_{H^{2\d}}^r\\
\leq&C_r\  L^r\ \mathbf{E}\( \int_0^t (t-s)^{-2(\delta+\d_0)} e^{-2c(t-s)}   (1+ \lv v(s)\rv_{L^{p_0}}^2 + \lv z_c(s)\rv_{L^{p_0}}^2) \dif s  \)^{r/2}. 
\endaligned
\end{equation*}
By the H\"older  inequality, it holds that
\begin{equation*}
\aligned
\mathbf{E}[\|z_c(t)\|_{H^{2\d}}^r]
\leq& C_r\  L^r\ \mathbf{E} \int_0^t e^{-c r(t-s)/2}  (1+ \lv v(s)\rv_{L^{p_0}}^r + \lv z_c(s)\rv_{L^{p_0}}^r) \dif s\\
&\ \ \ \ \ \ \ \ \ \ \ \ \ \ \ \ \ \ \ \ \ \ \ \ \ \ \ \ \ \ \ \ \ \ \  \( \int_0^t  (t-s)^{-2(\delta+\d_0) r/(r-2)} e^{-cr(t-s)/(r-2)} \dif s \)^{(r-2)/2}  . 
\endaligned
\end{equation*}
Since $0<\delta +\delta_0 <1/2-1/r$, there exists some constant $C_{r,\d,\d_0,L}$ such that
\begin{equation*}
\begin{aligned}
\sup_{s\in [0,t]}\mathbf{E}[\|z_c(s)\|^r_{H^{2\d}}] \leq&\  C_{r,\d,\d_0,L}\  c^{r(\d+\d_0-1/2)}   (1+\$ v\$_{L^{p_0},r}^r +  \sup_{s \in[0,t]}\mathbf{E}[\lv z_c(s)\rv_{H^{2\d}}^r] ).
\end{aligned}
\end{equation*}
By choosing $c$ large enough such that 
\begin{equation*}
C_{r,\d,\d_0,L}\  c^{r(\d+\d_0-1/2)} <1/2,
\end{equation*}
 there holds
\begin{equation*}
\sup_{s\in [0,t]}\mathbf{E}[\|z_c(s)\|^r_{H^{2\d}}] \leq 
1+\$ v\$_{L^{p_0},r}^r .
\end{equation*}
Since the right hand side is independent of $t$, there exists some constant $C_{p,\d,\d_0,L}$ such that \eqref{eq:z:bound1} holds for $c$ sufficiently large.

{Next, for any positive $\g$, $\d$, $r$ such that $0< \gamma+ \delta + \d_0< 1/2-2/r$, we define $\g_0=1/2-1/r-\delta -\d_0-\epsilon$, where we choose $\epsilon>0$ small enough such that $\gamma+\frac1r<\g_0$.} Then we prove that there exists some constant $C_{c,r,\d,\d_0,L}$ such that   for any $0\leq s \leq t $
\begin{equation}\label{eq:holder:target}
\mathbf{E}[\|z_c(t)-z_c(s)\|_{H^{2\d}}^r] \leq C_{c,r,\d,\d_0,L}  (t-s) ^{\gamma_0 r} (1+ \$v\$_{L^{p_0},r}^r).
\end{equation}

For any $0\leq s \leq t $ satisfying $|t-s|\leq1$, we have the following decomposition
\begin{equation*}
\aligned
\|z_c(t)-z_c(s)\|_{H^{2\d}}^r \leq
&  3^{r-1} \left\|\int_s^t  e^{-c(t-u)} S_{}(t-u) {\mathbb{P}} \circ G(v(u)+z_c(u)) \dif W_u\right\|_{H^{2\d}}^r\\
& + 3^{r-1} \left\|\int_{0}^s  (S(t-s)-I)e^{-c(s-u)} S(s-u) {\mathbb{P}} \circ G(v(u)+z_c(u)) \dif W_u\right\|_{H^{2\d}}^r\\
& + 3^{r-1} \left\|\int_{0}^s  (e^{-c(t-u)}- e^{-c(s-u)}) S(t-u) {\mathbb{P}}\circ G(v(u)+z_c(u)) \dif W_u\right\|_{H^{2\d}}^r\\
:=&3^{r-1}(I_1 + I_2+ I_3).
\endaligned
\end{equation*}

For the first term, by the B\"urkholder-Davis-Gundy inequality, there holds
\begin{equation*}
\aligned
\mathbf{E}[I_1]
\leq& C_r\mathbf{E} \(\int_s^t  \lv e^{-c(t-u)} S(t-u) G(v(u)+z(u))  \rv_{\mathcal{L}^2_0(U;H^{2\d})}^2 \dif u \)^{r/2}\\
\leq& C_r\mathbf{E} \(\int_s^t  \lv S(t-u) \rv^2_{\mathcal{L}(H^{-2\d_0}; H^{2\d})} \lv G(v(u)+z(u))  \rv_{\mathcal{L}^2_0(U; H^{-2\d_0})}^2\dif u \)^{r/2}.
\endaligned
\end{equation*}
By \eqref{eq:asg} and \eqref{eq:G:lip} and H$\ddot{\text{o}}$lder inequality, there holds
\begin{equation*}
\aligned
\mathbf{E} [I_1] \leq& C_r L^r\ \mathbf{E}\big(\int_s^t (t-u)^{-2(\delta+\d_0)}  (1+ \|v(u)\|_{L^{p_0}}^{2} + \|z(u)\|^2_{L^2}) \dif u\big)^{r/2}\\
\leq&C_r L^r\ \big(\int_s^t (t-u)^{-2r(\delta+\d_0)/(r-2)} \dif u \big)^{(r-2)/2} \mathbf{E} \int_s^t  (1+ \|v(u)\|^r_{L^{p_0}} + \|z(u)\|^r_{L^2}) \dif u\\
\leq&C_r L^r (t-s)^{r/2-r(\delta+\d_0)}  (1+ \$v\$_{L^{p_0},r}^r +\sup_{u \in [s,t]} \mathbf{E}[\|z(u)\|_{L^2}^r]).
\endaligned
\end{equation*}
By \eqref{eq:z:bound1}, we obtain that, for any  $\delta + \d_0 <1/2$, there exists some constant $C_{r,\d,\d_0,L}$ such that
\begin{equation}\label{eq:cont:1}
\mathbf{E}[I_1] \leq {C_{r,\d,\d_0,L}  
(t-s)^{\gamma_0 r}  (1+ \$v\$_{L^{p_0},r}^r ).}
\end{equation}

For the second term, by the B\"urkholder-Davis-Gundy inequality again, there holds
\begin{equation*}
\mathbf{E}[I_2] \leq  C_r \ \|S(t-s)-I\|_{\mathcal{L}(H^{2\d+\gamma_0}; H^{2\d})}^r \mathbf{E} \( \int_{0}^s e^{-2c(s-u)} \|S(s-u) G(v(u)+ z(u))\|_{\mathcal{L}^2_0(U; H^{2\d+\gamma_0})}^2 \dif u \)^{r/2}.
\end{equation*}
By analyticity of the semigroup $S_{}(t)$ and \eqref{eq:asg}, we have
\begin{equation*}
\aligned
\|S(t-s)-I\|_{\mathcal{L}(H^{2\d+2\gamma_0}; H^{2\d})}  \leq&\left \|\int_{0}^{t-s} (-A)^{1-\gamma_0}S(v) dv \right\|_{\mathcal{L}(H^{2\d})} \leq C_{\gamma_0} (t-s)^{\gamma_0},
\endaligned
\end{equation*}
and
\begin{equation*}
\begin{aligned}
\|S(s-u) G(v(u)+ z(u))\|_{\mathcal{L}^2_0(U; H^{2\d+2\gamma_0})} \leq& \|S(s-u)\|_{\mathcal{L}(H^{-2\d_0};H^{2\d+2\gamma_0})}\|G(v(u)+ z(u))\|_{\mathcal{L}^2_0(U; H^{-2\d_0})}\\
\leq& C_{\d,\d_0,\gamma_0} (s-u)^{-(\d+\d_0+\gamma_0)} \|G(v(u)+ z(u))\|_{\mathcal{L}^2_0(U; H^{-2\d_0})}.
\end{aligned}
\end{equation*}
By Hölder inequality and the condition that $\gamma_0+\d+\d_0<1/2-1/r$, we have
\[
\begin{aligned}
\mathbf{E}[I_2] \leq& C_{r,\d,\d_0,\gamma_0} (t-s)^{\gamma_0 r} \(\int_0^s e^{-\frac{2c(s-u)r}{r-2}}(s-u)^{-\frac{2(\d+\d_0+\gamma_0)r}{r-2}} \dif u\)^{\frac{r-2}{2}} \E  \|G(v(u)+ z(u))\|_{\mathcal{L}^2_0(U; H^{-2\d_0})}^r\\
\leq& C_{c,r,\d,\d_0,\gamma_0} (t-s)^{\gamma_0 r}(s\wedge1)^{-r(\d_0+\d+\gamma_0)+r/2- 1} \E  \|G(v(u)+ z(u))\|_{\mathcal{L}^2_0(U; H^{-2\d_0})}^r.
\end{aligned}
\]
By   \eqref{eq:z:bound1}, \eqref{eq:G:lg}, there exists some constants $C_{c,r,\d,\d_0,\gamma_0,L}$ such that
\begin{equation}\label{eq:cont:2}
\mathbf{E}[I_2] \leq C_{c,r,\d,\d_0,\gamma_0,L}(t-s)^{\gamma_0 r}(1+ \$v\$_{L^{p_0},r}^r ).
\end{equation} 

For the third term, we can prove similarly that 
\begin{equation*}
\begin{aligned}
\mathbf{E}[I_3] 
\leq& {|e^{-c(t-s)}- 1|^r} \left\|\int_{0}^s  e^{-c(s-u)} S(t-u) G(v(u)+z_c(u)) \dif W_u\right\|_{H^{2\d}}^r\\
\leq& (c(t-s)\wedge 1)^r \(\int_0^s e^{-2c(s-u)}\|S(t-u)G(v(u)+z_c(u))\|_{\mathcal{L}^2_0(U;H^{2\d})}^2 \dif u\)^{r/2}\\
 \leq&(c(t-s)\wedge 1)^r C_{\d,\d_0}(t-s)^{-(\d_0+\d)r} \(\int_0^s e^{-2c(s-u)r/(r-2)} \dif u\)^{(r-2)/2} C_{r,L}(1+ \$v\$_{L^{p_0},r}^r ).
\end{aligned}
\end{equation*}
Since $\gamma_0<1/2-1/r-\d-\d_0<1-\d-\d_0$, there exists some constant $C_{c,r,\d,\d_0,\gamma_0,L}$ such that
\begin{equation}\label{eq:cont:3}
\mathbf{E}[I_3] \leq C_{c,r,\d,\d_0,\gamma_0,L}(t-s)^{\gamma_0 r}(1+ \$v\$_{L^{p_0},r}^r ).
\end{equation} 

Combing \eqref{eq:cont:1}, \eqref{eq:cont:2} and \eqref{eq:cont:3},  we obtain \eqref{eq:holder:target}.

By Kolmogorov's continuity theorem, since  $\gamma_0 r>1$ and $\gamma<\g_0-\frac1r$, there holds
\[
\E[\|z_c(t)\|_{C_t^{{\g}}H^{2\d}}^r] \leq C_{c,r,\d,\d_0,\g,L} (1+ \$v\$_{L^{p_0},r}^r).
\]

\end{proof}

Let $v_1,v_2 \in  L^{r}(\Omega; C([0, \infty);L^2_{}))$ be two given $(\mathcal{F}_t)_{t \geq 0}$-adapted process for some $r>2$. Let $z_{c,i}$ be the solutions of the following stochastic partial differential equation for $i=1,2$,
\begin{equation*}
\dif  z_{c,i}(t) =  (\Delta - cI) z_{c,i}(t) \dif t+\nabla p_i(t)\dif t + G(v_i(s)+ z_{c,i}(t)) \dif W_t,\ \ \ \ z_{c,i}(0)=0,
\end{equation*}
where $G$ satisfies Assumption \ref{a:G:1}.  We can prove the following lemma exactly the same as the proof of Lemma \ref{l:z:holder}.
\bl\label{l:zv:continuous}
{Suppose Assumption \ref{a:G:1}  holds with some $(\delta_0,p_0)$ and suppose $v_i \in L^{r}(\Omega; C([0, \infty);L^2))$ with $r>2,i=1,2$.  Let  $\gamma $, $\d$ be  positive parameters such that
\[
\ 0<\gamma+ \delta + \d_0< 1/2-2/r.
\]

 There exists $c>0$ large enough and some constant $C_{c,r,\d,\d_0,\g,L}$  such that}
\begin{equation}
\sup_{t \geq 0}\mathbf{E}[\|z_{c,1}-z_{c,2}\|^r_{C_t^{\g}H^{2\d}}]  \leq C_{c,r,\d,\d_0,\g,L}  \$ v_1-v_2\$_{L^{p_0},r}^r .
\end{equation}
\el

\section{Construction of non-unique stationary solutions
}\label{sec:stationary}
In this section,  we construct an ergodic stationary solution to the Navier-Stokes equations \eqref{1} and Theorem \ref{thm:main:2} is established. Our approach proceeds in two key steps: Firstly, we build a solution to \eqref{1}  that satisfies ${\mathbf{E}}\| {u}(t)\|_{L^2}^2= e(t)$ for a prescribed energy profile $e(t)$ by using stochastic convex integration. Secondly, by imposing the condition that $e(t)$ is some constant, the existence of an (ergodic) stationary solution follows immediately from the standard Krylov-Bogoliubov argument.

\subsection{Notations on parameters and  iteration set up}\label{sec:setup}

  It seems that standard techniques of transforming the stochastic Navier-Stokes equation \eqref{1} into a random PDE break down. To overcome this challenge, in this paper  we apply a stochastic convex integration method, where each iteration step involves  additionally  solving a SPDEs driven by a nonlinear  It${\rm\hat{o}}$ noise. 
More precisely, at each step $q \in \mathbb{N}$, the pair $(v_q, \mathring{R}_q)$ solves the following equation:
\begin{equation}\label{induction ps}
\aligned
\partial_tv_{q}-c \bar{z}_{q}-\Delta v_{q} +\div((v_{q}+\bar{z}_{q})\otimes (v_{q}+\bar{z}_{q}))+\nabla p_{q}&=\div \mathring{R}_{q}\ ,
\\
\div \ v_{q}&=0\ ,\\
\endaligned
\end{equation} 
where $\bar{z}_{q}=\Pi_{\zeta_{q}}(z_{q})$ for some truncation parameter $\zeta_q>0$ and $z_{q}$ satisfies 
\begin{equation}\label{eq:stochastic:1}
\begin{aligned}
\dif  {z}_{q}-(\Delta-c I) {z}_{q} \dif t+{\nabla p_{z,q}}\dif t=&\ G(v_{q-1}+ {z}_{q}) \dif W_t\ , \\
{\div z_q=}&0\ ,\\
   z_{q}(0)=&0\ ,
\end{aligned}
\end{equation}
for $t \geq 0$ and $z_{q}(t)=0$ for $t <0$.  Here $v_{q-1}$ is the velocity field at the step $q-1$.
The dissipative constant $c$ is chosen to be large enough so that we can apply the results in Section \ref{sec:stochastic} (see the beginning of Section \ref{sec:iteration} and Section \ref{sec:iterative:cut}).

 Our main iteration works as follows. At each step $q$, we have in hand a pair $(v_q, \mathring{R}_q)$  which solves  \eqref{induction ps}. By solving \eqref{eq:stochastic:1} at step $q+1$, we obtain a new stochastic term $z_{q+1}$. Then,  we can construct a new velocity field $v_{q+1}$ and a new stress term $ \mathring{R}_{q+1}$ which satisfies \eqref{induction ps} at step $q+1$ with stochastic term $\bar{z}_{q+1}$ by the stochastic convex integration scheme. Repeating the above procedure gives us a sequence of pairs $(v_{q+1}, \mathring{R}_{q+1})$.  
 Here the previous pathwise analysis breaks down due to the lack of pathwise bound on the stochastic term. To circumvent this difficulty, we employ the moment estimates established in Section \ref{sec:stationary} to control the stochastic noise contributions.

Let $e: [0, +\infty) \to(0,\infty)$ be a smooth function satisfying $\bar e\geq e(t)\geq \underline{e}>0$, where $\underline{e}$ is some large enough constant determined later.  Moreover, we extend the energy profile $e(t)$ by 
$e(t)=e(0)$ for $ t \leq 0.$


{In Section \ref{sec:stationary} and Section \ref{sec:iteration}, Assumption \ref{a:G:1} holds for some $(p_0,\d_0)$.} We  assume that the positive parameters $\g$, $r_0$ satisfy the following condition 
\begin{equation}\label{def:gamma:p}
\ 2\g + \d_0 \in \(0,\frac{1}{2}-\frac{1}{r_0}\)\ .
\end{equation}
We also let the dissipative constant $c$ be large enough so that Lemma \ref{l:z:holder} and \ref{l:zv:continuous} holds. To be more specific, {for any $q\in\mN$}, there holds
\begin{align}
\sup_{t \in \mathbb{R}}\E[\|z_q\|_{C_t^{\g}H^{\g}}^{2r_0}] \leq C_{c,r_0,\g,\d_0,L}(1+ \$v_{q-1}\$_{L^{p_0},2r_0}^{2r_0}),\notag\\
\sup_{t \in \mathbb{R}}\E[\|z_q\|_{C_t^{\g}H^{\g}}^{4r_0}] \leq C_{c,r_0,\g,\d_0,L}(1+ \$v_{q-1}\$_{L^{p_0},4r_0}^{4r_0}),\notag\\
{\sup_{t \in \mathbb{R}}\E[\|z_{q+1}-z_q\|_{C_t^{\g}H^{\g}}^{2r_0}] \leq C_{c,r_0,\g,\d_0,L} \$v_q-v_{q-1}\$_{L^{p_0},2r_0}^{2r_0}.\label{zq:moment}}
\end{align}

 From now on, the parameters $p_0,\d_0,\gamma,r_0$ and $c$ are fixed. Let  $\a_0\geq 12r_0$ (cf. Lemma \ref{l:iterative:R:m}), $M_0$ (cf. Lemma \ref{l:M:bound}), $N_0$ (cf. \eqref{eq:N0}) and $\Xi$ (cf. \eqref{eq:eps} and \eqref{eq:Xi}) be constants which are chosen to be sufficiently large. 

We are now ready to introduce some parameters and terms which will be used in the convex integration scheme.
\begin{definition}\label{def:para:2}
We define 
 \begin{equation*}
 \eps_{q} = 2^{-q-3}\ \ ,\ \ q=-3,-2,-1,0,1.
 \end{equation*}
For $q \geq 2$, we define $\eps_{q}$ by
\begin{equation*}
\eps_{q} =\Xi^{-\a_0{M_0}^{q-2}}.
 \end{equation*}
  By choosing $\Xi$ to be large enough, there holds 
\begin{equation}\label{eq:eps}
\eps_{q} \leq  {\eps_{q-1}}/{2}\ \ ,\ \ \forall q\geq -1.
\end{equation}
 \end{definition}
 
 \begin{definition}\label{def:para:3}
 For $q \geq 1$, we define $\Lambda_{q}$ by
\begin{equation*}
\Lambda_q:= 1 +{\$\mathring{R}_{q}\$_{L^1,12 r_0}^{12} }+ \$v_{q}\$_{C^1_{t,x},10r_0}^{10}.
\end{equation*} 

Let $\ell_1=1/2$. 
For $q\geq 2$, let $\ell_{q}$ be such that 
\begin{equation*}
\ell_q^{\g}\ \Lambda_{q-1}:=  \eps_{q}^{N_0}.
\end{equation*}

For $q \geq 1$, we define $M_{q}(n,p)$  by 
\begin{equation*}
\begin{aligned}
M_{q}(n, p) := \ell_{q}^{-1} +\$\mathring{R}_{q}\$_{L^1, np}^{n}  + \$v_{q}\$_{C_{t,x}^1, np}^{n}.
\end{aligned}
\end{equation*}

For $q \geq 1$, we define 
\begin{equation*}
  \lambda_{q}:= {(\lfloor \ell_q^{-3p_0/(2-p_0)}+\ell_q^{-40}\rfloor+1)^8,}
 \end{equation*}
 where $\lfloor x\rfloor$ is the largest integer no greater than x.
 
Let $\zeta_1=2$. 
For $q\geq 2$, we define $$\zeta_q=\ell_q^{-1}.$$
\end{definition}

\begin{remark}
{In particular, from the fact that $\a_0\geq 12r_0$, there holds
\begin{equation}\label{eq:M:bound:0}
\Lambda_{q} \leq M_{q}(\a_0,m),\ \ \ m\geq 1,\ q\in\mN.
\end{equation}
and from the definition of $\lambda_q$, we have
\begin{align}
    \ell_{q}^{-19}\lambda_{q}^{-1/16}+\ell_{q}^{-23} \lambda_{q}^{-\frac98(\frac2{p_0}-1)} \lesssim \ell_{q},\ \ \ q \geq 1.\label{para:l}
\end{align}}
\end{remark}

\begin{definition}
For $q \geq 1$, we define  
\begin{equation}\label{def:bZ:q}
\bar{Z}_{q+1}:= {\sup_{t \in \mathbb{R}}\E[\|{z}_{q}- \bar{z}_{q} \|_{C_tL^2}^{2r_0}]}+ \sup_{t \in \mathbb{R}}\E[\|\bar{z}_{q+1}- \bar{z}_{q} \|_{C_tL^2}^{2r_0}] + \sup_{t \in \mathbb{R}}\ell_{q+1}^{2\g r_0}\ \E[ \|\bar{z}_{q} \|_{C^{\g}_{t}H^{\g}}^{2r_0}],
\end{equation}
where we recall that $\bar{z}_q = \Pi_{\zeta_q}z_q$.
\end{definition}

\begin{definition}
For $q \geq 1$, we define  
\begin{equation}\label{eq:def:theta}
\theta_{q}(t):=\frac1{3\cdot (2\pi)^3} \Big[ e(t)(1-\eps_{q-2})-\E\|v_{q}(t)+\bar{z}_q(t)\|_{L^2}^2\Big],
\end{equation}
and
\begin{equation*}
\delta E_{q}(t): = \big| e(t)(1-\eps_{q-3})-\E\|v_{q}(t)+\bar{z}_{q}(t)\|_{L^2}^2 \big|,
\end{equation*}
 where $\theta_q$ will appear as part of the amplitude in $w_{q+1}^{(p)}$ to meet the requirement of the energy profile for the solution (see \eqref{principle}).
\end{definition}

At each step $q$, the term $\Lambda_q$ is introduced to package all the terms concerning the previous steps in the building blocks $w_{q+1}$. The high frequency parameters $r_{\perp}$, $r_{||}$, $\lambda$ and $\mu$ which characterize the nature of intermittency are designed to 'cancel' the term $\Lambda_q$ so that the kinetic norm of the perturbation is small (see the proof of Lemma \ref{l:iterative:R:p0}). Moreover, we reveal the nature that $\Lambda_q$ is controlled by some iterative power which appears in the traditional presentation of convex integration scheme (see Lemma \ref{l:M:bound}). 

Let us introduce the following assumption on the iterations.
\begin{assumption}\label{a:iterative}
We assume  the following  statements hold:
\begin{enumerate}
\item 
There exists  some universal constant $c_R$, independent of $q$, such that
\begin{equation}\label{eq:R:p0}
\$\mathring{R}_q\$_{L^1,r_0}\leq  c_R \eps_{q-1}^{}.
\end{equation}
\item  There exists  some universal constant $c_v$, independent of $q$, such that
\begin{equation}\label{eq:v:p0}
{\$v_{q}\$_{L^2,2r_0}\leq  c_v\sum_{i=1}^q\eps_{i-4}^{1/2}}, \ \ \$v_{q}-v_{q-1}\$_{L^2,2r_0}\leq  c_v\eps_{q-4}^{1/2},\ \ \$v_q-v_{q-1}\$_{L^{p_0},2r_0}  {\leq c_v}  \ell_{q}^{\g} \Lambda_{q-1}^{1/2}.
\end{equation}
\item 
For any $m\geq 1$, there holds
\begin{equation}\label{eq:M:bound}
\Lambda_{q} \leq M_{q}(\a_0,m) \leq \eps_{q+1}^{-1}.
\end{equation}
\item 
For all $ t \geq0$, there exists $c_{\theta}$ and $c_E$, independent of $q$, such that
\begin{equation}\label{eq:energy}
0 < \theta_{q}(t) \leq   c_{\theta}\eps_{q-3},\ \ \   \d E_{q}(t) \leq c_E \eps_{q-2}^{}.
\end{equation}
\end{enumerate}
\end{assumption}

\begin{remark}
If a constant is independent of the iteration step $q$, this constant is universal. 
 We remark here that all the implicit constants in '$\lesssim$' are universal and independent of the constants $c_R$, $c_v$, $c_{\theta}$ and $c_E$. The constants $c_R$, $c_v$, $c_{\theta}$ and $c_E$ are chosen to be large enough according to the initial condition above and  the proof of Lemma \ref{l:iterative:v}, Lemma \ref{l:iterative:R:p0} and Lemma \ref{l:theta:E}, where we shall see that $c_v$, $c_{\theta}$ and $c_E$ depends on $c_R$ while $c_R$ is universal and independent of $c_v$, $c_{\theta}$ and $c_E$. 
\end{remark}

 We start our iteration  from $v_0= v_1\equiv0$.

\begin{proposition}
Let $v_0= v_1\equiv0$ on $\mathbb{R}$. Then  there exists $\mathring{R}_1$ such that Assumption \ref{a:iterative} holds for $q=1$.
\end{proposition}

\begin{proof}
Since $v_0= v_1\equiv0$ on $\mathbb{R}$, we obtain  \eqref{eq:v:p0} for $q=1$. Then we obtain $z_1$ by solving \eqref{eq:stochastic:1} and  have that
\begin{equation*}
 \bar{z}_{1}=\Pi_{\zeta_{1}}(z_{1}),\ \ \mathring{R}_1=\bar{z}_{1}\mathring\otimes \bar{z}_{1}-c\ \mathcal{R}\bar{z}_{1}.
\end{equation*}
 It is easy to see that $(v_1,\mathring{R}_1)$ solves the system  \eqref{induction ps}.
By \eqref{eR}, there holds
\begin{equation*}
\$\mathring{R}_1\$_{L^1,m} \leq \$\bar{z}_1\$_{L^2,2m}^2 + c \$\bar{z}_1\$_{L^1,m}\ ,\ m \in \mathbb{N}.
\end{equation*}
By Lemma \ref{l:truncate},  $\sup_{\omega \in \Omega}\|\bar{z}_1\|_{L^{\infty}} <\infty$. Hence,  by choosing $c_R>0$ large enough we obtain \eqref{eq:R:p0} for $q=1$. Moreover, we can choose $\Xi$ to be large enough such that
\begin{equation}\label{eq:Xi}
M_1(\a_0,m) \leq  {\ell_{1}^{-1} +\$\mathring{R}_{1}\$_{L^1, \alpha_0m}^{\alpha_0}}\leq \Xi^{\a_0}= \eps_{2}^{-1},\ \ \forall m \in\mathbb{N},
\end{equation}
 which {together with \eqref{eq:M:bound:0}} implies that \eqref{eq:M:bound} holds for $q=1$.

By letting $\underline{e}$ large enough and applying Lemma \ref{l:z:holder}, we have
 \begin{equation*}
\begin{aligned}
 \bar{e}   \geq& e(t)(1-\eps_{-1}) - \E\|v_1(t)+\bar{z}_1(t)\|_{L^2}^{2}\geq e(t)(1-\eps_{-2}) - (\E\|\bar{z}_1(t)\|_{L^2}^{2r_0})^{1/r_0}>0.
\end{aligned}
\end{equation*}
which implies that
\begin{align*}
    0<\theta_1(t)\leq \frac{\bar{e}}{3\cdot (2\pi)^3}    \leq c_{\theta}\eps_{-2},\ \   \d E_{1}(t)\leq\bar{e} \leq c_E \eps_{-1},
\end{align*}
by choosing sufficiently large     $c_{\theta}$ and $c_E$.
Therefore, Assumption \ref{a:iterative} holds when $q=1$. 
\end{proof}

The following Proposition is our main iteration:
{
\begin{proposition}\label{l:main}
 Assume that  Assumption \ref{a:iterative} holds up to step $q$, then it  holds  at step $q+1$. Moreover, there exists some universal constant $c_Z$ such that
\begin{equation}\label{eq:Z:q}
 \bar{Z}_{q+1} \leq c_Z \ell_{q}^{2\g r_0}\Lambda_{q-1}^{r_0}.
\end{equation} 
\end{proposition}
}
The proof of this proposition is provided in Section \ref{sec:iteration}. Assuming the validity of Proposition \ref{l:main}, we can prove Theorem \ref{thm:main:1} and Theorem \ref{thm:main:2}.

\subsection{Proof of Theorem \ref{thm:main:1}}
By interpolation, for any $\vartheta\in (0,1)$, we have
\begin{equation*}
\sum_{q\geq0}\$v_{q+1}-v_q\$_{H^{\vartheta},2r_0}\lesssim \sum_{q\geq0}\$v_{q+1}-v_q\$_{L^2,2r_0}^{1-\vartheta}\$v_{q+1}-v_q\$_{\C^1,2r_0}^{\vartheta}.
\end{equation*}
By  \eqref{eq:M:bound}, there holds
\begin{equation*}
\$v_{q+1}-v_q\$_{\C^1,2r_0} \lesssim \eps_{q+2}^{-1} + \eps_{q+1}^{-1} \lesssim \eps_{q+2}^{-1}.
\end{equation*}
By the second statement in \eqref{eq:v:p0}, we have 
\begin{equation*}
\$v_{q+1}-v_q\$_{L^2,2r_0}  \lesssim \eps_{q-3}^{1/2},
\end{equation*}
which implies that
\begin{equation*}
\begin{aligned}
\sum_{q\geq0}\$v_{q+1}- v_q\$_{H^{\vartheta},2r_0} \lesssim  \sum_{q\geq0} \eps_{q-3}^{(1-\vartheta)/2}\ \eps_{q+2}^{-\vartheta}.
\end{aligned}
\end{equation*}
Moreover, there holds $\eps_{q-3}=\eps_{q+2}^{M_0^{-5}}$.  Therefore, we obtain that for $\vartheta >0$ small enough
\begin{equation}\label{eq:dv:converge}
\sum_{q\geq0}\$v_{q+1}-v_q\$_{H^{\vartheta},2r_0} <\infty.
\end{equation}
Similarly we could change $H^{\vartheta}$ to $C^\vartheta_tL^2$. Hence, we have proved \eqref{est:u1}. 

As a consequence, a limit $v=\lim_{q\rightarrow\infty}v_q$ exists and lies in $L^{2r_0}(\Omega ; C([0,+\infty);H^{\vartheta})\cap C^\vartheta([0,+\infty);L^2))$.  Since $v_q$ is $(\mathcal{F}_t)_{t\ge0}$-adapted for every $q\in\mathbb{N}_{0}$, the limit $v$ is $(\mathcal{F}_t)_{t\geq0}$-adapted as well.
	
Since $v_q$ and $z_q$  satisfy the stochastic partial differential equation \eqref{induction ps} and \eqref{eq:stochastic:1},  we obtain that
\begin{equation}\label{eq:sol:q}
\begin{aligned}
\langle u_q(t),&\psi \rangle -\langle u_q(s),\psi \rangle + \int_{s}^{t}\langle \div((v_q+\bar{z}_q)\otimes (v_q+\bar{z}_q)),\psi \rangle \dif r\ = \\
&c\int_s^t \langle  \bar{z}_q-z_q ,\psi \rangle dr+\int_{s}^{t} \langle \Delta u_q, \psi\rangle \,\dif r +\langle \int_s^t G(v_{q-1}+z_q) \dif W_r,\psi\rangle+ \int_s^t \langle \div \mathring{R}_{q}, \psi \rangle,
\end{aligned}
\end{equation}
where $u_q=v_q + z_q$, $\psi$ is a divergence free smooth vector  field.
	
It follows from \eqref{eq:R:p0} that $\lim_{q\rightarrow\infty}\mathring{R}_q=0$ in $L^{r_0}(\Omega,C([0,\infty);L^1))$.  Moreover, it follows from  \eqref{zq:moment} and  \eqref{eq:v:p0} that 
\begin{equation*}
\sum_{q \geq 0}\${z}_{q+1}- {z}_q\$_{C_t^{\g}H^{\g},2r_0} \lesssim \sum_{q \geq -1} \$v_{q+1}-v_q\$_{L^2,2r_0}< \infty,
\end{equation*}
which implies that a limit $z= \lim_{q \to \infty} z_q$ exists in $L^{2r_0}(\Omega;C^{\g}([0, \infty); H^{\g})$ and the limit is also  $(\mathcal{F}_t)_{t\geq0}$-adapted.  {Moreover, by \eqref{eq:Z:q}, we have $z= \lim_{q \to \infty} \bar{z}_q$ in $L^{2r_0}(\Omega;C^{\g}([0, \infty); H^{\g})$.}
 
 By \eqref{eq:G:lip}, we have
\begin{equation}\label{eq:G:converge}
\begin{aligned}
\E[|\langle\int_s^t &\big(G(v_{q-1}+z_q) -G(v_q+z_{q+1})\big) \dif W_r,\psi\rangle|^2]\notag\\
& \lesssim  \int_s^t \E[\|G(v_{q-1}+z_q) -G(v_q+z_{q+1})  \|_{\mathcal{L}_0^2(U;H^{-2\d_0})}^2] \dif r \|\psi\|_{H^{2\d_0}}^2 \notag\\
&\lesssim_{\psi} \int_s^t \E[\|v_{q-1}-v_q\|_{L^2}^2+ \|z_{q}-z_{q+1}\|_{L^{2}}^2]\dif r .
\end{aligned}
\end{equation}
Therefore,  we can pass the limit in \eqref{eq:sol:q} and obtain that $u:=v+z$ is an $(\mathcal{F}_{t})_{t\geq 0}$-adapted analytically weak solution to \eqref{1}. Finally, \eqref{eq:K1} follows from \eqref{eq:energy}.
 
 \subsection{Proof of Theorem \ref{thm:main:2}}
Let $u$ be a probabilistically strong and analytically weak solution of the stochastic Navier-Stokes equation \eqref{1}   constructed in Theorem \ref{thm:main:1} with some constant energy profile,
\begin{equation*}
\mathbb{E}\|u\|_{L^2}^2=K_0.
\end{equation*}
Even though the solution $u$, as long as the Wiener process $W$ are defined on $[0, +\infty)$, we extend the couple $(u(t),W(t))$ from $[0, +\infty)$ to $\mathbb{R}$ by setting $(u(t),W(t))=(u(0),W(0))$, for $t<0$. We define $\nu_{\tau,T}$,  a measure on $\mathcal{T}$, by
\begin{align*}
\nu_{\tau,T}=\frac1T\int_0^T\mathcal{L}[S_{t+\tau}(u,W)]\dif t,\qquad T\geq 0.
\end{align*}
For $\vartheta>0$ small enough, by Theorem \ref{thm:main:1}, we have
\begin{equation}\label{eq:pf:thm1.7:1}
\begin{aligned}
\sup_{s\geq 0}\E\sup_{t\in[0,N]}\|u(t+s)\|^{2r_0}_{H^\vartheta}&\lesssim\sup_{s\geq 0}\sum_{i=0}^{N-1}\E\sup_{t\in[i,i+1]}\|u(t+s)\|_{H^\vartheta}^{2r_0}\\
& \lesssim N\sup_{s\geq0}\E\sup_{t\in[0,1]}\|u(t+s{+i})\|^{2r_0}_{H^\vartheta}\lesssim N,
\end{aligned}
\end{equation}
and similarly
\begin{equation}\label{eq:pf:thm1.7:2}
\begin{aligned}
\sup_{s\geq 0}\E\|u(\cdot+s)\|^{2r_0}_{C^\vartheta([0,N], L^2)}\lesssim  N.
\end{aligned}
\end{equation}
For $R_N:=MN^2$, $M,N\in\mathbb{N}$, we note that the set
\begin{align*}
K:=\bigcap_{N=1}^\infty\bigg\{(f,g);\,&\|f\|_{C^\vartheta_{[-N,N]} L^2}+\sup_{t\in[-N,N]}\|f(t)\|_{H^\vartheta} + \|g\|_{C^\vartheta_{[-N,N]} H^{\vartheta}}\leq R_N\bigg\}
\end{align*}
is relatively compact in $C(\mathbb{R};L^2_{\sigma})$. Moreover, there holds
\begin{equation*}
\begin{aligned}
& \(\frac{1}{T} \int_0^T \mathcal{L}[S_{t+\tau}(u,W)]\dif t\)(K^c) \leq \frac{1}{T}\int_0^T  \mathcal{L}[S_{t+\tau}(u,W)](K^c)\dif t\\
\leq&\frac{1}{T}\int_0^T \sum_{N=1}^{\infty}  \mathcal{L}[S_{t+\tau}(u,W)]\((f,g);\,\|f\|_{C^\vartheta_{[-N,N]} L^2}+\sup_{t\in[-N,N]}\|f(t)\|_{H^\vartheta} +\|g\|_{C^\vartheta_{[-N,N]} H^{\vartheta}} > R_N\)\dif t\\
\lesssim& \frac{1}{T}\int_0^T \sum_{N=1}^{\infty} {R_N}^{-2r_0}\(\sup_{t \geq 0} \E[\|u\|_{C^\vartheta_{[t,N+t]} L^2}^{2r_0}+\sup_{t\in[t,N+t]}\|u(t)\|_{H^\vartheta}^{2r_0}] + \E[\|W\|_{C^\vartheta_{[t,N+t]} H^{\vartheta}}^{2r_0}]  \)\dif t.
\end{aligned}
\end{equation*}
By  \eqref{eq:pf:thm1.7:1} and  \eqref{eq:pf:thm1.7:2}, we have
\begin{equation*}
\(\frac{1}{T} \int_0^T \mathcal{L}[S_{t+\tau}(u)]\dif t\)(K^c) \lesssim  \sum_{N=1}^{\infty} {R_N^{-2r_0}} N.
\end{equation*}
For any $\eps>0$, $M$ can be chosen large enough so that $\nu_{\tau,T}(K^c) < \eps$. This implies the tightness of the family of measures $\{\nu_{\tau,T-\tau}\}_{T >0}$ for every fixed $\tau>0$. Now, we can process as in \cite[Theorem 4.2]{HZZ23} to prove the existence of ergodic stationary solution. Finally, \eqref{eq:s56} follows directly from  \eqref{eq:pf:thm1.7:1} and  \eqref{eq:pf:thm1.7:2}.

\section{Proof of main iteration: Proposition \ref{l:main}}
\label{sec:iteration}
{
In this section, we focus on the proof of our main iteration in the stochastic convex integration method, namely, Proposition \ref{l:main}.   Firstly, Section \ref{sec:esti:noise} is devoted to the estimate of the noise term. Then, we consider the mollification step to avoid the loss of derivatives in Section \ref{sec:mollification}. In Section \ref{sec:perturbation}, we introduce the perturbation by applying the intermittent jets in Appendix \ref{s:B}. Section \ref{sec:est:wq+1} is devoted to provide some necessary estimate. In Section \ref{sec:w}, the new stress term at the step $q+1$ is introduced and the iterative estimate for this stress term is included. Finally, we show the estimate of energy  in Section \ref{sec:est:energy}.

\subsection{Estimate of stochastic term}\label{sec:esti:noise}
 Assume Assumption \ref{a:iterative} holds up to step $q$,  we solve \eqref{eq:stochastic:1} for a given $v_q$ to derive the new noise term  $z_{q+1}$ and the corresponding $\bar{z}_{q+1}$. This section focuses on estimating the difference between these  noise terms as defined in \eqref{def:bZ:q}.
\begin{lemma}\label{l:barZ:q}
Assume Assumption \ref{a:iterative} holds up to step $q$. There exists some universal constant $c_Z$ such that
\begin{equation*}
 \bar{Z}_{q+1} \leq c_Z \ell_{q}^{2\g r_0}\Lambda_{q-1}^{r_0}.
\end{equation*}
\end{lemma}
\begin{proof}
By Lemma \ref{l:truncate}, we have
\begin{equation*}
\begin{aligned}
&\sup_{t \in \mathbb{R}}\E[\|\bar{z}_{{q+1}}- \bar{z}_{q} \|_{C_tL^2}^{2r_0}] \\
\lesssim&  \sup_{t \in \mathbb{R}}\E[\|\Pi_{\zeta_{q+1}}{z}_{q+1}- \Pi_{\zeta_{q+1}}{z}_{q} \|_{C_tL^2}^{2r_0}] + \sup_{t \in \mathbb{R}}\E[\|\Pi_{\zeta_{q+1}}{z}_{q}- \Pi_{\zeta_{q}}{z}_{q} \|_{C_tL^2}^{2r_0}]  \\
\lesssim&  \sup_{t \in \mathbb{R}}\E[\|{z}_{q+1}- {z}_{q} \|_{C_tL^2}^{2r_0}] +  \sup_{t \in \mathbb{R}} \E[\|{z}_{q}- \Pi_{\zeta_{q}}{z}_{q} \|_{C_tL^2}^{2r_0}] .
\end{aligned}
\end{equation*}
{By \eqref{zq:moment} and \eqref{eq:v:p0} we have
\begin{equation}
\begin{aligned}
& \sup_{t \in \mathbb{R}}\E[\|{z}_{q+1}- {z}_{q} \|_{C_tL^2}^{2r_0}] 
\leq C_{c,r_0,\g,\d_0,L} \$v_q-v_{q-1}\$_{L^{p_0},2r_0}^{2r_0}\lesssim \ell_{q}^{2\g r_0} \Lambda_{q-1}^{r_0}.
\end{aligned}\label{eq:Z:q:1}
\end{equation}}

Now we need to estimate the term $\sup_{t \in \mathbb{R}}\E[\|{z}_{q}- \Pi_{\zeta_q}{z}_{q} \|_{C_tL^2}^{2r_0}]$, we have
\begin{equation*}
\begin{aligned}
\E[\sup_{s\in[t,t+1]}\|{z}_{q}(s)-& \Pi_{\zeta_q}{z}_{q}(s) \|_{L^2}^{2r_0} ] = \E[\sup_{s\in[t,t+1]}\|{z}_{q}(s)- \Pi_{\zeta_q}{z}_{q}(s) \|_{L^2}^{2r_0} \chi_{\{\sup_{s \in [t,t+1]}\|z_{q}(s)\|_{L^2}^2\leq {\zeta_q^2}\}}] \\
&+  \E[\sup_{s\in[t,t+1]}\|{z}_{q}(s)- \Pi_{\zeta_q}{z}_{q}(s) \|_{L^2}^{2r_0} \chi_{\{\sup_{s \in [t,t+1]}\|z_{q}(s)\|_{L^2}^2> {\zeta_q^2}\}}] .
\end{aligned}
\end{equation*}
For any $f= \sum_{k \in \mathbb{Z}^3} f_k e_k \in L^2$, with $\|f\|_{L^2}^2=\sum_{k \in \mathbb{Z}^3}f_k^2 \leq \zeta_q^2$, we have $f_k^2{\leq \zeta_q^2}$ for all $k \in \mathbb{Z}^3$. {Together with the definition of the projection $\Pi_{\zeta}$ in Section \ref{sec:operator}, we imply} that
\[
f- \Pi_{\zeta_q}f =\sum_{|k|^2 \leq \zeta_q} (f_k-f_{k,\zeta_q})e_k + \sum_{|k|^2>\zeta_q} f_k e_k = \sum_{|k|^2>\zeta_q} f_k e_k.
\]
Hence,  $\|z_{q}(s)\|^2_{L^2} \leq  \sup_{s \in [t,t+1]}\|z_{q}(s)\|_{L^2}^2\leq \zeta_q^2$, $\forall s \in [t,t+1]$, it gives us that
\begin{equation*}
\begin{aligned}
\|{z}_{q}(s)- \Pi_{\zeta_q}{z}_{q}(s) \|_{L^2}^{2r_0} &\chi_{\{\sup_{s \in [t,t+1]}\|z_{q}(s)\|_{L^2}^2\leq {\zeta_q^2}\}} \\
&= \(\sum_{|k|^2 >\zeta_q} \langle z_{q}(s), e_k \rangle^2\)^{r_0} \chi_{\{\sup_{s \in [t,t+1]}\|z_{q}(s)\|_{L^2}^2\leq {\zeta_q^2}\}},
\end{aligned}
\end{equation*}
which implies that
\begin{equation*}
\begin{aligned}
&\E[\sup_{s\in[t,t+1]}\|{z}_{q}(s)- \Pi_{\zeta_q}{z}_{q}(s) \|_{L^2}^{2r_0} \chi_{\{\sup_{s \in [t,t+1]}\|z_{q}(s)\|_{L^2}^2\leq \zeta_q^2\}}] \\
\leq& \zeta_q^{-2\g r_0}\ \E[\sup_{s\in[t,t+1]} \(\sum_{|k|^2 >\zeta_q} |k|^{4\g} \langle z_{q}(s), e_k \rangle^2\)^{r_0} ]\leq \ell_q^{2\g r_0}\ \E[\sup_{s \in [t,t+1]} \|z_{q}(s)\|_{H^{2\g}}^{2r_0}].
\end{aligned}
\end{equation*}
 On the other hand, by the H\"older inequality, we have 
\begin{equation*}
\begin{aligned}
 & \E[\sup_{s\in[t,t+1]}\|{z}_{q}(s)- \Pi_{\zeta_q}{z}_{q}(s) \|_{L^2}^{2r_0} \chi_{\{\sup_{s \in [t,t+1]}\|z_{q}(s)\|_{L^2}^2> \zeta_q^2\}}] \\
  \lesssim &  \E[\sup_{s\in[t,t+1]}\|{z}_{q}(s) \|_{L^2}^{2r_0} \chi_{\{\sup_{s \in [t,t+1]}\|z_{q}(s)\|_{L^2}^2> \zeta_q^2\}}] \\
  \lesssim&  \E[\sup_{s\in[t,t+1]}\|{z}_{q}(s) \|_{L^2}^{4r_0}]^{1/2} \E[\chi_{\{\sup_{s \in [t,t+1]}\|z_{q}(s)\|_{L^2}^2> \zeta_q^2\}}]^{1/2}.
  \end{aligned}
  \end{equation*}
  By Chebyshev's inequality, we have 
  \begin{equation*}
  \begin{aligned}
 \E[\sup_{s\in[t,t+1]}\|{z}_{q}(s)- \Pi_{\zeta_q}{z}_{q}(s) \|_{L^2}^{2r_0} \chi_{\{\sup_{s \in [t,t+1]}\|z_{q}(s)\|_{L^2}^2> \zeta_q^2\}}]  \lesssim& \zeta_q^{-2r_0} \E[\sup_{s\in[t,t+1]}\|{z}_{q}(s) \|_{L^2}^{4r_0}].
  \end{aligned}
\end{equation*}
 Therefore, we obtain that
\begin{equation*}
\begin{aligned}
&\sup_{t \in \mathbb{R}} \E[\|z_{q}-\Pi_{\zeta_q} z_{q}\|_{C_t L^2}^{2r_0}] \lesssim \ell_{q}^{2\g r_0}  \big(\sup_{t \in \mathbb{R}}&\E[\sup_{s \in [t,t+1]} \|z_{q}(s)\|_{H^{\g}}^{2r_0}]+ \sup_{t \in \mathbb{R}}\E[\sup_{s\in[t,t+1]}\|{z}_{q}(s) \|_{L^2}^{4r_0}]\big).
\end{aligned}
\end{equation*}
Since \eqref{def:gamma:p} holds, by \eqref{zq:moment} we obtain that
\begin{equation}\label{eq:Z:q:2}
\sup_{t \in \mathbb{R}}\E[\|{z}_{q}- \Pi_{\zeta_q}{z}_{q} \|_{C_tL^2}^{2r_0}] \lesssim \ \ell_{q}^{2\g r_0}\Lambda_{q-1}^{r_0},
\end{equation}
and
\begin{equation}\label{eq:Z:q:3}
\sup_{t \in \mathbb{R}}\E[ \ell_{q+1}^{2\gamma  r_0}\|z_{q} \|_{C^{\gamma}_{t}H^{\gamma}}^{2r_0}] \lesssim \ell_{q+1}^{2\gamma r_0}( 1 + \$v_{q-1}\$_{L^{\g_0},2r_0}^{2r_0}) \lesssim \ \ell_{q}^{2\gamma r_0}\Lambda_{q-1}^{r_0}.
\end{equation}
Combining \eqref{eq:Z:q:1}, \eqref{eq:Z:q:2} and \eqref{eq:Z:q:3}, there exists some universal constant $c_Z$ such that
\begin{equation*}
 \bar{Z}_{q+1} \leq  c_Z\ell_{q}^{2\g r_0}\Lambda_{q-1}^{r_0} .
 \end{equation*}
 \end{proof}

\subsection{Mollification}\label{sec:mollification}

 To avoid the loss of derivatives in the construction, we introduce the mollification parameter 
\begin{equation}\label{eq:parameters:l}
\ell:=\ell_{q+1}\in (0,1/2),
\end{equation}
 and replace $v_q$ by a mollified velocity field $v_{q,\ell}$ on $\mathbb{R}$.
 
 Let $\rho_\ell=\frac{1}{\ell^3}\rho(\frac{\cdot}{\ell})$ be a family of standard mollifiers on $\mathbb{R}^3$, and let $\varphi_\ell=\frac{1}{\ell}\varphi(\frac{\cdot}{\ell})$ be a family of  standard mollifiers with support in $(0,1)$. 
In the following, we shall adopt the notation  by
\begin{equation}\label{def:moll}
\bar{z}_{q,\ell}:= (\bar{z}_{q}*_x\rho_\ell)*_t\varphi_\ell,\ 
v_{q,\ell}:= (v_{q}*_x\rho_\ell)*_t\varphi_\ell,\ 
\mathring{R}_{q,\ell}:= (\mathring{R}_{q}*_x\rho_\ell)*_t\varphi_\ell.
\end{equation}
Since the mollifier $\varphi_\ell$ is supported on $(0,1)$, it is easy to see that $\bar{z}_{q,\ell}$ is $(\mathcal{F}_t)_{t\geq0}$-adapted and so are $v_{q,\ell}$ and $\mathring{R}_{q,\ell}$.

 Since \eqref{induction ps} holds at step $q$ on $\mathbb{R}$, it follows that $(v_{q,\ell_{}},\mathring{R}_{q,\ell_{}})$ satisfies
\begin{equation}\label{eq:mollification}
\aligned
\partial_tv_{q,\ell_{}} - c\bar{z}_{q,\ell_{}}-\Delta v_{q,\ell_{}}+\div((v_{q,\ell_{}}+\bar{z}_{q,\ell_{}})\otimes (v_{q,\ell_{}}+\bar{z}_{q,\ell_{}}))+\nabla p_{q,\ell_{}}&=\div (\mathring{R}_{q,\ell_{}}+\mathring{R}^{(q+1)}_{\textrm{mol}}),
\\\div v_{q,\ell_{}}&=0,
\endaligned
\end{equation}
 where
\begin{equation*}
\mathring{R}_{\textrm{mol}}^{(q+1)}=(v_{q,\ell}+\bar{z}_{q,\ell})\mathring{\otimes}(v_{q,\ell}+\bar{z}_{q,\ell})-((v_{q}+\bar{z}_{q})\mathring{\otimes}(v_{q}+\bar{z}_{q}))_{\ell},
\end{equation*}
 and $p_{q,\ell}$ is the associated pressure term.

\subsection{Perturbation at step $q+1$}\label{sec:perturbation}
We define $v_{q+1}$ by adding a perturbation term $w_{q+1}$ to the mollified term $v_{q,\ell}$,
\begin{equation}\label{eq:v:q+1}
v_{q+1} := v_{q,\ell}+ w_{q+1},
\end{equation}
and the pair $(v_{q+1}, \mathring{R}_{q+1})$ is designed to solve the following system, 
\begin{equation*}
\aligned
\partial_tv_{q+1}-c \bar{z}_{q+1}-\Delta v_{q+1} +\div((v_{q+1}+\bar{z}_{q+1})\otimes (v_{q+1}+\bar{z}_{q+1}))+\nabla p_{q+1}&=\div \mathring{R}_{q+1},
\\
\div \ v_{q+1}&=0.
\endaligned
\end{equation*}
where we recall that $\bar{z}_{q+1}=\Pi_{\zeta_{q+1}}(z_{q+1})$,  $z_{q+1}$ is defined  analogous to \eqref{eq:stochastic:1} with $v_{q-1}$ replaced by $v_q$.  $p_{q+1}$ is the associated pressure term.

 To define the perturbation $w_{q+1}$, we first define  
\[
 \theta_{q,\ell}:=\theta_q*_t\varphi_{\ell_{q+1}},
\]
where we recall that  $\theta_{q}$ is defined  in \eqref{eq:def:theta}.
By \eqref{eq:energy} we know $\theta_q>0$. Moreover,  we write
\begin{equation}\label{eq:g:bound:1}
\bar{\theta}_q:=\sup_{t \geq0}|\theta_{q}(t)|. 
\end{equation}

Let $\rho$ be defined by
\begin{equation}\label{eq:rho3}
\rho_{q+1}:=2\sqrt{\ell_{q+1}^2+|\mathring{R}_{q,\ell}|^2}+\theta_{q,\ell}.
\end{equation}
We define the principal part $w_{q+1}^{(p)}$ of the perturbation $w_{q+1}$ as
\begin{equation}\label{principle}
w_{q+1}^{(p)}:=\sum_{\xi\in\Lambda} a_{(\xi)}W_{(\xi)},
\end{equation}
so that whose low frequency part cancels $ \div\mathring{R}_{q,\ell}$. Here the amplitude $a_{(\xi)}:= a_{(\xi),q+1}$ at step $q+1$ is defined by
\begin{equation*}
a_{(\xi)}(\omega,t,x):=\rho_{q+1}(\omega,t,x)^{1/2}\gamma_\xi\left(\Id-\frac{\mathring{R}_{q,\ell}(\omega,t,x)}{\rho_{q+1}(\omega,t,x)}\right),
\end{equation*}
where $\gamma_{\xi}$ is defined in Lemma \ref{geometric} and $W_{(\xi)}$ is defined in  \eqref{intermittent} with all the parameters depending on $\lambda_{q+1}$:
\begin{equation}\label{eq:parameters:q+1}
r_{\perp}:=\lambda_{q+1}^{-7/8},\ r_{\|}:=\lambda_{q+1}^{-1/2},\ \lambda:=\lambda_{q+1},\ \mu:= \lambda_{q+1}^{5/4}.
\end{equation}
As a direct consequence of  Lemma \ref{geometric},  the following important identity holds
\begin{equation}\label{can}
w_{q+1}^{(p)}\otimes w_{q+1}^{(p)}+\mathring{R}_{q,\ell}=\sum_{\xi\in \Lambda}a_{(\xi)}^2 \mathbb{P}_{\neq0}(W_{(\xi)}\otimes W_{(\xi)})+\rho_{q+1} \Id.
\end{equation}
 Here we used the notation $\mathbb{P}_{\neq0} f := f -\frac{1}{(2\pi)^3}\int_{\mathbb{T}^3} f\dif x$.
For more details, we refer readers to \cite[Section 3]{HZZ22}.

It is clear that $W_{(\xi)}$ is not divergence free, we define the incompressibility corrector by
\begin{equation*}
w_{q+1}^{(c)}:=\sum_{\xi\in \Lambda}\textrm{curl}(\nabla a_{(\xi)}\times V_{(\xi)})+\nabla a_{(\xi)}\times \textrm{curl}V_{(\xi)}+a_{(\xi)}W_{(\xi)}^{(c)},
\end{equation*}
with $W_{(\xi)}^{(c)}$ and $V_{(\xi)}$ being given in \eqref{corrector} and all the parameters are defined in \eqref{eq:parameters:q+1}.

By a direct computation we deduce that
\begin{equation*}
	w_{q+1}^{(p)}+w_{q+1}^{(c)}=\sum_{\xi\in\Lambda}\textrm{curl}\,\textrm{curl}(a_{(\xi)}V_{(\xi)}),
\end{equation*}
hence
\begin{equation*}\div(w_{q+1}^{(p)}+w_{q+1}^{(c)})=0.\end{equation*}

Finally, we introduce the temporal corrector $w_{q+1}^{(t)}$ by
\begin{equation*}
w_{q+1}^{(t)}:=-\frac{1}{\mu}\sum_{\xi\in \Lambda}\mathbb{P}\mathbb{P}_{\neq0}\left(a_{(\xi)}^2\phi_{(\xi)}^2\psi_{(\xi)}^2\xi\right).
\end{equation*}

We define the perturbation term $w_{q+1}$ by
$$
w_{q+1} := w_{q+1}^{(p)} +  w_{q+1}^{(c)} +  w_{q+1}^{(t)}.
$$

By Gauss-Green Theorem, $ w_{q+1}^{(p)} +  w_{q+1}^{(c)} $ is mean zero. By a direct computation (see \cite[(7.20)]{BV19}),  we obtain
\begin{equation*}
\aligned
&\partial_t w_{q+1}^{(t)}+\sum_{\xi\in\Lambda}\mathbb{P}_{\neq0}\left(a_{(\xi)}^2\div(W_{(\xi)}\otimes W_{(\xi)})\right)
\\
&\qquad= -\frac{1}{\mu}\sum_{\xi\in\Lambda}\mathbb{P}\mathbb{P}_{\neq0}\partial_t\left(a_{(\xi)}^2\phi_{(\xi)}^2\psi_{(\xi)}^2\xi\right)
+\frac{1}{\mu}\sum_{\xi\in\Lambda}\mathbb{P}_{\neq0}\left( a^2_{(\xi)}\partial_t(\phi^2_{(\xi)}\psi^2_{(\xi)}\xi)\right)
\\&\qquad= (\Id-\mathbb{P})\frac{1}{\mu}\sum_{\xi\in\Lambda}\mathbb{P}_{\neq0}\partial_t\left(a_{(\xi)}^2\phi_{(\xi)}^2\psi_{(\xi)}^2\xi\right)
-\frac{1}{\mu}\sum_{\xi\in\Lambda}\mathbb{P}_{\neq0}\left(\partial_t a^2_{(\xi)}(\phi^2_{(\xi)}\psi^2_{(\xi)}\xi)\right).
\endaligned
\end{equation*}
Due to the Helmholtz decomposition, the first term on the right hand side can be absorbed in the pressure term which implies that $ w_{q+1}^{(t)}$ is also mean zero. Therefore, $w_{q+1}$ is divergence free and mean free. 
Since the coefficients $a_{(\xi)}$ are $(\mathcal{F}_t)_{t\geq0}$-adapted, we deduce that
$w_{q+1}^{(p)}$, $w_{q+1}^{(c)}$ and $w_{q+1}^{(t)}$ are $(\mathcal{F}_t)_{t\geq0}$-adapted. Then $w_{q+1}$ and $v_{q+1}$  are also $(\mathcal{F}_t)_{t\geq0}$-adapted.

\subsection{Estimate of $w_{q+1}$}\label{sec:est:wq+1}
Now we will provide some estimates for the perturbation term $w_{q+1}$.  As implied in \eqref{eq:parameters:l} and \eqref{eq:parameters:q+1}, all the parameters below depends on the information up to the step $q+1$. According to the proof in  \cite[Appendix B]{HZZ22},  for $N\geq1$ there holds
\begin{equation}\label{estimate aN}
\begin{aligned}
\|a_{(\xi)}\|_{C^N_{{[t,t+1]},x}}& \lesssim\ell^{-7-6N}(\|\mathring{R}_q\|_{C_{{[t-\ell,t+1]}}L^1}+1)^{N+1},
\end{aligned}
\end{equation}
and
\begin{equation}\label{estimate aN0}
\begin{aligned}
\|a_{(\xi)}\|_{C^0_{{[t,t+1]},x}}& \lesssim \ell^{-2}(\|\mathring{R}_q\|_{C_{{[t-\ell,t+1]}}L^1}+1)^{1/2}.
\end{aligned}
\end{equation}

\subsubsection{Estimate of $\|w_{q+1}\|_{C_t L^p}$}
For a general $L^p$-norm we apply (\ref{bounds}), \eqref{estimate aN} and (\ref{estimate aN0}) to deduce for  $p\in(1,\infty)$
\begin{equation}\label{est:principle:Lp}
\aligned
\|w_{q+1}^{(p)}\|_{C_{{[t,t+1]}}L^p}&\lesssim \sum_{\xi\in \Lambda}\|a_{(\xi)}\|_{C^0_{{[t,t+1]},x}}\|W_{(\xi)}\|_{C_{{[t,t+1]}}L^p}\lesssim  \ell^{-2}(\|\mathring{R}_q\|_{C_{{[t-\ell,t+1]}}L^1}+1)^{1/2}r_\perp^{2/p-1}r_\|^{1/p-1/2},
\endaligned
\end{equation}
\begin{equation}\label{est:correction:Lp}
\aligned
\|w_{q+1}^{(c)}\|_{C_{{[t,t+1]}}L^p}&\lesssim\sum_{\xi\in \Lambda}\left(\|a_{(\xi)}\|_{C^0_{{[t,t+1]},x}}\|W_{(\xi)}^{(c)}\|_{C_{{[t,t+1]}}L^p}+\|a_{(\xi)}\|_{C^2_{{[t,t+1]},x}}\|V_{(\xi)}\|_{C_{{[t,t+1]}}W^{1,p}}\right)
\\&\lesssim \ell^{-19}(\|\mathring{R}_q\|_{C_{{[t-\ell,t+1]}}L^1}+1)^3r_\perp^{2/p-1}r_\|^{1/p-1/2}\left(r_\perp r_\|^{-1}+\lambda_{}^{-1}\right),
\endaligned
\end{equation}
and
\begin{equation}\label{est:temporal:Lp}
\aligned
\|w_{q+1}^{(t)}\|_{C_{{[t,t+1]}}L^p}&\lesssim \mu^{-1}\sum_{\xi\in\Lambda}\|a_{(\xi)}\|_{C^0_{{[t,t+1]},x}}^2\|\phi_{(\xi)}\|_{L^{2p}}^2\|\psi_{(\xi)}\|_{C_{{[t,t+1]}}L^{2p}}^2\\
&\lesssim\ell^{-4}(\|\mathring{R}_q\|_{C_{{[t-\ell,t+1]}}L^1}+1)r_\perp^{2/p-1}r_\|^{1/p-1/2}(\mu^{-1}r_\perp^{-1}r_\|^{-1/2}).
\endaligned
\end{equation}
Therefore, we obtain that {for $p>1$}
\begin{equation}\label{eq:w:L1}
\|w_{q+1}\|_{C_{{[t,t+1]}}L^p}\lesssim  \ell^{-19} (r_\perp^{2/p-1}r_\|^{1/p-1/2})(1+\|\mathring{R}_q\|^{3}_{C_{{[t-\ell,t+1]}}L^1} ) (r_{\perp}r_{\|}^{-1} + \lambda^{-1} + \mu^{-1} r_{\perp}^{-1} r_{\|}^{-1/2}+1).
\end{equation}

\subsubsection{Estimate of $\|w_{q+1}\|_{C_{{[t,t+1]}} L^2}$}
When $p=2$, we  need some more delicate estimate for the $\|\cdot\|_{C_tL^2}$ norm. By applying \cite[Lemma B.1]{CL20}, we obtain that
\begin{equation}\label{eq:principle:Lp}
\begin{aligned}
&\|w_{q+1}^{(p)}\|_{C_{{[t,t+1]}}L^2}\lesssim \sum_{\xi\in \Lambda}\|a_{(\xi)} W_{(\xi)}\|_{C_{{[t,t+1]}}L^2}\\
 \lesssim& \sum_{\xi\in \Lambda}\|a_{(\xi)}\|_{C_{{[t,t+1]}}L^2}\|W_{(\xi)}\|_{C_{{[t,t+1]}}L^2} + (r_{\perp} \lambda)^{-1/2}\|a_{(\xi)}\|_{C^1_{{{[t,t+1]},x}}}\|W_{(\xi)}\|_{{C_{{[t,t+1]}}}L^2}\\  
\lesssim& (\ell + \|\mathring{R}_{q}\|_{C_{{[t-\ell,t+1]}}L^1} +\bar{\theta}_q)^{1/2} + (r_{\perp} \lambda)^{-1/2}\ell^{-13}(\|\mathring{R}_q\|_{C_{{[t-\ell,t+1]}}L^1}+1)^{2},
\end{aligned}
\end{equation}
where $\bar{\theta}_q$ is defined in \eqref{eq:g:bound:1}.

\begin{lemma}\label{l:iterative:v}
Assume Assumption \ref{a:iterative} holds up to step $q$. There exists some universal constant  $c_v$ such that Assumption \ref{a:iterative}(2) holds at step $q+1$, i.e. 
\begin{equation}\label{eq:a:dv:q+1}
{\$v_{q+1}\$_{L^2,2r_0}\leq  c_v\sum_{i=1}^{q+1}\eps_{i-4}^{1/2}}, \ \  \$v_{q+1}-v_{q}\$_{L^2,2r_0} \leq   c_v\eps_{q-3}^{1/2},\ \ \$v_{q+1}-v_{q}\$_{L^{p_0},2r_0}  {\leq c_v}  \ell_{q+1}^{\g} \Lambda_{q}^{1/2}.
\end{equation}
\end{lemma}
\begin{proof}
By \eqref{para:l}, \eqref{eq:principle:Lp}, \eqref{est:correction:Lp} and \eqref{est:temporal:Lp},  there holds
\begin{equation}
\begin{aligned}
\$w_{q+1}\$_{L^2,2r_0}\lesssim& {\ell_{q+1}^{-19}\lambda_{q+1}^{-1/16
}} (1+ \$\mathring{R}_{q}\$_{L^1,6r_0}^{3}) + (\ell_{q+1} + \$\mathring{R}_{q}\$_{L^1,r_0} +\bar{\theta}_{q})^{1/2}\\
\lesssim& {\ell_{q+1}} (1+ \$\mathring{R}_{q}\$_{L^1,6r_0}^{3}) + (\ell_{q+1} + \$\mathring{R}_{q}\$_{L^1,r_0} +\bar{\theta}_{q})^{1/2}.\label{bd:wq+1l2r}
\end{aligned}
\end{equation}
By Assumption \ref{a:iterative} and Lemma \ref{l:theta:E}, there holds
\begin{equation}
\begin{aligned}\label{bd:wq+1l2}
\$w_{q+1}\$_{L^2,2r_0}\lesssim&\ \ell_{q+1} \Lambda_{q} + \ell_{q+1}^{1/2}+  \eps_{q-1}^{1/2} +  c_{\theta}\eps_{q-3}^{1/2}.
\end{aligned}
\end{equation}
Moreover, we have
\begin{equation*}
\begin{aligned}
 \$ v_{q,\ell}- v_q \$_{L^2,2r_0} \lesssim \ell_{q+1} \$v_q\$_{\C^1,2r_0} \lesssim \ell_{q+1} \Lambda_{q}^{1/2}{\lesssim \eps_{q+1}}.
\end{aligned}
\end{equation*}
Hence, there exists some constant $c_v$ such that
\begin{equation*}
\$v_{q+1}-v_{q}\$_{L^2,2r_0} \leq \$w_{q+1}\$_{L^2,2r_0} + \$ v_{q,\ell}- v_q \$_{L^2,2r_0} \leq c_v \eps_{q-3}^{1/2}.
\end{equation*}
{We obtain the second term. Then the first term follows immediately.}

{By \eqref{para:l} and \eqref{eq:w:L1}  we have
 \begin{align}\label{bd:wq+p0}
     \$w_{q+1}\$_{L^{p_0},2r_0}\lesssim \ell_{q+1}^{-19} \lambda_{q+1}^{-\frac98(\frac2{p_0}-1)} \Lambda_{q}^{1/2}\lesssim \ell_{q+1}\Lambda_{q}^{1/2},
 \end{align}
 which implies that, {there exists some constant, still denoted by $c_v$, such that
\begin{equation*}
\$v_{q+1}-v_{q}\$_{L^{p_0},2r_0} \leq \$w_{q+1}\$_{L^{p_0},2r_0} + \$ v_{q,\ell}- v_q \$_{L^{p_0},2r_0} \leq c_v \ell_{q+1}^{\g} \Lambda_{q}^{1/2}.
\end{equation*}}
}
\end{proof}

\subsubsection{Estimate of $\|w_{q+1}\|_{C_{{[t,t+1]}}W^{1,p}}$}
For any $p\geq1$, there holds
\begin{equation}
\aligned
&\|w_{q+1}^{(p)}+w_{q+1}^{(c)}\|_{C_{{[t,t+1]}}W^{1,p}}\leq\sum_{\xi\in\Lambda}\|\textrm{curl\,}\textrm{curl}(a_{(\xi)}V_{(\xi)})\|_{C_{{[t,t+1]}}W^{1,p}}\\
&\quad\quad\lesssim \|a_{(\xi)}\|_{C_{{[t,t+1]},x}^3} \|V_{(\xi)}\|_{C_{{[t,t+1]}} L^p} + \|a_{(\xi)}\|_{C_{t,x}^2} \|V_{(\xi)}\|_{C_t W^{1,p}} \notag\\ 
&\quad\quad+  \|a_{(\xi)}\|_{C_{{[t,t+1]},x}^1} \|V_{(\xi)}\|_{C_{{[t,t+1]}} W^{2,p}} +  \|a_{(\xi)}\|_{C_{{[t,t+1]},x}^0} \|V_{(\xi)}\|_{C_{{[t,t+1]}} W^{3,p}}.
\endaligned
\end{equation}
By (\ref{bounds}), \eqref{estimate aN} and (\ref{estimate aN0}), we have
\begin{equation}\label{eq:w:pc:w1p}
\aligned
\|w_{q+1}^{(p)}+w_{q+1}^{(c)}\|_{C_{{[t,t+1]}}W^{1,p}}\lesssim& \ell^{-25} (1+ \|\mathring{R}_q\|_{C_{{[t-\ell,t+1]}}L^1})^4 r_{\perp}^{2/p-1} r_{\|}^{1/p-1/2} \lambda.
\endaligned
\end{equation}

For the termporal corrector, we have
\begin{equation}\label{corrector est2}
\aligned
&\|w_{q+1}^{(t)}\|_{C_{{[t,t+1]}}W^{1,p}}\lesssim {\mu}^{-1}\ \sum_{\xi \in \Lambda} \|a_{(\xi)}\|_{C_{{[t,t+1]},x}^0} \|a_{(\xi)}\|_{C_{{[t,t+1]},x}^1} \|\phi_{(\xi)}\|_{L^{2p}}^2 \|\psi_{(\xi)}\|_{C_{{[t,t+1]}}L^{2p}}^2\\
&\quad\quad+  \|a_{(\xi)}\|^2_{C_{{[t,t+1]},x}^0} \|\phi_{(\xi)}\|_{L^{2p}}  \|\nabla \phi_{(\xi)}\|_{L^{2p}}  \|\psi_{(\xi)}\|_{C_{{[t,t+1]}}L^{2p}}^2  \notag\\
 &\quad\quad+  \|a_{(\xi)}\|^2_{C_{{[t,t+1]},x}^0} \|\psi_{(\xi)}\|_{C_{{[t,t+1]}}L^{2p}}  \|\nabla \psi_{(\xi)}\|_{C_{{[t,t+1]}}L^{2p}}  \|\phi_{(\xi)}\|_{L^{2p}}^2.
\endaligned
\end{equation}
By (\ref{bounds}), \eqref{estimate aN} and (\ref{estimate aN0}) again, we have
\begin{equation}\label{eq:w:t:w1p}
\begin{aligned}
&\|w_{q+1}^{(t)}\|_{C_{{[t,t+1]}}W^{1,p}} \lesssim  \ell^{-15}(1+ \|\mathring{R}_q\|_{C_{{[t-\ell,t+1]}}L^1})^{3}{r_{\perp}^{2/p-2} r_{\|}^{1/p-1}}\mu^{-1} \lambda  .
\end{aligned}
\end{equation}

\subsubsection{Estimate of  $\|w_{q+1}\|_{C_{{[t,t+1]},x}^1}$}
 Applying \eqref{estimate aN} and \eqref{bounds} directly, we have
\begin{equation*}
\|w_{q+1}^{(p)}\|_{C_{{[t,t+1]},x}^1} \lesssim \ell^{-13} (1 + \|\mathring{R}_q\|_{C_{[a-\ell,b+\ell]}L^1}^2) r_{\perp}^{-1} r_{\|}^{-1/2} \lambda^2,
\end{equation*}
and
\begin{equation*}
\begin{aligned}
&\|w_{q+1}^{(c)}\|_{C_{{[t,t+1]},x}^1} \lesssim \|\textrm{curl}(\nabla a_{(\xi)} \times V_{(\xi)})\|_{C_{{[t,t+1]},x}^1} + \|\nabla a_{(\xi)} \times \textrm{curl}V_{(\xi)}\|_{C_{{[t,t+1]},x}^1} + \| a_{(\xi)} W^{(c)}_{(\xi)}\|_{C_{{[t,t+1]},x}^1} \\
\lesssim& \|a_{(\xi)}\|_{C_{{[t,t+1]},x}^3} (\|V_{(\xi)}\|_{ C_{{[t,t+1]},x}^2} + \|W_{(\xi)}^{(c)}\|_{C_{{[t,t+1]},x}^1}) \lesssim \ell^{-25}  (1 + \|\mathring{R}_q\|_{C_{{[t-\ell,t+1]}}L^1}^4) r_{\|}^{-3/2} \lambda^2.
\end{aligned}
\end{equation*}
For the temporal corrector $w_{q+1}^{(t)}$, by interpolation, we obtain that
\begin{equation*}
\begin{aligned}
&\|w_{q+1}^{(t)}\|_{C_{{[t,t+1]},x}^1} \lesssim \mu^{-1} \sum_{\xi \in \Lambda} \|a_{(\xi)}^2 \phi_{(\xi)}^2 \psi_{(\xi)}^2 \|_{C_{{[t,t+1]}} W^{1+\a,p}} +  \|a_{(\xi)}^2 \phi_{(\xi)}^2 \psi_{(\xi)}^2 \|_{C_{{[t,t+1]}}^1 W^{\a,p}} \\
\lesssim&\mu^{-1} \sum_{\xi \in \Lambda} \|a_{(\xi)}^2 \phi_{(\xi)}^2 \psi_{(\xi)}^2 \|_{C_{{[t,t+1]}} W^{1,p}}^{1-\a}  \|a_{(\xi)}^2 \phi_{(\xi)}^2 \psi_{(\xi)}^2 \|_{C_{{[t,t+1]}} W^{2,p}}^{\a} \\
&\ \ \ \ \ \ \ \ \ \ \ \ \ \ \ \ \ \ \ \ \ \ \ \  \ \ \ \ \ \ \ \ \ \ \ \ \ \ \ \ \ \ +  \|a_{(\xi)}^2 \phi_{(\xi)}^2 \psi_{(\xi)}^2 \|^{1-\a}_{C_{{[t,t+1]}}^1 L^{p}}\|a_{(\xi)}^2 \phi_{(\xi)}^2 \psi_{(\xi)}^2 \|_{C_{{[t,t+1]}}^1 W^{1,p}}^{\a} ,\\
\lesssim &{\ell^{-21}(1 + \|\mathring{R}_q\|_{C_{{[t-\ell,t+1]}}L^1}^3)r_{\perp}^{-1} r_{\|}^{-2} \lambda^{1+\alpha}.}
\end{aligned}
\end{equation*}
where $\a \in (0,1)$ and $p$ is large enough so that we can apply the Sobolev embedding $W^{\a,p}\hookrightarrow C^0$ to bound the $C_{t,x}^1$ norm with Sobolev norms. By applying \eqref{estimate aN} and \eqref{bounds}, there exists {some constant ${a}$} large enough such that 
{\begin{equation}\label{eq:w:c1}
\|w_{q+1}\|_{C_{{[t,t+1]},x}^1} \lesssim \ell_{q+1}^{-a}  (1 + \|\mathring{R}_q\|_{C_{{[t-\ell,t+1]}}L^1}^{4}). 
\end{equation}}

\subsection{Reynold stress at step $q+1$} \label{sec:w}
Substituting $v_{q+1}$ by \eqref{eq:v:q+1} and applying \eqref{eq:mollification} above, we can process similarly as in \cite{HZZ24} and obtain that
\begin{equation*}
\aligned
\mathring{R}_{q+1}=\mathring{R}_{\textrm{lin,z}}^{(q+1)}+ \mathring{R}^{(q+1)}_{\textrm{lin,w}}+&\mathring{R}_{\textrm{noise}}^{(q+1)}+ \mathring{R}_{\textrm{comm,1}}^{(q+1)} + \mathring{R}_{\textrm{comm,2}}^{(q+1)}\\
&+ \mathring{R}_{\textrm{mol}}^{(q+1)} + \mathring{R}^{(q+1)}_{\textrm{osc,x}} + \mathring{R}^{(q+1)}_{\textrm{osc,t}} + \mathring{R}^{(q+1)}_{\textrm{cor}}, 
\endaligned
\end{equation*}
where
\begin{equation}\label{eq:def:osc}
\aligned
\mathring{R}^{(q+1)}_{\textrm{osc,x}}:= \sum_{\xi \in \Lambda} \cB \(\nabla a^2_{(\xi)} , \mP_{\neq 0}\(W_{(\xi)} {\otimes} W_{(\xi)}\) \),\\
\mathring{R}^{(q+1)}_{\textrm{osc,t}}:=  -\frac{1}{\mu}\sum_{\xi \in \Lambda}  \mathcal{R}\(\mathbb{P}_{\neq0}\partial_t \(a_{(\xi)}^2\)\phi_{(\xi)}^2 \psi_{(\xi)}^2 \xi\),
\endaligned
\end{equation}
\begin{equation}\label{eq:def:w}
\aligned
\mathring{R}^{(q+1)}_{\textrm{cor}}: =& (w_{q+1}^{(t)} + w_{q+1}^{(c)}) \mathring{\otimes}_s w_{q+1}^{(p)} +  (w_{q+1}^{(t)} + w_{q+1}^{(c)}) \mathring{\otimes}  (w_{q+1}^{(t)} + w_{q+1}^{(c)}), \\
\mathring{R}^{(q+1)}_{\textrm{lin,w}}:=& \mathcal{R}\(\partial_t  (w_{q+1}^{(p)} + w_{q+1}^{(c)})  -   \Delta w_{q+1}\),
\endaligned
\end{equation}

and 
\begin{equation}\label{eq:def:zv}
\begin{aligned}
\mathring{R}_{\textrm{lin,z}}^{(q+1)}:=& \mathcal{R}(-c\ \mP_{\neq 0}\mP(\bar{z}_{q+1}- \bar{z}_{q,\ell})), \\
\mathring{R}_{\textrm{noise}}^{(q+1)}: =& (\bar{z}_{q+1}-\bar{z}_{q,\ell}) \mathring{\otimes} (\bar{z}_{q+1}- \bar{z}_{q,\ell}),\\
\mathring{R}_{\textrm{comm,1}}^{(q+1)}:=& (\bar{z}_{q+1} - \bar{z}_{q,\ell}) \mathring{\otimes}_s (v_{q+1}+ \bar{z}_{q,\ell}),\\
\mathring{R}_{\textrm{comm,2}}^{(q+1)}:=&  w_{q+1} \mathring{\otimes}_s (v_{q,\ell} + \bar{z}_{q,\ell}),\\
\mathring{R}_{\textrm{mol}}^{(q+1)}:=&(v_{q,\ell}+\bar{z}_{q,\ell})\mathring{\otimes}(v_{q,\ell}+\bar{z}_{q,\ell})-((v_q+\bar{z}_q)\mathring{\otimes}(v_q+\bar{z}_q))_{\ell}.
\end{aligned}
\end{equation}
 Here  we use the notation  $a\mathring{\otimes}_s b$ by the trace-free part of the symmetric tensor $a\otimes b+b\otimes a$. $\mathcal{R}$ is the inverse divergence operator introduced in Section \ref{sec:div-1}.

\subsubsection{Estimate of the stress terms}
We estimate all the terms in the Reynolds stress  $\mathring{R}_{q+1}$ here.  All the implicit constants in '$\lesssim$' are independent of the iteration step $q$ and the dissipative constant $c$.

By \eqref{eq:def:osc} and applying \eqref{eB}, we have
\begin{equation*}
\begin{aligned}
&\|\mathring{R}_{osc,x}^{(q+1)}\|_{C_{{[t,t+1]}} L^1} \lesssim \|\nabla a^2_{(\xi)}\|_{C_{{[t,t+1]},x}^1} \lv\mathcal{R}\(\mP_{\neq 0} W_{(\xi)} \otimes W_{(\xi)}\) \rv_{C_{{[t,t+1]}}L^1}\\
\lesssim& \|a_{(\xi)}\|_{C_{{[t,t+1]},x}^1} \| a_{(\xi)}\|_{C_{{[t,t+1]},x}^2}  \lv\mathcal{R}\(\mP_{\geq r_{\perp} \lambda/2} W_{(\xi)} \otimes W_{(\xi)}\) \rv_{C_{{[t,t+1]}}L^1}.
\end{aligned}
\end{equation*}
By \eqref{eR}, \eqref{bounds} and \eqref{estimate aN}, we have
\begin{equation*}
\|\mathring{R}_{osc,x}^{(q+1)}\|_{C_{{[t,t+1]}} L^1}  \lesssim \ell^{-32}(r_{\perp} \lambda)^{-1} (1+ \|\mathring{R}_q\|_{C_{{[t-\ell,t+1]}}L^1}^{5}).
\end{equation*}

By \eqref{eq:def:osc}, we have
\begin{equation*}
\begin{aligned}
&\|\mathring{R}_{osc,t}^{(q+1)}\|_{C_{{[t,t+1]}} L^1} \lesssim \mu^{-1} \|\mP_{\neq 0} \partial_t a_{(\xi)}^2 \phi_{(\xi)}^2 \psi_{(\xi)}^2 \xi\|_{C_{{[t,t+1]}}L^1}\\
\lesssim& \mu^{-1} \|\partial_t a_{(\xi)}\|_{C_{{[t,t+1]},x}^0}  \| a_{(\xi)}\|_{C_{{[t,t+1]},x}^0} \| \phi_{(\xi)}\|_{C_{{[t,t+1]}}L^2}^2  \| \psi_{(\xi)}\|_{C_{{[t,t+1]}}L^2}^2.
\end{aligned}
\end{equation*}
By \eqref{bounds} and \eqref{estimate aN}, we have
\begin{equation*}
\begin{aligned}
\|\mathring{R}_{osc,t}^{(q+1)}\|_{C_{{[t,t+1]}} L^1} \lesssim \ell^{-15} \mu^{-1}  (1+ \|\mathring{R}_q\|_{C_{{[t-\ell,t+1]}}L^1}^{5/2}).
\end{aligned}
\end{equation*}

By \eqref{eq:def:w}, we have
\begin{equation*}
\begin{aligned}
\|\mathring{R}^{(q+1)}_{\textrm{cor}}\|_{C_{{[t,t+1]}}L^1} \lesssim& \|(w_{q+1}^{(t)} + w_{q+1}^{(c)}) \mathring{\otimes}_s w_{q+1}^{(p)}\|_{C_{{[t,t+1]}}L^1} +  \|(w_{q+1}^{(t)} + w_{q+1}^{(c)}) \mathring{\otimes}  (w_{q+1}^{(t)} + w_{q+1}^{(c)})\|_{C_{{[t,t+1]}}L^1} \\
\lesssim& \|w_{q+1}^{(t)} + w_{q+1}^{(c)} \|_{C_{{[t,t+1]}}L^2} \| w_{q+1}^{(p)}\|_{C_{{[t,t+1]}}L^2} +  \|w_{q+1}^{(t)} + w_{q+1}^{(c)} \|_{C_{{[t,t+1]}}L^2}^2 .\\
\end{aligned}
\end{equation*}
By \eqref{est:principle:Lp}, \eqref{est:correction:Lp} and \eqref{est:temporal:Lp}, we have
\begin{equation}\label{eq:cor}
\begin{aligned}
\|\mathring{R}^{(q+1)}_{\textrm{cor}}\|_{C_{{[t,t+1]}}L^1} \lesssim \ell^{-21} \( {r_{\perp}}{r_{\|}}	^{-1} + \lambda^{-1} + \mu^{-1} r_{\perp}^{-1} r_{\|}^{-1/2} \)(1+ \|\mathring{R}_q\|_{C_{{[t-\ell,t+1]}}L^1}^{6}).
\end{aligned}
\end{equation}

By \eqref{eq:def:w} and \eqref{eR}, we have for some $\epsilon>0$ small enough,
\begin{equation*}
\begin{aligned}
&\|\mathring{R}^{(q+1)}_{\textrm{lin,w}}\|_{C_{{[t,t+1]}}L^1}\lesssim \|\mathcal{R}\partial_t  (w_{q+1}^{(p)} + w_{q+1}^{(c)})\|_{C_{{[t,t+1]}}L^{1+\epsilon}}  +  \|\mathcal{R}\Delta w_{q+1}\|_{C_{{[t,t+1]}}L^{1+\epsilon}}\\
\lesssim& \|\partial_t \sum_{\xi \in \Lambda} \mathcal{R}(\textrm{curl} \textrm{curl} (a_{(\xi)}  V_{(\xi)})) \|_{C_{{[t,t+1]}}L^{1+\epsilon}} + \|w_{q+1}\|_{C_{{[t,t+1]}} W^{1,{1+\epsilon}}}\\
\lesssim& \|a_{(\xi)}\|_{C_{{[t,t+1]},x}^2} \|\partial_tV_{(\xi)}\|_{C_{{[t,t+1]}} W^{1,{1+\epsilon}}}  + \|w_{q+1}\|_{C_{{[t,t+1]}} W^{1,{1+\epsilon}}}.
\end{aligned}
\end{equation*}
By \eqref{bounds} and \eqref{estimate aN}, there holds
\begin{equation}\label{eq:lin,w:1}
\begin{aligned}
&\|a_{(\xi)}\|_{C_{{[t,t+1]},x}^2} \|\partial_tV_{(\xi)}\|_{C_{{[t,t+1]}} W^{1,{1+\epsilon}}} 
\lesssim  \ell^{-19} (1+ \|\mathring{R}_q\|_{C_{{[t-\ell,t+1]}}L^1}^{3})  r_{\|}^{-1/2} r_{\perp}^2 \mu{\lambda^{3\epsilon}}.
\end{aligned}
\end{equation}
By \eqref{eq:w:pc:w1p} and \eqref{eq:w:t:w1p}, there holds
\begin{equation}\label{eq:lin,w:2}
\|w_{q+1}\|_{C_{{[t,t+1]}} W^{1,{1+\epsilon}}} \lesssim \ell^{-25}(1+ \|\mathring{R}_q\|_{C_{{[t-\ell,t+1]}}L^1}^{3}) (r_{\perp} r_{\|}^{1/2} \lambda+  \mu^{-1} \lambda + \mu^{-1}\lambda {r_{\perp}}{r_{\|}}^{-1}){\lambda^{3\epsilon}}.
\end{equation}
By \eqref{eq:lin,w:1} and \eqref{eq:lin,w:2}, we have
\begin{equation*}
\|\mathring{R}^{(q+1)}_{\textrm{lin,w}}\|_{C_{{[t,t+1]}}L^1}  \lesssim \ell^{-26}  (1+ \|\mathring{R}_q\|_{C_{{[t-\ell,t+1]}}L^1}^{3}) (r_{\perp} r_{\|}^{1/2} \lambda+ \mu^{-1} \lambda+ \mu^{-1}\lambda {r_{\perp}}{r_{\|}}^{-1} +  r_{\|}^{-1/2} r_{\perp}^2 \mu),
\end{equation*}
where we chose $\epsilon>0$ small enough such that $\lambda^{3\epsilon}\leq \ell^{-1}$.

By \eqref{eq:def:zv} and \eqref{eR} and standard mollification estimates, 
\begin{equation}
\begin{aligned}
\|\mathring{R}_{\textrm{lin,z}}^{(q+1)}\|_{C_{{[t,t+1]}}L^1} &\lesssim \|\mathcal{R}(-c\ \mP_{\neq 0}\mP(\bar{z}_{q+1}- \bar{z}_{q,\ell})) \|_{C_{{[t,t+1]}}L^1}\notag\\ 
&\lesssim   c \|\bar{z}_{q+1}- \bar{z}_{q} \|_{C_{{[t,t+1]}}L^2} + c \|\bar{z}_{q}- \bar{z}_{q,\ell} \|_{C_{{[t,t+1]}}L^2} \notag\\ 
&\lesssim c \|\bar{z}_{q+1}- \bar{z}_{q} \|_{C_{{[t,t+1]}}L^2} + c\ell^{\gamma}\| \bar{z}_{q} \|_{C^{\gamma}_{{[t-\ell,t+1]}}H^{\gamma}}.
\end{aligned}
\end{equation}

By \eqref{eq:def:zv} and standard convolution estimate, we have 
\begin{equation*}
\begin{aligned}
\|\mathring{R}_{\textrm{noise}}^{(q+1)}\|_{C_{{[t,t+1]}}L^1} \lesssim \|\bar{z}_{q+1}-\bar{z}_{q,\ell}\|^2_{C_{{[t,t+1]}}L^2}\lesssim \|\bar{z}_{q+1}- \bar{z}_{q} \|^2_{C_{{[t,t+1]}}L^2} +  \ell^{2 \g}\|\bar{z}_{q} \|_{C^{\g}_{{[t-\ell,t+1]}}H^{\g}}^2.
\end{aligned}
\end{equation*}

Similarly we have 
\begin{equation*}\label{eq:comm,1}
\begin{aligned}
&\|\mathring{R}_{\textrm{comm,1}}^{(q+1)}\|_{C_{{[t,t+1]}}{L^1}} \lesssim \|\bar{z}_{q+1} - \bar{z}_{q,\ell}\|_{C_{{[t,t+1]}}L^2} (\| v_{q+1}\|_{C_{{[t,t+1]}}L^2} + \|\bar{z}_{q,\ell}\|_{C_{{[t,t+1]}}L^2})\\
\lesssim& \( \|\bar{z}_{q+1}- \bar{z}_{q} \|_{C_{{[t-\ell,t+1]}}L^2} +  \ell_{}^{\g}\|\bar{z}_{q} \|_{C^{\g}_{{[t-\ell,t+1]}}H^{\g}}\) ( \|v_{q+1}\|_{C_{{[t,t+1]}}L^2} +  \|\bar{z}_q\|_{C_{{[t-\ell,t+1]}}L^2}).
\end{aligned}
\end{equation*}

By \eqref{eq:def:zv}, we have
\begin{equation*}
\begin{aligned}
&\|\mathring{R}_{\textrm{comm,2}}^{(q+1)}\|_{C_{{[t,t+1]}}{L^1}} \lesssim \|w_{q+1}\|_{C_{{[t,t+1]}}L^1} (\| v_{q,\ell}\|_{C_{{[t,t+1]}}L^\infty} + \|\bar{z}_{q,\ell}\|_{C_{{[t,t+1]}}L^\infty}).
\end{aligned}
\end{equation*}
By Young's convolution inequality  and Sobolev embedding, we have
\begin{equation*}
\| v_{q,\ell}\|_{C_{{[t,t+1]}}L^\infty} + \|\bar{z}_{q,\ell}\|_{C_{{[t,t+1]}}L^\infty} \lesssim \ell^{-4} (\| v_{q}\|_{C_{{[t-\ell,t+1]}}L^1} + \|\bar{z}_{q}\|_{C_{{[t-\ell,t+1]}}L^1}). 
\end{equation*}
Hence, by \eqref{eq:w:L1}, we obtain that
\begin{equation}\label{eq:comm,2}
\begin{aligned}
\|\mathring{R}_{\textrm{comm,2}}^{(q+1)}\|_{C_{{[t,t+1]}}{L^1}} &\lesssim \ell^{-23}\lambda_q^{d\epsilon} r_\perp^{}r_\|^{1/2} (1+\|\mathring{R}_q\|^{3}_{C_{{[t-\ell,t+1]}}L^1} )(\| v_{q}\|_{C_{{[t-\ell,t+1]}}L^1} + \|\bar{z}_{q}\|_{C_{{[t-\ell,t+1]}}L^1})\notag\\
&\lesssim \ell^{-24}r_\perp^{}r_\|^{1/2} (1+\|\mathring{R}_q\|^{3}_{C_{{[t-\ell,t+1]}}L^1} )(\| v_{q}\|_{C_{{[t-\ell,t+1]}}L^1} + \|\bar{z}_{q}\|_{C_{{[t-\ell,t+1]}}L^1}),
\end{aligned}
\end{equation}
where we chose $\epsilon>0$small enough to absorb the universl  constant.

By \eqref{eq:def:zv} and standard mollification estimates, we have
\begin{equation*}
\begin{aligned}
&\|\mathring{R}_{\textrm{mol}}^{(q+1)}\|_{C_{{[t,t+1]}}L^1}\lesssim \|(v_{q,\ell}+\bar{z}_{q,\ell})\mathring{\otimes}(v_{q,\ell}+\bar{z}_{q,\ell})-((v_q+\bar{z}_q)\mathring{\otimes}(v_q+\bar{z}_q))_{\ell}\|_{C_{{[t,t+1]}}L^1}\\
\lesssim& \ell^{\gamma} (\|v_q\|_{C_{{[t-\ell,t+1]}}^{\g}H^{\g}} + \|\bar{z}_{q} \|_{C^{\gamma}_{{[t-\ell,t+1]}}H^{\gamma}} )(\|v_q\|_{C_{{[t-\ell,t+1]}}L^2}+\|\bar{z}_q\|_{C_{{[t-\ell,t+1]}}L^2}).
\end{aligned}
\end{equation*}

{In summary, we have the following bounds.}
\begin{equation}\label{Rq+1:l1:r}
\begin{aligned}
&\|\mathring{R}_{q+1} \|_{C_tL^1} \lesssim\ \ell_{q+1}^{-32} \lambda_{q+1}^{-1/8}(1+ \|\mathring{R}_q\|_{C_{[t-\ell,t+1]}L^1}^{6})(1+\| v_{q}\|_{C_{[t-\ell,t+1]}L^{2}} + \|\bar{z}_{q}\|_{C_{[t-\ell,t+1]} L^{2}})\\
& + \(\|\bar{z}_{q+1}- \bar{z}_{q} \|_{C_{[t-\ell,t+1]}L^2}+  \ell_{q+1}^{\g}\|\bar{z}_{q} \|_{C^{\g}_{[t-\ell,t+1]}H^{\g}} \)+ \(\|\bar{z}_{q+1}-\bar{z}_{q} \|^2_{C_{t}L^2} +  \ell_{q+1}^{2\g}\|\bar{z}_{q} \|_{C^{\g}_{[t-\ell,t+1]}H^{\g}}^2\)\\
&+ \( \|\bar{z}_{q+1}- \bar{z}_{q} \|_{C_{[t-\ell,t+1]}L^2} +  \ell_{q+1}^{\g}\|\bar{z}_{q} \|_{C^{\g}_{[t-\ell,t+1]}H^{\g}}\) ( \|v_{q+1}\|_{C_{t}L^2} +  \|\bar{z}_q\|_{C_{[t-\ell,t+1]}L^2})\\
&+\ell_{q+1}^{\g} (\|v_q\|_{C_{[t-\ell,t+1],x}^1} + \|\bar{z}_{q} \|_{C^{\g}_{[t-\ell,t+1]}H^{\g}} )(\|v_q\|_{C_{[t-\ell,t+1]}L^2}+\|\bar{z}_q\|_{C_{[t-\ell,t+1]}L^2}).
\end{aligned}
\end{equation}

Moreover, {we have the following direct bounds which are convenient to   estimate the general $m$-th moment.}
\begin{equation*}
\begin{aligned}
\|\mathring{R}_{\textrm{lin,z}}^{(q+1)}\|_{C_{{[t,t+1]}}L^1}  \lesssim \|\bar{z}_{q+1}\|_{C_{{[t,t+1]}}L^2}+ \|\bar{z}_q\|_{C_{{[t-\ell,t+1]}}L^2}.
\end{aligned}
\end{equation*}
\begin{equation*}
\|\mathring{R}_{\textrm{noise}}^{(q+1)}\|_{C_{{[t,t+1]}}L^1} \lesssim \|\bar{z}_{q+1}\|_{C_{{[t,t+1]}}L^2}^2 + \|\bar{z}_q\|_{C_{{[t-\ell,t+1]}} L^2}^2.
\end{equation*}
\begin{equation*}
\|\mathring{R}_{\textrm{comm,1}}^{(q+1)}\|_{C_{{[t,t+1]}}{L^1}} \lesssim \|\bar{z}_{q+1}\|_{C_{{[t,t+1]}}L^2}^2+ \|\bar{z}_q\|_{C_{{[t-\ell,t+1]}}L^2}^2+\|v_{q+1}\|_{C_{{[t,t+1]}}L^2}^2.
\end{equation*}
\begin{equation*}
\|\mathring{R}_{\textrm{mol}}^{(q+1)}\|_{C_{{[t,t+1]}}L^1} \lesssim \|v_q\|_{C_{{[t-\ell,t+1]}}L^2}^2 + \|\bar{z}_q\|_{C_{{[t-\ell,t+1]}}L^2}^2.
\end{equation*}

{In summary, we have}
\begin{equation}
\begin{aligned}\label{Rq+1:l1:m}
&\|\mathring{R}_{q+1} \|_{C_tL^1} \lesssim (1+ \|\mathring{R}_q\|_{C_{[t-\ell,t+1]}L^1}^6)(1+\| v_{q}\|_{C_{[t-\ell,t+1]]}L^{2}} + \|\bar{z}_{q}\|_{C_{[t-\ell,t+1]} L^{2}})\\
& + \|\bar{z}_{q+1}\|_{C_{[t-\ell,t+1]}L^2}^{2} + \|\bar{z}_{q}\|_{C_{[t-\ell,t+1]}L^2}^{2}+ \|v_{q+1}\|_{C_{[t-\ell,t+1]}L^2}^{2}+ \|v_{q}\|_{C_{[t-\ell,t+1]}L^2}^{2}.
\end{aligned}
\end{equation}
\subsubsection{Iterative estimates on the stress}
\begin{lemma}\label{l:iterative:R:m}
There exists some universal constant  $\a_0\geq 12r_0$ such that

\begin{equation}\label{eq:a:R:m:q+1}
\ell_{q+1}^{-1}+\$\mathring{R}_{q+1}\$_{L^1,m} + \$v_{q+1}\$_{C_{t,x}^1,m} \leq  \ell_{q+1}^{-\a_0} M_{q}(\a_0,m),
\end{equation}
for all $m \geq 1$.
\end{lemma}
\begin{proof}
By  Lemma \ref{l:truncate}, there holds
\begin{equation*}
\sup_{t \in \mathbb{R}}\|\bar{z}_{q+1}\|_{C_tL^2} + \|\bar{z}_{q}\|_{C_tL^2}\lesssim \zeta_{q+1}^{4}= \ell_{q+1}^{-4}.
\end{equation*}
Then by \eqref{Rq+1:l1:m} there holds 
\begin{equation}\label{eq:R:largebound}
\begin{aligned}
\$\mathring{R}_{q+1} \$_{L^1,m} &\lesssim
(1+ \$\mathring{R}_q\$_{L^1,12m}^{6})(1+\$ v_{q}\$_{L^{2},2m} +\ell_{q+1}^{-4}) +\ell_{q+1}^{-8}+ \$v_{q+1}\$_{L^2,2m}^{2} +   \$v_{q}\$_{L^2,2m}^{2}.
\end{aligned}
\end{equation}

By \eqref{eq:w:c1}, there holds for $a$ large enough
\begin{equation*}
\begin{aligned}
\$v_{q+1}\$_{C_{t,x}^1,2m} \lesssim \$w_{q+1}\$_{C_{t,x}^1,2m} + \$v_{q,\ell}\$_{C_{t,x}^1,2m} \lesssim \ell_{q+1}^{-a} (1+ \$\mathring{R}_{q}\$_{L^1, 4m}^4 + \$v_{q}\$_{C_{t,x}^1,2m}),
\end{aligned}
\end{equation*}
which implies
\begin{equation*}
\$v_{q+1}\$_{C_{t,x}^1,2m} \lesssim \ell_{q+1}^{-a} M_q(a,m).
\end{equation*}

Going back to  the right hand side of \eqref{eq:R:largebound},  there exists some constant $\a_0$ large enough such that \eqref{eq:a:R:m:q+1} holds.

\end{proof}

\begin{lemma}\label{l:M:bound}
Assume Assumption \ref{a:iterative} holds up to step $q$.  Then  Assumption \ref{a:iterative}(3) holds at step $q+1$, i.e.   for any $m \geq 1$, there holds
\begin{equation*}
\Lambda_{q+1} \leq M_{q+1}(\a_0,m) \leq  \eps_{q+2}^{-1}.
\end{equation*}
\end{lemma}
\begin{proof}
 By \eqref{eq:M:bound:0}, we only need to estimate the second term above.
For any  $m\geq 1$, by Lemma \ref{l:iterative:R:m}, we have 
\begin{equation*}
\begin{aligned}
&M_{q+1}(\a_0,m) = \$\mathring{R}_{q+1}\$_{L^1,\a_0m}^{\a_0} + \$v_{q+1}\$_{C_{t,x}^1,\a_0 m}^{\a_0} +\ell_{q+1}^{-1}\\
\leq&\ (\$\mathring{R}_{q+1}\$_{L^1,\a_0m} + \$v_{q+1}\$_{C_{t,x}^1,\a_0m}+\ell_{q+1}^{-1})^{\a_0}\leq \ \ell_{q+1}^{-\a_0^2} M_{q}(\a_0, \a_0 m)^{\a_0}.
\end{aligned}
\end{equation*}
By \eqref{eq:M:bound} there holds
\begin{equation*}
\ell_{q+1}^{-\g} = \Lambda_{q}\ \eps_{q+1}^{-N_0} \leq M_{q}(\a_0,m)\eps_{q+1}^{-N_0}\leq \eps_{q+1}^{-N_0-1},
\end{equation*}
which implies that there exists some universal constant $M_0$ such that
\begin{equation*}
 M_{q+1}(\a_0,m)   \leq  \eps_{q+1}^{-(N_0+1)\a_0^2/\g-\a_0} \leq \eps_{q+1}^{-M_0}= (\Xi^{-\a_0 M_0^{q-1}})^{-M_0}= \eps_{q+2}^{-1}\ .
\end{equation*}
\end{proof}

\begin{lemma}\label{l:iterative:R:p0}
Assume Assumption \ref{a:iterative} holds up to step $q$. There exists some universal constant $c_R$ such that Assumption \ref{a:iterative}(1) holds at step $q+1$, i.e. 
\begin{equation}\label{eq:a:R:p:q+1}
\$\mathring{R}_{q+1}\$_{L^1,r_0}  \leq c_R  \eps_q^{}.
\end{equation}
\end{lemma}
\begin{proof}
By \eqref{Rq+1:l1:r} we obtain
\begin{equation*}
\begin{aligned}
&\$\mathring{R}_{q+1} \$_{L^1,r_0} \lesssim \ell_{q+1}^{-32} \lambda_{q+1}^{-1/8}(1+ \$\mathring{R}_q\$_{L^1,12r_0}^{12}+\$ v_{q}\$_{L^{2},2r_0}^2 + \$\bar{z}_{q}\$_{L^{2},2r_0}^2) + \bar{Z}_{q+1}^{1/r_0}\\
& +\bar{Z}_{q+1}^{1/(2r_0)} (1+ \$v_{q+1}\$_{L^2,2r_0}+ \$\bar{z}_{q}\$_{L^2,2r_0})+\ell_{q+1}^{\g}(1+\$v_{q}\$_{C_{t,x}^1,2r_0}^2+ \$\bar{z}_{q}\$_{C^{\g}_tH^{\g},2r_0}^2).
\end{aligned}
\end{equation*}
{By Lemma \ref{l:truncate}, \eqref{zq:moment} and \eqref{eq:v:p0} we obtain that
\begin{equation*}
\E[\|\bar{z}_{q}\|_{C_{[t,t+1]} L^{2}}^{2r_0}] \leq \E[\|{z}_{q} \|_{C^{\g}_{[t,t+1]}H^{\g}}^{2r_0}] \lesssim 1+\$v_{q-1}\$_{L^2,2r_0}^{2r_0}\lesssim1.
\end{equation*}}
Then by \eqref{eq:v:p0} and Lemma \ref{l:barZ:q} we obtain that
\begin{equation*}
\begin{aligned}
&\$\mathring{R}_{q+1} \$_{L^1,r_0} \lesssim\ \ell_{q+1}^{\g}\Lambda_q+ \bar{Z}_{q+1}^{1/r_0} +\bar{Z}_{q+1}^{1/(2r_0)}\lesssim\ {\ell_{q+1}^{\g}\Lambda_q+\ell_{q}^{\g}\Lambda_{q-1}\lesssim \eps_{q+1}+\eps_{q}}.
\end{aligned}
\end{equation*}
Therefore, there exists some constant $c_R$ such that \eqref{eq:a:R:p:q+1} holds.
\end{proof}

\subsection{Estimate of energy}\label{sec:est:energy}
We conclude this section by providing the estimate of energy.
\begin{lemma}\label{l:theta:E}
Assume Assumption \ref{a:iterative} holds up to step $q$.  There exists some universal constant $c_{\theta}$ such that Assumption \ref{a:iterative}(4) holds at step $q+1$, i.e. 
\begin{equation*}
0 < \theta_{q+1}(t) \leq c_{\theta} \eps_{q-2},\ \ \ \delta E_{q+1}(t) \leq c_E \eps_{q-1}^{}.
\end{equation*}
\end{lemma}
\begin{proof}
 Recalling the definition of $v_{q+1}$, we have
\begin{equation*}
\begin{aligned}
 &\mathbf{E}[\|v_{q+1}(t) + \bar{z}_{q+1}(t)\|_{L^2}^2]=  \mathbf{E}[\|v_{q}(t) + \bar{z}_{q}(t)+ v_{q,\ell}(t) -v_{q}(t) + \bar{z}_{q+1}(t)-\bar{z}_q(t)+ w_{q+1}(t)\|^2_{L^2}]\\
=&  \mathbf{E}[\|v_{q}(t) + \bar{z}_{q}(t)\|_{L^2}^2]+  \mathbf{E}[\|w_{q+1}(t)\|^2_{L^2}]+ \mathbf{E}[\|v_{q,\ell}(t) -v_{q}(t)\|_{L^2}^2] +\mathbf{E}[ \|\bar{z}_{q+1}(t)-\bar{z}_q(t)\|_{L^2}^2]\\
& + 2\mathbf{E}[\|(v_q(t)+ \bar{z}_q(t)) w_{q+1}(t)\|_{L^1}] +2 \mathbf{E}[\|(v_q(t)+\bar{z}_q(t))(v_{ q,\ell}(t)- v_q(t))\|_{L^1}] \\
&+ 2\mathbf{E}[\|(v_q(t)+ \bar{z}_q(t)) (\bar{z}_{q+1}{(t)}-\bar{z}_q(t))\|_{L^1}] +2 \mathbf{E}[\|(v_{ q,\ell}(t)-v_q(t)) w_{q+1}(t)\|_{L^1}]\\
&+ 2\mathbf{E}[\|(\bar{z}_{q+1}(t) - \bar{z}_q(t)) w_{q+1}(t)\|_{L^1}] +2  \mathbf{E}[\|(\bar{z}_{q+1}(t)-\bar{z}_q(t)) (v_{q,\ell}(t)-v_q(t))\|_{L^1}].
\end{aligned}
\end{equation*}
Hence, there holds
\begin{equation*}
\begin{aligned}
&\delta E_{q+1}(t)^2 \lesssim |e(t)(1-\eps_{q-2}) - \mathbf{E}[\|v_{q}(t) + \bar{z}_{q}(t)\|_{L^2}^2]- \mathbf{E}[\|w_{q+1}{(t)}\|^2_{L^2}] |^2\\
&+\mathbf{E}[\|v_{q,\ell}(t) -v_{q}(t)\|_{L^2}^2]^2+ \mathbf{E}[\|\bar{z}_{q+1}(t)-\bar{z}_q(t)\|_{L^2}^2]^2+ \mathbf{E}[\|(v_q(t)+ \bar{z}_q(t))w_{q+1}(t)\|_{L^1}] ^2\\
&+ \mathbf{E}[\|v_q(t)+\bar{z}_q(t)\|_{L^2}^2]\mathbf{E}[\|v_{ q,\ell}(t)- v_q(t)\|_{L^2}^2]+ \mathbf{E}[\|v_q(t)+ \bar{z}_q(t)\|_{L^2}^2] \mathbf{E}[\|\bar{z}_{q+1}{(t)}- \bar{z}_q(t)\|_{L^2}^2] \\
&+ \mathbf{E}[\|v_{ q,\ell}(t)-v_q(t)\|_{L^2}^2] \mathbf{E}[ {\|}w_{q+1}(t)\|_{L^2}^2] + \mathbf{E}[\|\bar{z}_{q+1}(t)-\bar{z}_q(t)\|_{L^2}^2] \mathbf{E}[\| w_{q+1}(t)\|_{L^2}^2]\\
& +  \mathbf{E}[\|\bar{z}_{q+1}(t)-\bar{z}_q(t)\|_{L^2}^2] \mathbf{E}[ \|v_{q,\ell}(t)-v_q(t)\|_{L^2}^2].
\end{aligned}
\end{equation*}
By a standard mollification estimate, we have
\begin{equation}\label{eq:dE:vql}
\sup_{t \geq0}\mathbf{E}[\|v_{q,\ell}(t) -v_{q}(t)\|_{L^2}^2]\lesssim \ell_{q+1}^2 \$v_q\$_{\C^1,2r_0}^2 \lesssim \ell_{q+1}^2 \Lambda_{q}.
\end{equation}
By  Lemma \ref{l:barZ:q}, we have
\begin{equation}\label{eq:dE:z}
\sup_{t \geq0}\mathbf{E}[\|\bar{z}_{q+1}(t)-\bar{z}_q(t)\|_{L^2}^2] \lesssim \ell_{q}^{\g} \Lambda_{q-1}.
\end{equation}
Moreover, by \eqref{zq:moment}, \eqref{eq:v:p0}, \eqref{bd:wq+1l2} and Lemma \ref{l:truncate}  we have
{\begin{equation}\label{eq:dE:wvz}
\begin{aligned}
 \mathbf{E}[ \|w_{q+1}(t)\|_{L^2}^2+
 \|v_q(t)+\bar{z}_q(t)\|_{L^2}^2] \lesssim& c_\theta+\$ v_q\$_{L^2,2}^2+\$ z_{q}\$_{L^2,2}^2\\
 \lesssim& c_\theta+\$ v_q\$_{L^2,2r_0}^2+\$v_{q-1}\$_{L^2,2r_0}^2\lesssim {c}_v+c_\theta.
 \end{aligned}
\end{equation}}
Combining \eqref{eq:dE:vql}, \eqref{eq:dE:z} and \eqref{eq:dE:wvz}, we obtain that
\begin{equation*}
\begin{aligned}
\delta E_{q+1}(t)^2 \lesssim&  |e(t)(1-\eps_{q-2}) - \mathbf{E}[\|v_{q}(t) +  \bar{z}_{q}(t)\|_{L^2}^2]- \mathbf{E}[\|w_{q+1}(t)\|^2_{L^2}] |^2\\
&+(c_{\theta}+c_v+ c_{R}+1)(\ell_{q+1}^2 \Lambda_{q}+ \ell_{q}^{\g} \Lambda_{q-1}) + \mathbf{E}[\|(v_q(t)+ \bar{z}_q(t))w_{q+1}(t)\|_{L^1}]^2.
\end{aligned}
\end{equation*}
By taking $N_0$ large enough, so that
{\begin{equation}\label{eq:N0}
(c_{\theta}+c_v+ c_{R}+1)
(\ell_{q+1}^2 \Lambda_{q}+ \ell_{q}^{\g} \Lambda_{q-1}) \leq  2(c_{\theta}+c_v+ c_{R}+1) \eps_q^{N_0} \leq \eps_q^2,
\end{equation}}
and the following holds:
\begin{equation}\label{eq:dE:1}
\begin{aligned}
\delta E_{q+1}(t)^2 \lesssim\ & |e(t)(1-\eps_{q-2}) - \mathbf{E}[\|v_{q}(t) +  \bar{z}_{q}(t)\|_{L^2}^2]- \mathbf{E}[\|w_{q+1}(t)\|^2_{L^2}] |^2  \\
&+ \mathbf{E}[\|(v_q(t)+ \bar{z}_q(t))w_{q+1}(t)\|_{L^1}]^2 +  \eps_{q}^{2}.
\end{aligned}
\end{equation}

To estimate the term $ \mathbf{E}[\|(v_q(t)+ z_q(t))w_{q+1}(t)\|_{L^1}]$, we apply  Lemma \ref{l:truncate}, \eqref{para:l} and \eqref{bd:wq+p0} to obtain
{\begin{equation*}
\begin{aligned}
\mathbf{E}[\|(v_q(t)&+ \bar{z}_q(t))w_{q+1}(t)\|_{L^1}]  \lesssim  \$v_q+ \bar{z}_q\$_{L^\infty,2}\$w_{q+1}\$_{L^{p_0},2r_0}\\
\lesssim& \ell_{q+1}^{-19} \lambda_{q+1}^{-\frac98(\frac2{p_0}-1)} \Lambda_{q}^{1/2}(\Lambda_{q}^{1/2}+\ell_{q+1}^{-4})\lesssim  \ell_{q+1}^{-23} \lambda_{q+1}^{-\frac98(\frac2{p_0}-1)} \Lambda_{q}\lesssim \ell_{q+1}^{\g}\Lambda_q.
\end{aligned}
\end{equation*}}

For the second term in \eqref{eq:dE:1}, we have
\begin{equation*}
\begin{aligned}
 &|e(t)(1-\eps_{q-2}) - \mathbf{E}[\|v_{q}(t) + \bar{z}_{q}(t)\|_{L^2}^2]- \mathbf{E}[\|w_{q+1}\|^2_{L^2}] | \\
 \leq& \E[\|(w_{q+1}^{(t)} +w_{q+1}^{(c)})\otimes_s w_{q+1}^{(p)}\|_{L^1}] + \E[\|w_{q+1}^{(t)} +w_{q+1}^{(c)}\|_{L^2}^2]\\
 &+|e(t)(1-\eps_{q-2}) - \mathbf{E}[\|v_{q}(t) + \bar{z}_{q}(t)\|_{L^2}^2]- \mathbf{E}[\|w^{(p)}_{q+1}\|^2_{L^2}] |.
\end{aligned}
\end{equation*}
By \eqref{para:l} and \eqref{eq:cor}, there holds
\begin{equation*}
 \E[\|(w_{q+1}^{(t)} +w_{q+1}^{(c)})\otimes_s w_{q+1}^{(p)}\|_{L^1}] + \E[\|w_{q+1}^{(t)} +w_{q+1}^{(c)}\|_{L^2}^2] \lesssim  {\ell_{q+1}^{-21}\lambda_{q+1}^{-1/8}\Lambda_q\lesssim   \ell_{q+1}^{\g}} \Lambda_{q}.
\end{equation*}
By taking the trace in \eqref{can}, by \eqref{eq:rho3} we have
\begin{align*}
|  w_{q+1}^{(p)}|^2-3\theta_{q}=6\sqrt{\ell_{q+1}^2+|\mathring{R}_{q,\ell}|^2}+3(\theta_{q,\ell}-\theta_q)+\sum_{\xi\in \Lambda}a_{(\xi)}^2\mathbb{P}_{\neq0}|W_{(\xi)}|^2.
\end{align*}
Hence, there holds
\begin{equation}\label{eq:gg ps}
\begin{aligned}
\mathbf{E}|\|  w_{q+1}^{(p)}\|_{L^2}^2-3\theta_{q}(2\pi)^3|\leq& 6\cdot(2\pi)^3\ell_{q+1}+6\mathbf{E}\|\mathring{R}_{q,\ell_{}}\|_{L^1}+3\cdot(2\pi)^3|\theta_{q,\ell_{}}-\theta_q|\\
&+\mathbf{E}\sum_{\xi\in \Lambda}\Big|\int a_{(\xi)}^2\mathbb{P}_{\neq0}|W_{(\xi)}|^2{\dif x}\Big|.
\end{aligned}
\end{equation}
By a standard mollification estimate and \eqref{zq:moment} we have
\begin{equation*}
\begin{aligned}
|\theta_{q,\ell}(t)&- \theta_q(t)| \lesssim \ell_{q+1} \|e\|_{C_t^1} + \left|\E[\|v_{q} + \bar{z}_q\|_{L^2}^2] *_t \varphi_{\ell_{q+1}}(t)-\E[\|v_{q}(t)  + \bar{z}_q(t)\|_{L^2}^2]\right|\\
\lesssim& \ell_{q+1} + \ell_{q+1}^{\gamma}  \E[\|v_{q} + \bar{z}_q\|_{C_t^{\gamma}L^2}^2] \lesssim{\ell_{q+1} + \ell_{q+1}^{\gamma}(\$v_{q}\$_{C^{1}_{t,x},2}^2 +\$\bar{z}_q\$_{C^{\gamma}_{t}L^2,2r_0}^2)
\lesssim  \ell_{q+1}^{\gamma} \Lambda_q.}
\end{aligned}
\end{equation*}
For the last term in \eqref{eq:gg ps},  we have
\begin{equation*}
\begin{aligned}
&\sum_{\xi\in\Lambda}\Big|\int_{\mathbb{T}^{3}} a_{(\xi)}^2\mathbb{P}_{\neq0}|W_{(\xi)}|^2 \dif x\Big|
=\sum_{\xi\in\Lambda}\Big|\int_{\mathbb{T}^{3}} a_{(\xi)}^2{\mathbb{P}}_{\geq r_\perp\lambda_{q+1}/2}|W_{(\xi)}|^2\dif x\Big|\\
=&\sum_{\xi\in\Lambda}\Big|\int_{\mathbb{T}^{3}} |\nabla|^Na_{(\xi)}^2|\nabla|^{-N}\mathbb{P}_{\geq r_\perp\lambda_{q+1}/2}|W_{(\xi)}|^2 \dif x\Big|\lesssim\|a_{(\xi)}^2\|_{C_{t,x}^N}(r_\perp\lambda_{q+1})^{-N}\||W_{(\xi)}|^2\|_{C_tL^2}.
\end{aligned}
\end{equation*}
By \eqref{estimate aN}, \eqref{estimate aN0}, there holds
\begin{equation*}
\begin{aligned}
\sum_{\xi\in\Lambda}\Big|\int_{\mathbb{T}^{3}} a_{(\xi)}^2\mathbb{P}_{\neq0}|W_{(\xi)}|^2 \dif x\Big|
&\lesssim\ell^{-6N-9}(\|\mathring{R}_q\|_{C_{[t-\ell,t+1]}L^1}+1)^{N+3/2}(r_\perp\lambda_{q+1})^{-N}r_\perp^{-1}r_{\|}^{-\frac12}\\
&\lesssim \ell^{-6N-9}(\|\mathring{R}_q\|_{C_{[t-\ell,t+1]}L^1}+1)^{N+3/2}\lambda_{q+1}^{\frac{9-N}8}.
\end{aligned}
\end{equation*}
Hence, by taking {expectation and $N$ large enough and \eqref{para:l} we have
\begin{align}\label{eq:dE:aux1}
    \mathbf{E}\sum_{\xi\in \Lambda}\Big|\int_{\mathbb{T}^3} a_{(\xi)}^2\mathbb{P}_{\neq0}|W_{(\xi)}|^2 \dif x\Big|\lesssim \ell^{-6N-9}_{q+1}\Lambda_q^N\ell_{q+1}^{10(N-9)}\lesssim \Lambda_q^N\ell_{q+1}^{N}\lesssim \eps_{q+1}.
\end{align}}

By Lemma \ref{l:truncate}, and \eqref{eq:R:p0}, there holds 
\begin{equation*}
\mathbf{E}|\|  w_{q+1}^{(p)}\|_{L^2}^2-3\theta_{q+1}(2\pi)^3|\lesssim c_R \eps_{q-1}^{}.
\end{equation*}

Therefore, combining all the bounds above together, there exists some constant $c_E$ such that
\begin{equation*}
\d E_{q+1}(t) \leq c_E \eps_{q-1}.
\end{equation*}
 
For $\theta_{q+1}$, there holds
\begin{equation*}
\begin{aligned}
|3\cdot(2\pi)^3{\theta}_{q+1}(t)- e(t)(\eps_{q-2}-\eps_{q-1})|= \d E_{q+1}(t).
\end{aligned}
\end{equation*} 
Let $\bar{e}>\underline{e}>2 c_E$. By \eqref{eq:eps}, there exists some constant $c_{\theta}$ such that
\begin{equation*}
0<\underline{e}(\eps_{q-2}-\eps_{q-1}) -  c_E\eps_{q-1} \leq 3\cdot(2\pi)^3{\theta}_{q+1}(t) \leq \bar{e}(\eps_{q-2}-\eps_{q-1})+ c_E\eps_{q-1} \leq 3\cdot(2\pi)^3c_{\theta} \eps_{q-2}.
\end{equation*}
\end{proof}

\section{Construction of non-unique solutions with Cauchy problem}
\label{sec:iterative:cut}
This section is devoted to prove Theorem \ref{thm:main:3}  where the convex integration follows  a similar way as in Section \ref{sec:iteration} with a little bit different since we need to  keep the initial condition. 

 As in Section \ref{sec:stationary} and Section \ref{sec:iteration}, Assumption \ref{a:G:1} holds for some fixed parameters $(p_0,\d_0)$, the parameters $\g$, $r_0$ satisfy condition \eqref{def:gamma:p} and the dissipative constant $c$ is large enough so that a similar result to \eqref{zq:moment} holds:
\begin{align}
\sup_{t \in \mathbb{R}}\E[\|z_q\|_{C_t^{\g}H^{\g}}^{2r_0}] \leq C_{c,r_0,\g,\d_0,L}(1+ \$v_{q-1}\$_{L^{p_0},2r_0}^{2r_0}+\$ \z\$_{L^{p_0},2r_0}^{2r_0}),\notag\\
\sup_{t \in \mathbb{R}}\E[\|z_q\|_{C_t^{\g}H^{\g}}^{4r_0}] \leq C_{c,r_0,\g,\d_0,L}(1+ \$v_{q-1}\$_{L^{p_0},4r_0}^{4r_0}+\$ \z\$_{L^{p_0},4r_0}^{4r_0}),\notag\\
{\sup_{t \in \mathbb{R}}\E[\|z_{q+1}-z_q\|_{C_t^{\g}H^{\g}}^{2r_0}] \leq C_{c,r_0,\g,\d_0,L} \$v_q-v_{q-1}\$_{L^{p_0},2r_0}^{2r_0}.\label{zq:moment:1}}
\end{align}

Moreover, as in Section \ref{sec:iteration}, the constants $\a_0$ (cf. Lemma \ref{l:M:bound:i}), $M_0$ (cf. Lemma \ref{l:M:bound:i}), $N_0$ (cf. \eqref{l:iterative:R:p0:i}) and $\Xi$ (cf. Proposition \ref{p:iterative}) below are chosen to be sufficiently large.
 
 We define $u_0\in L^2_\sigma$ by the initial condition.  Then  for a dissipative constant $c>0$, we define $\z(t)$ by
\begin{align}
{\z}(t):= e^{-ct}S(t)u_0 .\label{def:dotz}
\end{align}

At each step $q \in \mathbb{N}$, we {assume that}  the pair $(v_q, \mathring{R}_q)$ solves the following equation:
\begin{equation}\label{induction ps:i}
\aligned
\partial_tv_{q}-c \bar{z}_{q}-\Delta v_{q} +\div((v_{q}+\bar{z}_{q})\otimes (v_{q}+\bar{z}_{q}))+\nabla p_{q}\dif t&=\div \mathring{R}_{q},\ 
\\
\div \ v_{q}&=0,\ \\
{v_q(0)}&=0,
\endaligned\ \ \ 
\end{equation} 
where we define  $\bar{z}_q(t):= \Pi_{\zeta_{q}} z_{q}(t) + \z(t)$. Here $\zeta_q>0$ is the truncation parameter, 
 $z_{q}(t)$  is the solution of the following stochastic partial differential equation:{
\begin{equation}\label{eq:sc:i}
\begin{aligned}
\dif  z_q(t) -  (\Delta - cI) z_q(t) \dif t + \nabla p_{z,q}\dif t =& G(v_{q-1}(t)+ \z(t)+ z_q(t)) \dif W_t,\\
\div z_q=&0,\\
z_q(0)=&0,
\end{aligned}
\end{equation}
for $t \geq 0$.} To handle the mollification around $t=0$, we define $z_{q}(t)=0$ and  $\z(t)=u_0$ for $t<0$. 

\begin{definition}\label{def:para:2:i}
For $q \geq 1$, we define $\Lambda_{q}$ and $M_{q}(n,p)$ as in Definition \ref{def:para:3}. Let $\eps_{1}=1$. For $q \geq 2$, we define $\eps_{q}$ by
\begin{equation*}
\eps_{q} =\Xi^{-\a_0{M_0}^{q-2}}.
 \end{equation*}

Let $\ell_1=1/2$
. For $q\geq 2$, let $\ell_{q}$ be such that 
\begin{equation}\label{def:ell:i}
\ell_q^{\g} \Lambda_{q-1}^{2}:=  \eps_{q}^{N_0}.
\end{equation}
 {Then we define \begin{align*}
  \lambda_q=   (\lfloor  \ell_q^{-10\g/3}+\ell_q^{-3p_0/(2-p_0)}+\ell_q^{-40}\rfloor+1)^8,
\end{align*}
and it follows that
\begin{align}
  \ell_{q+1}^{-19}\lambda_{q+1}^{-\frac{3\g}{4}}+  \ell_{q}^{-19}\lambda_{q}^{-1/16}+\ell_{q}^{-23} \lambda_{q}^{-\frac98(\frac2{p_0}-1)} \lesssim \ell_{q}.\label{para:2}
\end{align}
}
\end{definition}

\begin{definition}
To deal with the stochastic terms, we define 
{$\k_{q}:= \ell_{q}^{1/4}$.} 
For $q \geq 1$,  we  define $\bar{Z}_{q+1}$ by
\begin{equation*}
\bar{Z}_{q+1}:=  \sup_{t \geq 0}\E[\|{z}_{q}+ \z- \bar{z}_{q} \|_{C_tL^2}^{2r_0}] +\sup_{t\geq0}\E[\|\bar{z}_{q+1}- \bar{z}_{q} \|_{C_tL^2}^{2r_0}] + \sup_{t \geq \k_{q+1}}\ell_{q+1}^{2\g r_0}\ \E[ \|\bar{z}_{q} \|_{C^{\g}_{t}H^{\g}}^{2r_0}].
\end{equation*}

 To prove the non-uniqueness of constructed solutions, we define $\{\theta_q\}_{q \in \mathbb{N}}$ to be some sequence where $\theta_q$ is a constant for all $q \in \mathbb{N}$ and 
\[
\sum_{q \in \mathbb{N}}\theta_q^{1/2}< \infty.
\] 
\end{definition}

Let us introduce the following assumption.
\begin{assumption}\label{a:iterative:i}
 {We assume}  the following results holds.
\begin{enumerate}
\item 
There exists some universal constant $c_R$ such that
\begin{equation}\label{eq:R:i}
\$\mathring{R}_q\$_{L^1,r_0,2\k_{q}}\leq  c_R \eps_{q-1}^{},\ \ \ \$\mathring{R}_q\$_{L^1,r_0}\leq  {c}_R \Lambda_{q-1},
\end{equation}
\item
There exists some universal constant {$c_v$ such that
\begin{equation}\label{eq:dv:1:i}
\$v_{q}-v_{q-1}\$_{L^2,2r_0,3\k_{q}} \leq   c_v ( \theta_{q-1}^{1/2}+\eps_{q-2}^{1/2})\ ,\ \${v_{q}-v_{q-1}}\$_{L^2,2r_0,[0,3\k_{q}]} \leq {c}_v \Lambda_{q-2}^{1/2},
\end{equation}
}
and
\begin{equation}\label{eq:dv:2:i}
\$v_q-v_{q-1}\$_{L^{p_0},2r_0}+ {\$v_{q}-v_{q-1}\$_{L^{6/(3+2\gamma)},2r_0}}\leq c_v \eps_{q},\ \ {\$v_q\$_{L^{p_0},2r_0} \leq c_v \sum_{i=1}^q\eps_{i}}.
\end{equation}
{Moreover, it holds that
\begin{align*}
    \$v_{q}\$_{L^2,2r_0}\leq c_v \Lambda_{q-1}^{1/2}.
\end{align*}
}
\item
For any $t \geq 1$, there exists some universal constant, still denoted by $c_v$, such that
\begin{equation}\label{eq:dv:3:i}
|\E[\|v_{q}(t)\|_{L^2}^2]-\E[\|v_{q-1}(t)\|_{L^2}^2]- 3\cdot (2\pi)^3 \theta_{q-1}\big| < c_v \eps_{q-2}.
\end{equation}
\item
For any $m\geq 1$, there holds
\begin{equation*}
\Lambda_{q} \leq M_{q}(\a_0,m) \leq \eps_{q+1}^{-1}.
\end{equation*}
\end{enumerate}
\end{assumption}

\begin{remark}
As mentioned in Section \ref{sec:iteration}, a constant is universal if it is independent of the iteration step $q$. We remark here that all the implicit constants in '$\lesssim$' are universal and are independent of the constants $c_R$ and $c_v$. The constants $c_R$ and $c_v$ are chosen to be large enough according to the initial condition and the proof of Lemma \ref{l:iterative:R:p0:i}, Lemma \ref{l:iterative:v:i} and Lemma \ref{l:iterative:v:i:0}. Moreover $c_v$ depends on $c_R$ and $c_R$ is some universal constant.
\end{remark}

Our iteration starts from $v_0= v_1\equiv0$. By a similar argument  as in Proposition \ref{l:main},  we require the constant $\Xi$ to be large enough and have the following result.}

\begin{proposition}\label{p:iterative}
Let $v_0= v_1\equiv0$ on $\mathbb{R}$, {and ${z}_1$ is obtained by solving \eqref{eq:sc:i}. Then  there exists $\mathring{R}_1$ such that} Assumption \ref{a:iterative:i} holds for $q=1$.
\end{proposition}

The following is our main iteration:
\begin{proposition}\label{prop:2}
 Assume that  Assumption \ref{a:iterative:i} holds up to step $q$, then it also holds at step $q+1$. Moreover, there exists some universal constant $c_Z$ such that
\begin{equation*}
 \bar{Z}_{q+1} \leq c_Z \ell_{q}^{\g r_0}\Lambda_{q-1}^{r_0}.
\end{equation*}
\end{proposition}

\subsection{Proof of Proposition \ref{prop:2}}
{In this section, we prove the main iterative step Proposition \ref{prop:2} using an argument analogous to that of Proposition \ref{l:main}. The key additional consideration here is the treatment of the initial condition, while the remaining constructions and estimates follow similarly to the previous case and we kindly omitted for brevity.} 

\subsubsection{Estimate on the noise}
In this section, we aim to estimate the noise term $\bar{Z}_{q+1}$.  We start with the following estimate on the initial part $\z$ which is defined in \eqref{def:dotz}.
\begin{lemma}\label{l:z0:i}
For any {$t\geq0$} and $\g\in (0,1)$, there exists some constant $C_{c,\g}$ such that
\begin{equation*}
\|\z\|_{C^{\g}_{t}H^{2\g}}  \leq C_{c,\g}\ t^{-2\g} \|u_0\|_{L^2}.
\end{equation*}
\end{lemma}
\begin{proof}
By properties of analytic semigroups, for $t>s$, there exists some constant {$C_{\g}$ such that} 
\begin{equation*}
\begin{aligned}
&\|e^{-ct}S(t)u_0 - e^{-cs}S(s)u_0\|_{H^{2\g}}\\
 \leq& \|e^{-ct}S(t)u_0 - e^{-ct}S(s)u_0\|_{H^{2\g}}+ \|e^{-ct}S(s)u_0 - e^{-cs}S(s)u_0\|_{H^{2\g}}\\
\leq& C_{\g} \int_0^{t-s} \|(-A)^{1-\g}S(r)\|_{\mathcal{L}(L^2;L^2)} {\dif r}  \|(-A)^{2\g}S(s)u_0\|_{L^{2}} + c(t-s)\|S(s)u_0\|_{H^{2\gamma}}.
\end{aligned}
\end{equation*}
By \eqref{eq:asg}, we obtain that
\[
\|S(t)u_0 - S(s)u_0\|_{H^{2\g}} \leq C_{c,\g} (t-s)^{\g} s^{-2\g}\|u_0\|_{L^2},\ \ \ t>s,
\]
which implies the result.
\end{proof}

\begin{lemma}\label{l:barZ:q:i}
Assume Assumption \ref{a:iterative:i} holds up to step $q$. There exists some universal constant $c_Z$  such that 
\begin{equation}\label{eq:Z:i:t}
 \bar{Z}_{q+1} \leq  c_Z \ell_{q}^{\g r_0} \Lambda_{q-1}^{r_0}.
  \end{equation}
\end{lemma}
\begin{proof}
Notice that there holds $\bar{z}_{q+1}- \bar{z}_{q}  = \Pi_{\zeta_{q+1}}z_{q+1} - \Pi_{\zeta_{q}}z_{q}$ and ${z}_{q}+ \z- \bar{z}_{q} = z_{q}-  \Pi_{\zeta_{q}}z_{q} $. By processing exactly the same calculation as the proof of Lemma \ref{l:barZ:q}, we obtain that 
\begin{equation}\label{eq:Z:i:1}
\sup_{t \geq0}\E[\|{z}_{q}+ \z- \bar{z}_{q} \|_{C_tL^2}^{2r_0}] +\E[\|\bar{z}_{q+1}- \bar{z}_{q} \|_{C_tL^2}^{2r_0}]  + \ell_{q+1}^{2\g  r_0}\ \E[\| {z}_{q} \|_{C^{\g}_{t}H^{\g}}^{2r_0}]  \lesssim \ \ell_{q}^{\g r_0} \Lambda_{q-1}^{r_0}.
\end{equation}
For the second part of $\bar{Z}_{q+1}$, there holds 

\begin{equation}\label{eq:Z:i:2}
\ell_{q+1}^{\g}\ \E[ \|\bar{z}_{q} \|_{C^{\g}_{t}H^{\g}}^{2r_0}]^{1/(2r_0)} \leq \ell_{q+1}^{\g}\ (\E[ \|{z}_{q} \|_{C^{\g}_{t}H^{\g}}^{2r_0}]^{1/(2r_0)} + \|\z\|_{C_t^{\g}H^{2\g}} ).
\end{equation}

By Lemma \ref{l:z0:i}, we have 

\begin{equation}\label{eq:Z:i:3}
\sup_{s \geq \k_{q+1}}\ell_{q+1}^{\g} \|\z\|_{C_s^{\g}H^{2\g}}  \lesssim \ell_{q+1}^{\g/2} \ \|u_0\|_{L^2}.
\end{equation}

Combining \eqref{eq:Z:i:1}, \eqref{eq:Z:i:2} and \eqref{eq:Z:i:3}, we obtain \eqref{eq:Z:i:t}.
\end{proof}

\subsubsection{Estimate on $v_{q+1}$}
As processed in Section \ref{sec:mollification}, we introduce the mollification parameter 
\begin{equation*}
\ell:=\ell_{q+1}\in (0,1/2),
\end{equation*}
 and mollify $v_q$ on $\R$ {as \eqref{def:moll}}.  Since \eqref{induction ps:i} holds at step $q$ on $\R$, it follows that $(v_{q,\ell_{}},\mathring{R}_{q,\ell_{}})$ satisfies
\begin{equation*}
\aligned
\partial_tv_{q,\ell_{}} - c\bar{z}_{q,\ell_{}}-\Delta v_{q,\ell_{}}+\div((v_{q,\ell_{}}+\bar{z}_{q,\ell_{}})\otimes (v_{q,\ell_{}}+\bar{z}_{q,\ell_{}}))+\nabla p_{q,\ell_{}}&=\div (\mathring{R}_{q,\ell_{}}+\mathring{R}^{(q+1)}_{\textrm{mol}}),
\\\div v_{q,\ell_{}}&=0,\\
{v_{q,\ell}(0)}&=0,
\endaligned
\end{equation*}
 where
\begin{equation*}
\mathring{R}_{\textrm{mol}}^{(q+1)}=(v_{q,\ell}+\bar{z}_{q,\ell})\mathring{\otimes}(v_{q,\ell}+\bar{z}_{q,\ell})-((v_{q}+\bar{z}_{q})\mathring{\otimes}(v_{q}+\bar{z}_{q}))_{\ell}.
\end{equation*}

Let $\zeta_q$, $r_{\perp,q}$, $r_{\|,q}$, $\mu_q$ and $\lambda_q$ be as in Definition \ref{def:para:3}. 
 In this section, {the amplitude function is defined similarly as that defined in Section \ref{sec:perturbation} except that $\theta_q$ is the constant defined above, i.e.
 \begin{equation*}
\rho_{q+1}:=2\sqrt{\ell_{q+1}^2+|\mathring{R}_{q,\ell}|^2}+\theta_{q},
\end{equation*} 
and  
\begin{equation*}
a_{(\xi)}:=\rho_{q+1}^{1/2}\gamma_\xi\left(\Id-\frac{\mathring{R}_{q,\ell}}{\rho_{q+1}}\right),
\end{equation*}
where $\gamma_{\xi}$ is defined in Lemma \ref{geometric}. Following a calculation analogous to \cite[Appendix B]{HZZ22}, the estimates in \eqref{estimate aN} and \eqref{estimate aN0} remain valid. 
 }

 Then we define 
 ${w}^{(p)}_{q+1}$, ${w}^{(c)}_{q+1}$ and ${w}^{(t)}_{q+1}$ in the same manner as in Section \ref{sec:perturbation}.

To keep the initial condition, a cut-off function $\chi_q(t)$ is introduced as follows:
\begin{definition}
 We define the smooth time cut-off function $\chi_q$ be such that
\begin{equation*}
\begin{aligned}
\chi_q(t) :=
\begin{cases}
0\ \ \ \ &\text{when}\ t \leq \k_{q}\\
1\ \ \ \ &\text{when}\ t \geq 2\k_{q}\\
\end{cases}\ \ \ \text{and}\ \ \ \|\chi_q(t)\|_{C^n_t} \lesssim_n \k_{q}^{-n}. 
\end{aligned}
\end{equation*}
\end{definition}

Let $v_{q+1}:= v_{q,\ell} + \tilde{w}_{q+1}$, where the perturbation $\tilde{w}_{q+1}:=\tilde{w}^{(p)}_{q+1} + \tilde{w}^{(c)}_{q+1}+\tilde{w}^{(t)}_{q+1}$ is defined by
\begin{equation*}
\begin{aligned}
\tilde{w}^{(p)}_{q+1}:= \chi_{q+1} w^{(p)}_{q+1},\ \ \tilde{w}^{(c)}_{q+1}:= \chi_{q+1} w^{(c)}_{q+1},\ \ \tilde{w}^{(t)}_{q+1}:=&\ \chi_{q+1}^2 w^{(t)}_{q+1}.
\end{aligned}
\end{equation*}

\begin{lemma}\label{l:iterative:v:i}
Assume Assumption \ref{a:iterative:i} holds up to step $q$. There exists some universal constant $c_v$ such that Assumption \ref{a:iterative:i}(3) holds at step $q+1$, i.e.{
\begin{equation}\label{l:iterative:v:i:1}
 \$v_{q+1}-v_{q}\$_{L^2,2r_0,3\k_{q}} \leq   c_v ( \theta_{q}^{1/2}+\eps_{q-1}^{1/2}),\ \ \$v_{q+1}-v_q\$_{L^2,2r_0,[0,3\k_{q}]} \leq {c}_v  \Lambda_{q-1}^{1/2},
\end{equation}}
Moreover, there exists some constant, still denoted by $c_v$, such that
\begin{equation}\label{l:iterative:v:i:2}
 \$v_{q+1}-v_q\$_{L^{p_0},2r_0} +{\$v_{q+1}-v_q\$_{L^{6/(3+2\gamma)},2r_0}} \leq c_v\eps_{q+1},\ \  {\$v_{q+1}\$_{L^{p_0},2r_0} \leq c_v \sum_{i=1}^{q+1}\eps_{i}}.
\end{equation}
and 
\begin{equation}\label{l:iterative:v:i:3}
\begin{aligned}
    \$v_{q+1}\$_{L^2,2r_0}\leq c_v(1+\$\mathring{R}_{q}\$_{L^1,r_0}^{1/2}+ \$v_{q}\$_{L^2,2r_0})\leq c_v \Lambda_q^{1/2}.
\end{aligned}
\end{equation}
\end{lemma}
\begin{proof}
By \eqref{eq:principle:Lp}, \eqref{est:correction:Lp} and \eqref{est:temporal:Lp}, {a similar calculation as \eqref{bd:wq+1l2r} implies that} there holds
\begin{equation*}
\begin{aligned}
\$w_{q+1}\$_{L^2,2r_0,3\k_{q}}\lesssim&\ \ell_{q+1} (1+ \$\mathring{R}_{q}\$_{L^1,6r_0,3\k_{q}-\ell_{q+1}}^{3}) + (\ell_{q+1} + \$\mathring{R}_{q}\$_{L^1,r_0,3\k_{q}-\ell_{q+1}} +{\theta}_{q})^{1/2}.
\end{aligned}
\end{equation*}
Since $\k_{q}> \ell_{q+1}$, there holds
\begin{equation}\label{bd:wq+1:l2:1}
\begin{aligned}
\$w_{q+1}\$_{L^2,2r_0,3\k_{q}}\lesssim&\ \ell_{q+1}^{\g} (1+ \$\mathring{R}_{q}\$_{L^1,6r_0,2\k_{q}}^{3}) + (\ell_{q+1} + \$\mathring{R}_{q}\$_{L^1,r_0,2\k_{q}} +{\theta}_{q})^{1/2}\\
\lesssim&\ \ell_{q+1}^{\g} \Lambda_{q} + (\ell_{q+1} + c_R\eps_{q-1}+{\theta}_{q})^{1/2}{\lesssim  \eps_{q-1}^{1/2}+ \theta_{q}^{1/2}}.
\end{aligned}
\end{equation}
Moreover, we have
\begin{equation}\label{eq:vql:i}
\begin{aligned}
 \$ v_{q,\ell}- v_q \$_{L^2,2r_0} \lesssim \ell_{q+1} \$v_q\$_{\C^1,2r_0} \lesssim \ell_{q+1} \Lambda_{q}{\lesssim  \eps_{q-1}^{1/2}}.
\end{aligned}
\end{equation}
Hence, there exists some constant $c_v$ such that
\begin{equation*}
\$v_{q+1}-v_{q}\$_{L^2,2r_0,3\k_{q}} \leq \$w_{q+1}\$_{L^2,2r_0,3\k_{q}} + \$ v_{q,\ell}- v_q \$_{L^2,2r_0} \leq c_v (\eps_{q-1}^{1/2}+ \theta_{q}^{1/2}).
\end{equation*}
Therefore, we have proved the first statement in \eqref{l:iterative:v:i:1}. 

To prove the second statement in \eqref{l:iterative:v:i:1}, we can prove similarly that
\begin{equation}\label{bd:wq+1:l2:0}
\begin{aligned}
\$w_{q+1}\$_{L^2,2r_0,[0,3\k_{q}]}\lesssim&\ \ell_{q+1}^{\g} (1+ \$\mathring{R}_{q}\$_{L^1,6r_0,[0,3\k_{q}]}^{3}) + (\ell_{q+1} + \$\mathring{R}_{q}\$_{L^1,r_0,[0,3\k_q]} +{\theta}_{q})^{1/2}\\
 \lesssim&\ {\ell_{q+1}^{\g} \Lambda_{q} + (\ell_{q+1} + \$\mathring{R}_{q}\$_{L^1,r_0} +{\theta}_{q})^{1/2}\lesssim \Lambda_{q-1}^{1/2}.}
\end{aligned}
\end{equation}
{
Hence, there exists some universal constant, still denoted by ${c}_v$,  such that
\[
\begin{aligned}
\$v_{q+1}-v_{q}\$_{L^2,2r_0,[0,3\k_q]} \leq \$w_{q+1}\$_{L^2,2r_0,[0,3\k_{q}]} + \$ v_{q,\ell}- v_q \$_{L^2,2r_0} \leq {c}_v  \Lambda_{q-1}^{1/2}.
\end{aligned}
\]
}
Moreover, by \eqref{bd:wq+1:l2:1}, \eqref{eq:vql:i} and \eqref{bd:wq+1:l2:0}, there holds 
\[
 \$v_{q+1}\$_{L^2,2r_0} \lesssim \ell_{q+1}^{\g} \Lambda_{q} + (\ell_{q+1} + \$\mathring{R}_{q}\$_{L^1,r_0} +{\theta}_{q})^{1/2}+ \$v_{q}\$_{L^2,2r_0}.
\]
By \eqref{eq:R:i}, we obtain \eqref{l:iterative:v:i:3}.

Then, the proof of \eqref{l:iterative:v:i:2} follows directly from  
\eqref{para:2} and \eqref{eq:w:L1} 
 \begin{align*}
     \$w_{q+1}\$_{L^{p_0},2r_0}+ \$w_{q+1}\$_{L^{6/(3+2\g)},2r_0}\lesssim \ell_{q+1}^{-19}( \lambda_{q+1}^{-\frac98(\frac2{p_0}-1)} +\lambda_{q+1}^{-\frac{3\g}{4}})\Lambda_{q}^{1/2}\lesssim \ell_{q+1}\Lambda_{q}^{1/2},
 \end{align*}
 which implies that
\begin{equation*}
\begin{aligned}
\$v_{q+1}&-v_{q}\$_{L^{p_0},2r_0} +\$v_{q+1}-v_{q}\$_{L^{6/(3+2\g)},2r_0} \notag\\ 
&\leq \$w_{q+1}\$_{L^{p_0},2r_0} + \$w_{q+1}\$_{L^{6/(3+2\g)},2r_0}+\ell_{q+1} \$v_q \$_{C_{t,x}^1,2r_0} \lesssim \ell_{q+1}^{\g} \Lambda_{q}^{1/2}\lesssim \eps_{q+1}.
\end{aligned}
\end{equation*}
\end{proof}

\subsubsection{Estimate on $\mathring{R}_{q+1}$}
Similar to the calculation in Section \ref{sec:w}, we obtain that
\begin{equation*}
\aligned
\mathring{R}_{q+1}=\mathring{{R}}_{\textrm{cut}}^{(q+1)}+ \mathring{{R}}_{\textrm{noise}}^{(q+1)}&+ \mathring{{R}}_{\textrm{comm,1}}^{(q+1)}+ \mathring{{R}}_{\textrm{comm,2}}^{(q+1)} + \mathring{{R}}_{\textrm{mol}}^{(q+1)} \\
&+ \mathring{{R}}^{(q+1)}_{\textrm{osc,x}} + \mathring{{R}}^{(q+1)}_{\textrm{osc,t}} + \mathring{{R}}^{(q+1)}_{\textrm{cor}} + \mathring{{R}}^{(q+1)}_{\textrm{lin,w}} +\mathring{{R}}_{\textrm{lin,z}}^{(q+1)},
\endaligned
\end{equation*}
where
\begin{equation*}
\aligned
\mathring{{R}}^{(q+1)}_{\textrm{osc,x}}:=&\chi^2_{q+1}\sum_{\xi \in \Lambda} \cB \(\nabla a^2_{(\xi)} , \mP_{\neq 0}\(W_{(\xi)} {\otimes} W_{(\xi)}\) \),\\
\mathring{{R}}^{(q+1)}_{\textrm{osc,t}}:=& -\frac{1}{\mu}\chi^2_{q+1}\sum_{\xi \in \Lambda}  \mathcal{R}\(\mP_{\neq0} \partial_t \(a_{(\xi)}^2\)\phi_{(\xi)}^2 \psi_{(\xi)}^2 \xi\),\\
\mathring{{R}}^{(q+1)}_{\textrm{cut}}:=& \frac{1}{\mu}   (\chi_{q+1}^2)'\sum_{\xi \in \Lambda} \mP\mP_{\neq0} a_{(\xi)}^2 \phi_{(\xi)}^2 \psi_{(\xi)}^2 \xi  + \chi_{q+1}' ({w}_{q+1}^{(p)} + {w}_{q+1}^{(c)})+ (1-\chi^2_{q+1})\mathring{R}_{q,\ell} ,
\endaligned
\end{equation*}

\begin{equation*}
\aligned
\mathring{{R}}^{(q+1)}_{\textrm{cor}} :=& (\tilde{w}_{q+1}^{(t)} + \tilde{w}_{q+1}^{(c)}) \mathring{\otimes}_s \tilde{w}_{q+1}^{(p)} +  (\tilde{w}_{q+1}^{(t)} + \tilde{w}_{q+1}^{(c)}) \mathring{\otimes}  (\tilde{w}_{q+1}^{(t)} + \tilde{w}_{q+1}^{(c)}),\\
\mathring{{R}}^{(q+1)}_{\textrm{lin,w}}:=& \mathcal{R}\(\chi_{q+1}\partial_t  ({w}_{q+1}^{(p)} + {w}_{q+1}^{(c)}) -   \Delta \tilde{w}_{q+1}\),
\endaligned
\end{equation*}
and 
\begin{equation*}
\aligned
\mathring{R}_{\textrm{lin,z}}^{(q+1)}:=& \mathcal{R}(-c\ \mP_{\neq 0}\mP(\bar{z}_{q+1}- \bar{z}_{q,\ell})), \\
\mathring{R}_{\textrm{noise}}^{(q+1)}: =& (\bar{z}_{q+1}-\bar{z}_{q,l}) \mathring{\otimes} (\bar{z}_{q+1}- \bar{z}_{q,l}),\\
\mathring{R}_{\textrm{comm,1}}^{(q+1)}:=& (\bar{z}_{q+1} - \bar{z}_{q,\ell}) \mathring{\otimes}_s (v_{q+1}+ \bar{z}_{q,\ell}),\\
\mathring{R}_{\textrm{comm,2}}^{(q+1)}:=&  \tilde{w}_{q+1} \mathring{\otimes}_s (v_{q,\ell} + \bar{z}_{q,\ell}),\\
\mathring{R}_{\textrm{mol}}^{(q+1)}:=&(v_{q,\ell}+\bar{z}_{q,\ell})\mathring{\otimes}(v_{q,\ell}+\bar{z}_{q,\ell})-((v_q+\bar{z}_q)\mathring{\otimes}(v_q+\bar{z}_q))_{\ell}.
\endaligned
\end{equation*}

Similar to Section \ref{sec:w}, the following estimates for the Reynolds stress hold, {since the cut off functions $\chi_{q+1}$ is bounded}. {To prove that $\|\mathring{R}_{q+1}(t)\|_{L^1}$ is small for $t >\k_{q+1}$, we require the following estimate which could be obtained by a similar argument as \eqref{Rq+1:l1:r}:}
\begin{equation}\label{Rq+1:l1:i}
\begin{aligned}
&\|\mathring{R}_{q+1} -\mathring{R}^{(q+1)}_{\textrm{cut}}\|_{C_tL^1} \lesssim\ \ell_{q+1}^{-32} \lambda_{q+1}^{-1/8}(1+ \|\mathring{R}_q\|_{C_{[t-\ell,t+1]}L^1}^{6})(1+\| v_{q}\|_{C_{[t-\ell,t+1]}L^{2}} + \|\bar{z}_{q}\|_{C_{[t-\ell,t+1]} L^{2}})\\
& + \(\|\bar{z}_{q+1}- \bar{z}_{q} \|_{C_{[t-\ell,t+1]}L^2}+  \ell_{q+1}^{\g}\|\bar{z}_{q} \|_{C^{\g}_{[t-\ell,t+1]}H^{\g}} \)+ \(\|\bar{z}_{q+1}-\bar{z}_{q} \|^2_{C_{t}L^2} +  \ell_{q+1}^{2\g}\|\bar{z}_{q} \|_{C^{\g}_{[t-\ell,t+1]}H^{\g}}^2\)\\
&+ \( \|\bar{z}_{q+1}- \bar{z}_{q} \|_{C_{[t-\ell,t+1]}L^2} +  \ell_{q+1}^{\g}\|\bar{z}_{q} \|_{C^{\g}_{[t-\ell,t+1]}H^{\g}}\) ( \|v_{q+1}\|_{C_{t}L^2} +  \|\bar{z}_q\|_{C_{[t-\ell,t+1]}L^2})\\
&+\ell_{q+1}^{\g} (\|v_q\|_{C_{[t-\ell,t+1],x}^1} + \|\bar{z}_{q} \|_{C^{\g}_{[t-\ell,t+1]}H^{\g}} )(\|v_q\|_{C_{[t-\ell,t+1]}L^2}+\|\bar{z}_q\|_{C_{[t-\ell,t+1]}L^2}).
\end{aligned}
\end{equation}
{To bound any $m$-th moment for $\|\mathring{R}_{q+1}(t)\|_{L^1}$, we can process as in  \eqref{Rq+1:l1:m} and obtain the following estimate:
\begin{equation}\label{Rq+1:l1:m:i}
\begin{aligned}
\|\mathring{R}_{q+1} -\mathring{R}^{(q+1)}_{\textrm{cut}}\|_{C_tL^1} \lesssim& 1+ \|\mathring{R}_q\|_{C_{[t-\ell,t+1]}L^1}^{{12}} + \|\bar{z}_{q+1}\|_{C_{[t-\ell,t+1]}L^2}^{2} + \|\bar{z}_{q}\|_{C_{[t-\ell,t+1]}L^2}^{2}\\
& + \|v_{q+1}\|_{C_{[t-\ell,t+1]}L^2}^{2}+ {\|v_{q}\|_{C_{[t-\ell,t+1]}L^2}^{2}}.
\end{aligned}
\end{equation}
}

{It remains to bound the term $\|\mathring{R}^{(q+1)}_{\textrm{cut}}\|_{C_{{[t,t+1]}}L^1}$.} By \eqref{bounds} and \eqref{estimate aN0}, there holds
\begin{equation*}
\big\|  \frac{1}{\mu}   (\chi_{q+1}^2)'\sum_{\xi \in \Lambda} \mP\mP_{\neq0} a_{(\xi)}^2 \phi_{(\xi)}^2 \psi_{(\xi)}^2 \xi  \big\|_{C_{{[t,t+1]}}L^1}\lesssim \chi_{[\k_{q+1},2\k_{q+1}]}({[t,t+1]})  \ell_{q+1}^{\g}   (\|\mathring{R}_q\|_{C_{{[t-\ell,t+1]}}L^1}+1).
\end{equation*}
{By \eqref{est:principle:Lp} and \eqref{est:correction:Lp} with $p=1+\epsilon$ for some $\epsilon>0$ small enough, there holds
\begin{equation}
\begin{aligned}
\| \chi_{q+1}'({w}_{q+1}^{(p)} + {w}_{q+1}^{(c)})\|_{C_{{[t,t+1]}}L^1}&\lesssim \chi_{[\k_{q+1},2\k_{q+1}]}({[t,t+1]}) \ell_{q+1}^{-20}\lambda_{q+1}^{-\frac18+3\epsilon}  (\|\mathring{R}_q\|_{C_{{[t-\ell,t+1]}}L^1}^2+1)\notag\\
&\lesssim \chi_{[\k_{q+1},2\k_{q+1}]}({[t,t+1]})  \ell_{q+1}^{\g}  (\|\mathring{R}_q\|_{C_{{[t-\ell,t+1]}}L^1}^2+1),
\end{aligned}
\end{equation}
where we chose $\epsilon>0$ such that $\lambda_{q+1}^{3\epsilon}<\ell_{q+1}^\alpha$.}

Therefore, there exists some constant $M_{\textrm{cut}}$ such that
\begin{equation}\label{eq:cut:i}
\begin{aligned}
\|\mathring{{R}}^{(q+1)}_{\textrm{cut}} \|_{C_{{[t,t+1]}}L^1} \leq &M_{\textrm{cut}} \ \chi_{[\k_{q+1},2\k_{q+1}]}({[t,t+1]})  \ell_{q+1}^{\g}  (\|\mathring{R}_q\|_{C_{{[t-\ell,t+1]}}L^1}^2+1) \\
&+ \chi_{[0,2\k_{q+1}]}({[t,t+1]}) \|\mathring{{R}}_{q}\|_{C_{{[t-\ell,t+1]}}L^1}.
\end{aligned}
\end{equation}

\begin{lemma}\label{l:M:bound:i}
Assume Assumption \ref{a:iterative:i} holds up to step $q$. Then Assumption \ref{a:iterative:i}(5) holds at step $q+1$, i.e.  for any $m \geq 1$, there holds
\begin{equation*}
\Lambda_{q+1} \leq M_{q+1}(\a_0,m) \leq  \eps_{q+2}^{-1}.
\end{equation*}
\end{lemma}

{\begin{proof}
We can prove exactly the same as Lemma \ref{l:M:bound} where we require the constants $\a_0$ and $M_0$ to be large enough.
\end{proof}}

\begin{lemma}\label{l:iterative:R:p0:i}
Assume Assumption \ref{a:iterative:i} holds up to step $q$. There exists some universal constants $c_R$ such that
\begin{equation}\label{eq:a:R:p:q+1:i:1}
\$\mathring{R}_{q+1}\$_{L^1,r_0,2\k_{q+1}}  \leq c_R  \eps_q^{}.
\end{equation}
In particular, we have
\begin{equation}\label{eq:a:R:p:q+1:i:2}
\$\mathring{R}_{q+1}\$_{L^1,r_0}  \leq c_{R} \Lambda_{q} .
\end{equation}
\end{lemma}
\begin{proof}
By \eqref{Rq+1:l1:i}, Lemma \ref{l:barZ:q:i}, there exists some universal constant $c_R$ such that 
By \eqref{Rq+1:l1:r} we obtain 
\begin{equation*}
\begin{aligned}
&\$\mathring{R}_{q+1} - \mathring{R}^{(q+1)}_{\textrm{cut}}\$_{L^1,r_0} \lesssim \ell_{q+1}^{-32} \lambda_{q+1}^{-1/8}(1+ \$\mathring{R}_q\$_{L^1,12r_0}^{12}+\$ v_{q}\$_{L^{2},2r_0}^2 + \$\bar{z}_{q}\$_{L^{2},2r_0}^2) + \bar{Z}_{q+1}^{1/r_0}\\
& +\bar{Z}_{q+1}^{1/(2r_0)} (1+ \$v_{q+1}\$_{L^2,2r_0}+ \$\bar{z}_{q}\$_{L^2,2r_0})+\ell_{q+1}^{\g}(1+\$v_{q}\$_{C_{t,x}^1,2r_0}^2+ \$\bar{z}_{q}\$_{C^{\g}_tH^{\g},2r_0}^2).
\end{aligned}
\end{equation*}
By Lemma \ref{l:truncate} and \eqref{zq:moment:1} we obtain that
\begin{equation}\label{bd:barzq:l2}
\E[\|\bar{z}_{q}\|_{C_{[t,t+1]} L^{2}}^{2r_0}] \leq \E[\|{z}_{q} \|_{C^{\g}_{[t,t+1]}H^{\g}}^{2r_0}] \lesssim \$v_{q-1}\$_{L^{p_0},2r_0}^{2r_0}{+ \|u_0\|_{L^2}^{2r_0}\lesssim1}.
\end{equation}

Moreover, by the first inequality in \eqref{l:iterative:v:i:3} and \eqref{eq:R:i}, there holds
\[
\$v_{q+1}\$_{L^2,2r_0} \leq c_v(1+   \Lambda_{q-1}^{1/2} +   \$v_{q}\$_{L^2,2r_0}) .
\]
Applying \eqref{l:iterative:v:i:3} again with $q$ instead of $q+1$, we obtain that
\[
\$v_{q+1}\$_{L^2,2r_0} \lesssim 1+   \Lambda_{q-1}^{1/2}    \lesssim  \Lambda_{q-1}^{1/2}.
\]
 Then by  Lemma \ref{l:barZ:q:i} we obtain that
\begin{equation*}
\begin{aligned}
&\$\mathring{R}_{q+1}- \mathring{R}^{(q+1)}_{\textrm{cut}} \$_{L^1,r_0,{2\kappa_{q+1}}} \lesssim\ \ell_{q+1}^{\g} 
\Lambda_q+ \bar{Z}_{q+1}^{1/r_0} +\bar{Z}_{q+1}^{1/(2r_0)}\Lambda_{q-1}^{1/2}\lesssim\ \ell_{q+1}^{\g}\Lambda_q+\ell_{q}^{\g/2}\Lambda_{q-1}\lesssim \eps_{q}.
\end{aligned}
\end{equation*}

Since $\mathring{R}^{(q+1)}_{\textrm{cut}}=0$ when $t \geq 2\k_{q+1}$, we obtain \eqref{eq:a:R:p:q+1:i:1} {for some $c_R>0$}. 

{
By the same calculation as \eqref{bd:barzq:l2} there holds
\[
\sup_{t \geq 0}\E[\|\bar{z}_{q+1}\|_{C_{[t-\ell,t+1]}L^2}^{2r_0}]\lesssim \$v_q\$_{L^{p_0},2r_0} ^{2r_0}+ \|u_0\|_{L^2}^{2r_0}\lesssim 1.
\]
Combining this with \eqref{bd:wq+1:l2:1}, \eqref{bd:wq+1:l2:0}, \eqref{Rq+1:l1:m:i} and \eqref{eq:cut:i}, we obtain \eqref{eq:a:R:p:q+1:i:2}:
\begin{equation*}
\begin{aligned}
\$\mathring{R}_{q+1}\$_{L^1,r_0} \lesssim& 1+ \$\mathring{R}_q\$_{L^1,12r_0}^{12} + \$w_{q+1}\$_{L^2,2r_0}^{2}+  \$v_{q}\$_{C_{t,x}^1,2r_0}^{2}\leq c_R\Lambda_q.
\end{aligned}
\end{equation*}
}
 \end{proof}

\subsubsection{Estimate on the energy}

\begin{lemma}\label{l:iterative:v:i:0}
Assume Assumption \ref{a:iterative:i} holds up to step $q$. There exists some universal constant, still denoted by $c_v$ such that Assumption \ref{a:iterative:i}(4) holds at step $q+1$, i.e.  for any $t \geq 1$,
\begin{equation*}
|\E[\|v_{q+1}(t)\|_{L^2}^2]-\E[\|v_{q}(t)\|_{L^2}^2]- 3\cdot(2\pi)^3\theta_q\big| < c_v \eps_{q-1}.
\end{equation*}
\end{lemma}
\begin{proof}
We have the following decomposition: {for $t\geq1$}
\begin{equation*}
\begin{aligned}
&\big| \E[\|v_{q+1}{(t)}\|^2_{L^2}-\|v_{q}{(t)}\|_{L^2}^2] -3\cdot(2\pi)^3 \theta_q \big|\\
\leq& \big| \E[\|w_{q+1}^{(p)}{(t)}\|_{L^2}^2 -3\cdot(2\pi)^3 \theta_q \big| + 2\E[\|w_{q+1}^{(p)}{(t)}(w_{q+1}^{(t)}+w_{q+1}^{(c)}){(t)}\|_{L^1}] + \E[\|(w_{q+1}^{(t)}+w_{q+1}^{(c)}){(t)}\|_{L^2}^2]\\
&+ 2\E[\|w_{q+1}{(t)}v_{q,\ell}{(t)}\|_{L^1}] +\left| \E[\|v_{q,\ell}{(t)}\|_{L^2}^2-\|v_q{(t)}\|_{L^2}^2]\right|.
\end{aligned}
\end{equation*}
By \eqref{est:principle:Lp}, \eqref{est:correction:Lp} and \eqref{est:temporal:Lp}, we obtain that (similar to the proof of \eqref{eq:cor})
\begin{equation}\label{eq:v:3}
\E[\|w_{q+1}^{(p)}{(t)}(w_{q+1}^{(t)}+w_{q+1}^{(c)}){(t)}\|_{L^1}] + \E[\|(w_{q+1}^{(t)}+w_{q+1}^{(c)}){(t)}\|_{L^2}^2]{\lesssim \ell_{q+1}^{-32} \lambda_{q+1}^{-1/8} \Lambda_q} \lesssim \ell_{q+1}^{\g} \Lambda_q.
\end{equation}
{By \eqref{eq:w:L1} we obtain that 
\begin{align*}
    \E[\|w_{q+1}{(t)}v_{q,\ell}{(t)}\|_{L^1}]\lesssim \ell_{q+1}^{-19} \lambda_{q+1}^{-1/8}\Lambda_q\lesssim \ell_{q+1}^{\g} \Lambda_q.
\end{align*}}
By basic mollification estimate, we have
\begin{equation}\label{eq:v:1}
\begin{aligned}
 \left|\E[\|v_{q,\ell}{(t)}\|_{L^2}^2-\|v_q{(t)}\|_{L^2}^2]\right| 
\lesssim& \ell_{q+1}\E[\|v_q\|_{C_{t,x}^1}^2] \lesssim \ell_{q+1}^{\g} \Lambda_q.
\end{aligned}
\end{equation}
Since $\mathring{R}_{q,\ell}$ is trace free, by \eqref{can}, \eqref{eq:dE:aux1}  we have
\begin{equation}\label{eq:v:2}
\begin{aligned}
|\E\|w_{q+1}^{(p)}(t)\|_{L^2}^2- 3\cdot(2\pi)^3\theta_{q}|& \leq 6 \cdot (2\pi)^3\ell_{q+1}+ 6\E\|\mathring{R}_{q,\ell}(t)\|_{L^1} + \mathbf{E}\sum_{\xi\in \Lambda}\Big|\int a_{(\xi)}^2\mathbb{P}_{\neq0}|W_{(\xi)}|^2 \dif x\Big|\\
&{\lesssim \ell_{q+1}+\eps_{q-1}}.
\end{aligned}
\end{equation}
Combining \eqref{eq:v:3}, \eqref{eq:v:1} and \eqref{eq:v:2}, we can to derive the result {by the definition of $\eps_q$}.
\end{proof}

\subsection{Proof of Theorem \ref{thm:main:3}} 
\begin{proof}
    We first notice that by \eqref{eq:dv:1:i}
\begin{equation*}
\begin{aligned}
\E[\int_0^T \|v_{q+1}(t)-v_q(t)\|_{L^2}^2 \dif t] &\leq  \int_{3\k_{q+1}}^T \E [\|v_{q+1}(t)-v_q(t)\|_{L^2}^2] \dif t +  \int_{0}^{3\k_{q+1}}\E[ \|v_{q+1}(t)-v_q(t)\|_{L^2}^2 ]\dif t \\
&\leq T (\eps_{q-1} + \theta_{q})+\ell_{q+1}^{1/4} \Lambda_{q-1}^{}\lesssim \eps_{q-1}
+\theta_{q}.
\end{aligned}
\end{equation*}

Hence, we obtain that
\begin{equation*}
\sum_{q\in \mathbb{N}} \E[ \int_0^T \|v_{q+1}(t)-v_q(t)\|_{L^2}^2 \dif t] {\lesssim \sum_{q \in \mathbb{N}}\eps_{q-1} + \theta_{q}}<\infty,
\end{equation*}
which implies that  $v \in L^2_{\textrm{loc}}([0,\infty);L^2) $ almost surely. By \eqref{eq:dv:2:i},  it's easy to prove that $v \in C([0,\infty);L^{p_0}\cap L^{6/(3+2\gamma)})$  with a bound independent of $\{\theta_q\}$. Similar to the proof of \eqref{eq:dv:converge}, we obtain that
\begin{equation}\label{eq:thm:pf:1}
v \in C([0,\infty);L^{p_0}\cap L^{6/(3+2\gamma)})\cap C((0,\infty);W^{\vartheta,p_0}) \cap  L^2_{\textrm{loc}}([0,\infty);L^2) ,
\end{equation}
 almost surely for some $\vartheta>0$ sufficiently small. 

 Similar to the proof of Theorem \ref{thm:main:1} , ${\bar{z}}_{q}$ converges to some ${z}$ in $C([0,\infty);H^{2\g})$ {with a bound independent of $\{\theta_q\}$}. Combining this with \eqref{eq:thm:pf:1}, we have
\[
{u:=v+z} \in C([0,\infty);L^{p_0}) \cap  L^2_{\textrm{loc}}([0,\infty);L^2) \cap  C((0,\infty);W^{\vartheta, p_0}) .
\]

 By \eqref{eq:R:i}, there holds 
\begin{equation}
\begin{aligned}
\E[\int_0^T \|R_q(t)\|_{L^1}^{r_0} \dif t] &\lesssim \int_{2\k_{q}}^T \E[\|R_q(t)\|_{L^1}^{r_0}] \dif t + \int_0^{2\k_{q}} \E[\|R_q(t)\|_{L^1}^{r_0}] \dif t\notag\\ 
&\lesssim \eps_{q-1}^{r_0} + \k_{q}\Lambda_{q-1}^{r_0}\lesssim \eps_{q-1} + \ell_{q}^{1/4}\Lambda_{q-1}^{r_0} \lesssim\eps_{q-1}.
\end{aligned}
\end{equation}

Therefore, we can process similarly as in the proof of Theorem \ref{thm:main:1} and obtain that  $u$ is a probabilistically strong and analytically weak solution solution to the Navier-Stokes equation \eqref{1} on $[0, \infty)$ with initial condition $u_0$.

Now, we will prove the non-uniqueness of the solutions. Since $v_1=0$, we have
\begin{equation*}
\begin{aligned}
\left|\E[ \|v(t)+{z}(t)\|_{L^2}^2]- 3\cdot (2\pi)^3\theta_1 \right| =&  \big|  \E[\|{z}_{}(t)\|_{L^2}^2 + 2\|v_{}(t) {z}(t)\|_{L^1}+\|v_{}(t)\|_{L^2}^2] - 3\cdot (2\pi)^3\theta_1 \big|\\
\leq&   \E[\|{z}_{}(t)\|_{L^2}^2] +  2 \E[\|v_{}(t){z}(t)\|_{L^1}] +3\cdot (2\pi)^3\sum_{q \geq 2}\theta_q \\
&+ \sum_{q\geq 1} \big|  \E[\|v_{q+1}(t)\|_{L^2}^2-\|v_{q}(t)\|_{L^2}^2] - 3\cdot (2\pi)^3 \theta_q\big|.
\end{aligned}
\end{equation*}

To bound the first two terms, we apply H\"older inequality and obtain that
\[
\begin{aligned}
 \E[\|{z}_{}(t)\|_{L^2}^2] +  2 \E[\|v_{}(t){z}(t)\|_{L^1}] \leq  \E[\|{z}_{}(t)\|_{L^2}^2] +  2 \E[\|v_{}(t)\|_{L^{6/(3+2\g)}}^2]^{1/2}\E[\|{z}(t)\|_{H^{\g}}^2]^{1/2}<\infty,
 \end{aligned}
\]
where the constant is independent of $\{\theta_q\}$.

Combining all these with \eqref{eq:dv:3:i}, there exists some constant $C_{vz}$ such that
\begin{equation*}
\left|\E[ \|v(t)+z(t)\|_2^2]- 3\cdot (2\pi)^3\theta_1 \right|^2 \leq C_{vz} + 3\cdot (2\pi)^3\sum_{q \geq 2}\theta_q := K_0,
\end{equation*}
where the constant $K_0$ is independent of $\theta_1$. Therefore, we can pick two different values $K_1$ and $K_2$ for $\theta_1$ such that $|K_1-K_2|>2K_0$ so that we obtain two different values for $\E[ \|v(t)+z(t)\|_{L^2}^2]$, $t\geq 1$, which concludes the proof.

\end{proof}

The proof of Corollary \ref{coro:noninlaw} follows by the same argument as  \cite[Corollary 1.2]{HZZ23a}.

\noindent{\bf Acknowledgment.}  We are very grateful to Prof. Xiangchan Zhu for proposing this  problem and for her invaluable suggestions, which greatly enhanced our work.

\appendix
\renewcommand{\appendixname}{Appendix~\Alph{section}}
\renewcommand{\theequation}{A.\arabic{equation}}

\section{Intermittent jets}
\label{s:B}

In this  part we recall the construction of  intermittent jets from \cite[Section 7.4]{BV19}.
We point out that the construction is entirely deterministic, that is, none of the functions below depends on $\omega$.
Let us begin with  the following geometric lemma which can be found in  \cite[Lemma 6.6]{BV19}.

\bl\label{geometric}
Denote by $\overline{B_{1/2}}(\mathrm{Id})$ the closed ball of radius $1/2$ around the identity matrix $\mathrm{Id}$, in the space of $3\times 3$ symmetric matrices. There
exists $\Lambda\subset \mathbb{S}^2\cap \mathbb{Q}^3$ such that for each $\xi\in \Lambda$ there exists a  $C^\infty$-function $\gamma_\xi:\overline{B_{1/2}}(\mathrm{Id})\rightarrow\mathbb{R}$ such that
\begin{equation*}
R=\sum_{\xi\in\Lambda}\gamma_\xi^2(R)(\xi\otimes \xi)
\end{equation*}
for every symmetric matrix satisfying $|R-\mathrm{Id}|\leq 1/2$.
For $C_\Lambda=8|\Lambda|(1+8\pi^3)^{1/2}$, where $|\Lambda|$ is the cardinality of the set $\Lambda$, we define
the constant
\begin{equation*}
M=C_\Lambda\sup_{\xi\in \Lambda}(\|\gamma_\xi\|_{C^0}+\sum_{|j|\leq N}\|D^j\gamma_\xi\|_{C^0}).
\end{equation*}
For each $\xi\in \Lambda$ let us define $A_\xi\in \mathbb{S}^2\cap \mathbb{Q}^3$ to be an orthogonal vector to $\xi$. Then for each $\xi\in\Lambda$ we have that $\{\xi, A_\xi, \xi\times A_\xi\}\subset \mathbb{S}^2\cap \mathbb{Q}^3$ form an orthonormal basis for $\mathbb{R}^3$.
We label by $n_*$ the smallest natural such that
\begin{equation*}\{n_*\xi, n_*A_\xi, n_*\xi\times A_\xi\}\subset \mathbb{Z}^3\end{equation*}
for every $\xi\in \Lambda$.
\el

Let $\Phi:\mathbb{R}^2\rightarrow\mathbb{R}$ be a smooth function with support in a ball of radius $1$. We normalize $\Phi$ such that
$\phi=-\Delta \Phi$ obeys
\begin{equation}\label{eq:phi}
\frac{1}{4\pi^2}\int_{\mathbb{R}^2}\phi^2(x_1,x_2)\dif x_1\dif x_2=1.
\end{equation}
By definition we know $\int_{\mathbb{R}^2}\phi  \dif x=0$. Define $\psi:\mathbb{R}\rightarrow\mathbb{R}$ to be a smooth, mean zero function with support in the ball of radius $1$ satisfying
\begin{equation}\label{eq:psi}
\frac{1}{2\pi}\int_{\mathbb{R}}\psi^2(x_3)\dif x_3=1.
\end{equation}
For parameters $r_\perp, r_\|>0$ such that
\begin{equation*}
r_\perp\ll r_\|\ll1,\end{equation*}
we define the rescaled cut-off functions
\begin{equation*}\phi_{r_\perp}(x_1,x_2)=\frac{1}{r_\perp}\phi\left(\frac{x_1}{r_\perp},\frac{x_2}{r_\perp}\right),\quad
\Phi_{r_\perp}(x_1,x_2)=\frac{1}{r_\perp}\Phi\left(\frac{x_1}{r_\perp},\frac{x_2}{r_\perp}\right),\quad \psi_{r_\|}(x_3)=\frac{1}{r_\|^{1/2}}\psi\left(\frac{x_3}{r_\|}\right).\end{equation*}
We periodize $\phi_{r_\perp}, \Phi_{r_\perp}$ and $\psi_{r_\|}$ so that they are viewed as periodic functions on $\mathbb{T}^2, \mathbb{T}^2$ and $\mathbb{T}$ respectively.

Consider a large real number $\lambda$ such that 
\begin{equation*}
\lambda r_\perp\in\mathbb{N},
\end{equation*}
 and a large time oscillation parameter $\mu>0$. For every $\xi\in \Lambda$ we introduce
\begin{equation*}\aligned
\psi_{(\xi)}(t,x)&:=\psi_{\xi,r_\perp,r_\|,\lambda,\mu}(t,x):=\psi_{r_{\|}}(n_*r_\perp\lambda(x\cdot \xi+\mu t))
\\ \Phi_{(\xi)}(x)&:=\Phi_{\xi,r_\perp,\lambda}(x):=\Phi_{r_{\perp}}(n_*r_\perp\lambda(x-\alpha_\xi)\cdot A_\xi, n_*r_\perp\lambda(x-\alpha_\xi)\cdot(\xi\times A_\xi))\\
\phi_{(\xi)}(x)&:=\phi_{\xi,r_\perp,\lambda}(x):=\phi_{r_{\perp}}(n_*r_\perp\lambda(x-\alpha_\xi)\cdot A_\xi, n_*r_\perp\lambda(x-\alpha_\xi)\cdot(\xi\times A_\xi)),
\endaligned\end{equation*}
where $\alpha_\xi\in\mathbb{R}^3$ are shifts to ensure that $\{\Phi_{(\xi)}\}_{\xi\in\Lambda}$ have mutually disjoint support.
 We refer to \cite[Section 7.4]{BV19} for the specific choice of such  $\alpha_\xi$.

The intermittent jets $W_{(\xi)}:\mathbb{R}\times\mathbb{T}^3 \rightarrow\mathbb{R}^3$ are defined as in \cite[Section 7.4]{BV19}.
\begin{equation}\label{intermittent}
W_{(\xi)}(t,x):=W_{\xi,r_\perp,r_\|,\lambda,\mu}(t,x):=\xi\psi_{(\xi)}(t,x)\phi_{(\xi)}(x).\end{equation}
By the choice of $\alpha_\xi$ we have that
\begin{equation*}
W_{(\xi)}\otimes W_{(\xi')}\equiv0, \textrm{for} \xi\neq \xi'\in\Lambda,
\end{equation*}
and by the normalizations \eqref{eq:phi} and \eqref{eq:psi} we obtain
$$
\frac1{(2\pi)^3}\int_{\mathbb{T}^3}W_{(\xi)}(t,x)\otimes W_{(\xi)}(t,x)\dif x=\xi\otimes\xi.
$$
These facts combined with Lemma \ref{geometric} imply that
\begin{equation*}
\frac1{(2\pi)^3}\sum_{\xi\in\Lambda}\gamma_\xi^2(R)\int_{\mathbb{T}^3}W_{(\xi)}(t,x)\otimes W_{(\xi)}(t,x)\dif x=R,
\end{equation*}
for every symmetric matrix $R$ satisfying $|R-\textrm{Id}|\leq 1/2$. Since $W_{(\xi)}$ are not divergence free, we  introduce the corrector term
\begin{equation}\label{corrector}
W_{(\xi)}^{(c)}:=\frac{1}{n_*^2\lambda^2}\nabla \psi_{(\xi)}\times \textrm{curl}(\Phi_{(\xi)}\xi)
=\textrm{curl\,curl\,} V_{(\xi)}-W_{(\xi)}.
\end{equation}
with
\begin{equation*}
V_{(\xi)}(t,x):=\frac{1}{n_*^2\lambda^2}\xi\psi_{(\xi)}(t,x)\Phi_{(\xi)}(x).
\end{equation*}
Thus we have
\begin{equation*}
\div\left(W_{(\xi)}+W_{(\xi)}^{(c)}\right)\equiv0.
\end{equation*}

Finally, we  recall the key   bounds from \cite[Section 7.4]{BV19}. For $N, M\geq0$ and $p\in [1,\infty]$ the following holds
\begin{equation}\label{bounds}\aligned
&\|\nabla^N\partial_t^M\psi_{(\xi)}\|_{C_{t}L^p}\lesssim r^{1/p-1/2}_\|\left(\frac{r_\perp\lambda}{r_\|}\right)^N
\left(\frac{r_\perp\lambda \mu}{r_\|}\right)^M,\\
&\|\nabla^N\phi_{(\xi)}\|_{L^p}+\|\nabla^N\Phi_{(\xi)}\|_{L^p}\lesssim r^{2/p-1}_\perp\lambda^N,\\
&\|\nabla^N\partial_t^MW_{(\xi)}\|_{C_{t}L^p}+\frac{r_\|}{r_\perp}\|\nabla^N\partial_t^MW_{(\xi)}^{(c)}\|_{C_{t}L^p}+\lambda^2\|\nabla^N\partial_t^MV_{(\xi)}
\|_{C_{t}L^p}\\&\lesssim r^{2/p-1}_\perp r^{1/p-1/2}_\|\lambda^N\left(\frac{r_\perp\lambda\mu}{r_\|}\right)^M,
\endaligned\end{equation}
where the implicit constants may depend on $p, N$ and $M$, but are independent of $\lambda, r_\perp, r_\|, \mu$.

\renewcommand{\theequation}{B.\arabic{equation}}

\def\cprime{$'$} \def\ocirc#1{\ifmmode\setbox0=\hbox{$#1$}\dimen0=\ht0
  \advance\dimen0 by1pt\rlap{\hbox to\wd0{\hss\raise\dimen0
  \hbox{\hskip.2em$\scriptscriptstyle\circ$}\hss}}#1\else {\accent"17 #1}\fi}

\end{document}